\documentclass{amsart}

\usepackage{amssymb}
\usepackage{amsmath}
\usepackage{amsthm}
\usepackage{amsfonts}

\usepackage{a4wide}

\newtheorem{lemma}{Lemma}
\newtheorem{proposition}{Proposition}
\newtheorem{theorem}{Theorem}

\theoremstyle{definition}
\newtheorem{definition}{Definition}
\newtheorem{example}{Example}

\newtheorem{remark}{Remark}

\numberwithin{equation}{section}

\begin{document}

%\linenumbers

\title[Coclass trees of CF-groups and related BCF-groups]
{New perspectives of the power-commutator-structure: \\ Coclass trees of CF-groups and related BCF-groups
}

\author{Daniel C. Mayer}
\address{Naglergasse 53\\8010 Graz\\Austria}
\email{algebraic.number.theory@algebra.at}
\urladdr{http://www.algebra.at}

\thanks{Research supported by the Austrian Science Fund (FWF): projects J0497-PHY, P26008-N25, and by EUREA}



\date{Friday, 31 December 2021}

%--------------------------------------------------------------------------------

\begin{abstract}
Let \(e\ge 2\) be an integer.
Among the finite \(3\)-groups \(G\) with
bicyclic commutator quotient \(G/G^\prime\simeq C_{3^e}\times C_3\),
having one non-elementary component with logarithmic exponent \(e\),
there exists a unique pair of coclass trees
with distinguished rank distribution \(\varrho\sim (2,2,3;3)\).
One tree \(\mathcal{T}^{e}(M^{(e)}_1)\)
consists of CF-groups with coclass \(e\),
and the other tree \(\mathcal{T}^{e+1}(\mathbb{M}^{(e+1)}_1)\)
consists of BCF-groups with coclass \(e+1\).
It is proved that,
due to a chain of periodic bifurcations,
the vertices of all pairs \((\mathcal{T}^{e},\mathcal{T}^{e+1})\) with \(e\ge 3\)
can be constructed as \(p\)-descendants of the
single root \(M^{(3)}_1\) of order \(3^6\) by means of the
\(p\)-group generation algorithm by Newman and O'Brien.
\end{abstract}

\maketitle

%\newpage
%--------------------------------------------------------------------------------

\section{Introduction}
\label{s:Intro}

\noindent
We consider finite \(3\)-groups \(G\)
whose commutator quotient is \textit{bicyclic with one non-elementary component}, that is,
\(G/G^\prime\simeq C_{3^e}\times C_3\) with logarithmic exponent \(e\ge 2\).
By the Burnside basis theorem,
\(G=\langle x,y\rangle\) is two-generated,
and we stick to the convention that
\(w=x^{3^{e-1}}\), \(w^{3}\in G^\prime\) and \(y^{3}\in G^\prime\)
for the generators.

%--------------------------------------------------------------------------------

%\noindent
For such groups, we have introduced several invariants
\cite{Ma2021a}
in terms of their maximal normal subgroups \(H_1,\ldots,H_3;H_4\) of index \((G:H_i)=3\),
where the puncture at the fourth component is motivated by the distinction of
the bicyclic quotient \(H_4/G^\prime\simeq C_{3^{e-1}}\times C_3\),
as opposed to the cyclic quotients \(H_i/G^\prime\simeq C_{3^e}\) for \(1\le i\le 3\).
We have the logarithmic \textit{abelian quotient invariants} (AQI),
\[\alpha(G)=\bigl(H_1/H_1^\prime,\ldots,H_3/H_3^\prime;H_4/H_4^\prime\bigr),\]
their \textit{rank distribution}
\(\varrho(G)=\bigl(\mathrm{rank}_3(H_1/H_1^\prime),\ldots,\mathrm{rank}_3(H_3/H_3^\prime);\mathrm{rank}_3(H_4/H_4^\prime)\bigr)\),
and, on the other hand,
the \textit{punctured transfer kernel type} (pTKT),
\[\varkappa(G)=\bigl(\ker(T_1),\ldots,\ker(T_3);\ker(T_4)\bigr),\]
where \(T_i:\,G/G^\prime\to H_i/H_i^\prime\) denotes the Artin transfer homomorphism from \(G\) to \(H_i\).
AQI and pTKT are combined in the \textit{Artin pattern} \(\mathrm{AP}(G)=(\alpha(G),\varkappa(G))\) of \(G\).
Since there are only five possibilities for the kernels, the pTKT is abbreviated in the form
\(\varkappa(G)=\bigl(\varkappa_1,\ldots,\varkappa_3;\varkappa_4\bigr)\), where
\[\varkappa_i=
\begin{cases}
0\text{ if } \ker(T_i)=\langle w,y\rangle/G^\prime \text{ (complete 3-elementary subgroup)}, \\
j\text{ if } \ker(T_i)=\langle w^{j-1}\cdot y\rangle/G^\prime,\ 1\le j\le 3, \\
4\text{ if } \ker(T_i)=\langle w\rangle/G^\prime \text{ (distinguished third power generator)}.
\end{cases}\]

%--------------------------------------------------------------------------------

\noindent
Our special attention is devoted to \textit{CF-groups} \(G\)
for which the factors \(\gamma_i(G)/\gamma_{i+1}(G)\), \(i\ge 2\),
of the descending central series \((\gamma_i(G))_{i\ge 1}\) are cyclic of order \(3\)
(\textit{cyclic factor groups}).
Since \(\gamma_2(G)=\langle s_2,\gamma_3(G)\rangle\),
the second factor is always cyclic, but since
\(\gamma_3(G)=\langle s_3,t_3,\gamma_4(G)\rangle\),
the third factor is usually bicyclic
(\textit{BCF} --- \textit{bicyclic or cyclic factor groups}),
and there must exist some relation between \(s_3\) and \(t_3\) in a CF-group,
for instance, either one of the two commutators is trivial or \(s_3=t_3\).
(Here, \(s_2=\lbrack y,x\rbrack\), \(s_3=\lbrack s_2,x\rbrack\), and \(t_3=\lbrack s_2,y\rbrack\)
denote essential commutators.)

%--------------------------------------------------------------------------------

%\noindent
Even more specifically,
our focus will lie on \textit{coclass trees}
\cite[p. 89]{Ma2018}
whose vertices \(G\) share the common rank distribution
\(\varrho(G)\sim (2,2,3;3)\), that is,
trees of CF-groups with mainline of type \(\mathrm{a}.1\), \(\varkappa(G)=(0,0,0;0)\),
and trees of BCF-groups with mainline of type \(\mathrm{d}.10\), \(\varkappa(G)\sim (1,1,0;2)\).

%--------------------------------------------------------------------------------

%\noindent
After a summary of foundations in \S\
\ref{s:Basics},
we begin with simple laws for all mainlines of CF-coclass trees in \S\
\ref{s:LawsCF}.
Then we extend the investigations to chains of periodic bifurcations in \S\
\ref{s:LawsBCF},
where the complete system of all CF- and BCF-coclass trees
with \(e\ge 3\)
is shown to arise from a single root.

%\newpage
%--------------------------------------------------------------------------------

\section{Basic definitions and conventions}
\label{s:Basics}

\noindent
The lower exponent \(p\)-central series of a finite \(3\)-group \(G\) will always be denoted by
\((P_i(G))_{i\ge 0}\).

\begin{definition}
\label{dfn:Parent}
Let \(D\) be non-trivial finite \(3\)-group
with nilpotency class \(c=\mathrm{cl}(D)\ge 1\)
and lower exponent \(p\)-class \(c_p=\mathrm{cl}_p(D)\ge 1\), i.e.,
\(\gamma_{c}(D)>\gamma_{c+1}(D)=1\) and \(P_{c_p-1}(D)>P_{c_p}(D)=1\).
Then the quotient \(A=\pi(D)=D/\gamma_{c}(D)\) is called
the \textit{parent} of \(D\)
and the quotient \(A_p=\pi_p(D)=D/P_{c_p-1}(D)\) is called
the \(p\)-\textit{parent} of \(D\).
Conversely, \(D\) is called
an \textit{immediate descendant} of \(A\) and
an \textit{immediate \(p\)-descendant} of \(A_p\).
By the \textit{root path}, respectively \(p\)-\textit{root path}, of \(D\)
we understand the sequence \((\pi^i(D))_{i\ge 0}\), respectively  \((\pi_p^i(D))_{i\ge 0}\),
of its iterated parents, respectively \(p\)-parents.
\end{definition}

%--------------------------------------------------------------------------------

\begin{definition}
\label{dfn:EndoAndExo}
The \textit{propagation} from a \(p\)-parent \(A\) to a \(p\)-descendant \(D\) is called
\textbf{endo-genetic} if the commutator quotient remains unchanged,
that is, \(D/D^\prime\simeq A/A^\prime\).
Otherwise the propagation is called \textbf{exo-genetic}.
\end{definition}

\noindent
The propagation from non-trivial parent \(A\) to non-abelian descendant \(D\) is always endo-genetic,
because \(A=D/\gamma_c(D)\) with \(c\ge 2\), and thus
\(D/\gamma_2(D)\simeq (D/\gamma_c(D))/(\gamma_2(D)/\gamma_c(D))\simeq A/\gamma_2(A)\).

%--------------------------------------------------------------------------------

\begin{definition}
\label{dfn:DescendantTree}
The \textit{descendant tree} \(\mathcal{T}(R)\),
respectively \(p\)-\textit{descendant tree} \(\mathcal{T}_p(R)\),
with a finite \(3\)-group \(R\) as its root consists of
the following vertices and directed edges:
the vertices are all isomorphism classes of finite \(3\)-groups \(D\)
whose root path, respectively \(p\)-root path, contains \(R\),
and the directed edges are all pairs \((D,A)\), also denoted by \(D\to A\),
of immediate descendants \(D\) and parents \(A=\pi(D)\), respectively \(p\)-parents \(A=\pi_p(D)\),
among the vertices of the tree.
A descendant tree
whose vertices are subject to certain restrictive conditions
is called a \textit{pruned tree}.
\end{definition}

%--------------------------------------------------------------------------------

\begin{definition}
\label{dfn:CoclassTree}
A pruned tree
which contains a unique infinite \textit{main line} and
all of whose vertices share a common coclass \(r\)
is called a \textit{coclass tree}.
If the root is \(R\), the tree is denoted by \(\mathcal{T}^r(R)\).
\end{definition}

\noindent
The step size of all edges in a coclass tree is necessarily \(s=1\).
\textit{Depth-pruned} branches of a coclass tree become periodic,
beginning with a minimal \textit{periodic root} on the main line
\cite[Thm. 3.1]{Ma2018}.

%--------------------------------------------------------------------------------

\begin{definition}
\label{dfn:TreeType}
By a \textit{tree of type} \(\mathrm{X}\)
we understand a coclass tree
whose mainline consists 
of vertices with (punctured) transfer kernel type \(\mathrm{X}\)
\cite[Tbl. 1--2, pp. 3--4]{Ma2021a}.
(For instance \(\mathrm{X}=\mathrm{a}.1\) or \(\mathrm{d}.10\).)
\end{definition}

%--------------------------------------------------------------------------------

\noindent
We introduce an ostensive terminology
in order to illuminate three distinct situations
with crucial differences in the construction
by means of the \(p\)-group generation algorithm
\cite{Nm1977,Ob1990,HEO2005,GNO2006}.

\begin{definition}
\label{dfn:ShockWave}
A vertex \(V\) on a coclass tree \(\mathcal{T}^r\),
with \(V/V^\prime\simeq C_{3^e}\times C_3\) and \(r\in\lbrace e,e+1\rbrace\),
lies
\begin{itemize}
\item
\textbf{behind the shock wave}, if \(\mathrm{cl}(V)<r\),
\item
\textbf{on the shock wave}, if \(\mathrm{cl}(V)=r\),
\item
\textbf{ahead of the shock wave}, if \(\mathrm{cl}(V)>r\).
\end{itemize}
\end{definition}

\noindent
The behavior ahead of the shock wave will turn out to be
\textit{regular} with endo-genetic propagation,
dominated by the \textit{commutator structure}.
In contrast, we shall see that
the behavior behind the shock wave
is \textit{irregular} with exo-genetic propagation,
due to a dominance of the \textit{power structure}.
A \textit{singular} behavior can be observed on the shock wave,
where the propagation is mixed, partially endo-genetic and partially exo-genetic,
and \textit{periodic bifurcations} arise,
because both, the commutator structure and the power structure,
exert a combined impact.

%--------------------------------------------------------------------------------

%\noindent
In order to identify the isomorphism class of a finite \(3\)-group \(G\),
several ways are possible.

\noindent
Either the group is characterized by its \textit{absolute identifier}
\(\mathrm{SmallGroup}(o,i)\), or briefly \(\langle o,i\rangle\),
in the SmallGroups database
\cite{BEO2005},
where \(o=\#G\) denotes the order of \(G\),
bounded by \(o\le 3^8\),
and \(i\) is a positive integer.
The short form in angle brackets is returned by the Magma statement
\texttt{IdentifyGroup()}
\cite{BCP1997,BCFS2021,MAGMA2021},
provided that \(o\le 3^6\).
When the order \(o=3^e\) or the logarithmic order \(e\)
is given along a scale on the left hand side of a figure
illustrating a descendant tree of finite \(3\)-groups,
then we omit the order \(o\) in the absolute identifier
\(\langle o,i\rangle\) and simply write \(\langle i\rangle\).

Or \(G\) is constructed by means of the Magma statement
\texttt{Descendants(}\(P\)\texttt{:StepSizes:=\lbrack}\(s\)\texttt{\rbrack)}
as an immediate step size-\(s\) \(p\)-descendant of a \(p\)-parent group \(P\)
and characterized by a \textit{relative identifier}
\(G=P-\#s;j\) with \(1\le s\le n(P)\) and \(1\le j\le N_s(P)\),
where \(n(P)\) denotes the nuclear rank of the \(p\)-parent \(P\)
and \(N_s(P)\) is the number of immediate step size-\(s\) \(p\)-descendants of \(P\).

Finally, there is always the possibility to give a
\textit{power commutator (pc-) presentation} for \(G\).

%\newpage
%--------------------------------------------------------------------------------

\section{Laws for coclass trees of CF-groups}
\label{s:LawsCF}

\noindent
We separate our main statements into three parts:
uniqueness, invariants, and construction.

\begin{proposition}
\label{prp:UniquenessCF}
For each logarithmic exponent \(e\ge 2\),
there exists a unique coclass tree \(\mathcal{T}^e(M^{(e)}_1)\ni V\)
with fixed coclass \(\mathrm{cc}(V)=e\),
fixed commutator quotient \(V/V^\prime\simeq C_{3^e}\times C_3\),
and fixed rank distribution \(\varrho(V)\sim (2,2,3;3)\).
Its mainline \((M^{(e)}_i)_{i\ge 1}\)
is of type \(\mathrm{a}.1\), \(\varkappa(M^{(e)}_i)=(000;0)\).
The tree consists entirely of metabelian \(\mathrm{CF}\)-groups.
The branches are of depth \(3\). (See Figure
\ref{fig:Tree21a1AscioneA}.)
\end{proposition}

\begin{proof}
For each commutator quotient \(C_{3^e}\times C_3\) with log exponent \(e\ge 2\),
there exists a \textit{finite} number \(N\) of coclass trees \(\mathcal{T}^r(R^{(r)}_j)\)
with roots \(R^{(r)}_j\), \(1\le j\le N\),
and two minimal possible values \(e\le r\le e+1\) for the coclass.
Descendant vertices \(V\) of each root share invariants with the root,
e.g. the rank distribution \(\varrho(V)\).
The roots with \(r=e+1\) are non-CF groups (called \textit{BCF-groups} in
\cite{Ne1989},
i.e. groups with \textit{bicyclic or cyclic factors} of the lower central series),
and the others with \(r=e\) are CF-groups.
There are only two trees with rank distribution \(\varrho(V)\sim (2,2,3;3)\),
a BCF-tree of type \(\mathrm{d}.10\) and a CF-tree of type \(\mathrm{a}.1\).
The latter is the \textit{unique} tree with root \(R^{(e)}_j=M^{(e)}_1\),
recursively determined by
\(M^{(3)}_1=\langle 729,7\rangle\) and \(M^{(2)}_1=\langle 243,17\rangle\),
according to Theorems
\ref{thm:ConstructionMain}
and
\ref{thm:ConstructionOffside}.
Its branches are periodic of length \(2\) without pre-period,
and all of its vertices are metabelian CF-groups,
since the vertices of the first two branches are metabelian CF-groups.
\end{proof}

%\newpage
%--------------------------------------------------------------------------------

\subsection{Vertices on the mainline (with depth \(0\))}
\label{ss:MainlineCF}

\begin{proposition}
\label{prp:InvariantsCF}
For \(e\ge 3\),
invariants of vertices on the mainline \((M^{(e)}_i)_{i\ge 1}\)
of the coclass tree \(\mathcal{T}^e(M^{(e)}_1)\)
are given as follows:
\begin{equation}
\label{eqn:InvariantsCF}
\begin{aligned}
\text{logarithmic order } \mathrm{lo}(M^{(e)}_i)=e+i+2,
\text{ nilpotency class } \mathrm{cl}(M^{(e)}_i) &= i+2, \text{ for } i\ge 1, \\
p\text{-class } \mathrm{cl}_p(M^{(e)}_i)=
\begin{cases}
i+2 & \text{ if } i>e-2, \\
e   & \text{ if } i\le e-2,
\end{cases}
p\text{-coclass } \mathrm{cc}_p(M^{(e)}_i) &=
\begin{cases}
e   & \text{ if } i>e-2, \\
i+2 & \text{ if } i\le e-2.
\end{cases}
\end{aligned}
\end{equation}
\end{proposition}

\begin{proof}
Proposition
\ref{prp:InvariantsCF}
remains true when the mainline vertex \(M^{(e)}_i\)
is replaced by any proper descendant vertex
\(V^{(e)}_i\) with \(i\ge 2\).
All coclass trees under investigation
start at a root of class \(\mathrm{cl}(M^{(e)}_1)=3=1+2\),
for each \(e\ge 2\).
Thus, proper descendants possess nilpotency class \(\mathrm{cl}(V^{(e)}_i)=i+2\ge 4\).
By definition, all vertices \(V\)
of the coclass tree \(\mathcal{T}^e(M^{(e)}_1)\)
share the common coclass \(\mathrm{cc}(V)=e\).
Consequently, the logarithmic order is the sum
\(\mathrm{lo}(V^{(e)}_i)=\mathrm{cl}(V^{(e)}_i)+\mathrm{cc}(V^{(e)}_i)=i+2+e\).
Finally, the \textit{power structure} of all finite \(3\)-groups \(G\)
with commutator quotient \(G/G^\prime\simeq C_{3^e}\times C_3\) is responsible for
the constant \(p\)-class \(\mathrm{cl}_p(V^{(e)}_i)=e\),
independently of the class \(\mathrm{cl}(V^{(e)}_i)=i+2\le e\),
in the finite region on and behind the shock wave.
\end{proof}

%--------------------------------------------------------------------------------

\begin{theorem}
\label{thm:ConstructionMain}
Vertices on the mainline \((M^{(e)}_i)_{i\ge 1}\)
of the coclass tree \(\mathcal{T}^e(M^{(e)}_1)\)
can be constructed recursively,
according to three laws in dependence on the nilpotency class,
\begin{itemize}
\item
by \textbf{irregular exo}-genetic propagation
(behind the shock wave)
\begin{equation}
\label{eqn:Irregular}
M^{(e)}_i=M^{(e-1)}_i-\#1;1, \text{ for } e\ge 4,\ i\le e-3, \text{ i.e. }\mathrm{cl}(M^{(e)}_i)<e,
\end{equation}
\item
by \textbf{singular exo}-genetic propagation
(bifurcation on the shock wave)
\begin{equation}
\label{eqn:Singular}
M^{(e)}_{i}=M^{(e-1)}_{i-1}-\#2;1, \text{ for } e\ge 4,\ i=e-2, \text{ i.e. }\mathrm{cl}(M^{(e)}_i)=e,
\end{equation}
\item
by \textbf{regular endo}-genetic propagation
(ahead of the shock wave)
\begin{equation}
\label{eqn:Regular}
M^{(e)}_{i}=M^{(e)}_{i-1}-\#1;1, \text{ for } e\ge 3,\ i\ge e-1, \text{ i.e. }\mathrm{cl}(M^{(e)}_i)>e.
\end{equation}
\end{itemize}
\end{theorem}

%--------------------------------------------------------------------------------

\begin{figure}[ht]
\caption{Mainlines of CF-coclass trees and their various mechanisms of propagation}
\label{fig:MainlinesCF}

{\tiny

\setlength{\unitlength}{1cm}
\begin{picture}(14,10)(-4.5,-9)

% scale of nilpotency class
\put(-5,0.5){\makebox(0,0)[cb]{class}}

\put(-5,0){\line(0,-1){8}}
\multiput(-5.1,0)(0,-1){9}{\line(1,0){0.2}}

\put(-5.2,0){\makebox(0,0)[rc]{\(3\)}}
\put(-5.2,-1){\makebox(0,0)[rc]{\(4\)}}
\put(-5.2,-2){\makebox(0,0)[rc]{\(5\)}}
\put(-5.2,-3){\makebox(0,0)[rc]{\(6\)}}
\put(-5.2,-4){\makebox(0,0)[rc]{\(7\)}}
\put(-5.2,-5){\makebox(0,0)[rc]{\(8\)}}
\put(-5.2,-6){\makebox(0,0)[rc]{\(9\)}}
\put(-5.2,-7){\makebox(0,0)[rc]{\(10\)}}
\put(-5.2,-8){\makebox(0,0)[rc]{\(11\)}}

\put(-5,-8){\vector(0,-1){1}}

% header
\put(-3,0.5){\makebox(0,0)[cb]{\((9,3)\)}}
\put(-1,0.5){\makebox(0,0)[cb]{\((27,3)\)}}
\put(1,0.5){\makebox(0,0)[cb]{\((81,3)\)}}
\put(3,0.5){\makebox(0,0)[cb]{\((243,3)\)}}
\put(5,0.5){\makebox(0,0)[cb]{\((729,3)\)}}
\put(7,0.5){\makebox(0,0)[cb]{\((2187,3)\)}}
\put(9,0.5){\makebox(0,0)[cb]{\((6561,3)\)}}

% mainlines
\put(-3,0){\line(0,-1){3}}

\put(-1,0){\line(0,-1){3}}
\put(1,0){\line(0,-1){4}}
\put(3,0){\line(0,-1){5}}
\put(5,0){\line(0,-1){6}}
\put(7,0){\line(0,-1){7}}
\put(9,0){\line(0,-1){8}}

% irregular vertices (behind)
\multiput(1,0)(2,-1){5}{\circle{0.1}}
\multiput(3,0)(2,-1){4}{\circle{0.1}}
\multiput(5,0)(2,-1){3}{\circle{0.1}}
\multiput(7,0)(2,-1){2}{\circle{0.1}}
\multiput(9,0)(2,-1){1}{\circle{0.1}}

% singular vertices (shock wave)
\multiput(-1,0)(2,-1){6}{\circle*{0.1}}

% regular vertices (ahead of)
\multiput(-3,0)(0,-1){4}{\circle{0.1}}
\put(-3,-3){\vector(0,-1){1}}

\multiput(-1,-1)(2,-1){6}{\circle{0.1}}
\multiput(-1,-2)(2,-1){6}{\circle{0.1}}
\multiput(-1,-3)(2,-1){6}{\circle{0.1}}
\multiput(-1,-3)(2,-1){6}{\vector(0,-1){1}}

% periodic bifurcation (exo-genetic)
\multiput(0,-0.3)(2,-1){5}{\makebox(0,0)[ct]{bifurcation}}
\multiput(-1,0)(2,-1){5}{\vector(2,0){2}}
\multiput(-1,0)(2,-1){5}{\vector(2,-1){2}}
\put(9,-5){\line(2,-1){1}}
\put(10,-5.3){\makebox(0,0)[ct]{shock wave}}

% irregular propagation
\multiput(1.5,0)(2,0){4}{\vector(2,0){1}}
\multiput(3.5,-1)(2,0){3}{\vector(2,0){1}}
\multiput(5.5,-2)(2,0){2}{\vector(2,0){1}}
\multiput(7.5,-3)(2,0){1}{\vector(2,0){1}}

% symbols for vertices
\put(-2.9,0.1){\makebox(0,0)[lb]{\(M^{(2)}_1\)}}
\put(-2.9,-0.9){\makebox(0,0)[lb]{\(M^{(2)}_2\)}}
\put(-2.9,-1.9){\makebox(0,0)[lb]{\(M^{(2)}_3\)}}
\put(-2.9,-2.9){\makebox(0,0)[lb]{\(M^{(2)}_4\)}}

\put(-0.9,0.1){\makebox(0,0)[lb]{\(M^{(3)}_1\)}}
\put(-0.9,-0.9){\makebox(0,0)[lb]{\(M^{(3)}_2\)}}
\put(-0.9,-1.9){\makebox(0,0)[lb]{\(M^{(3)}_3\)}}
\put(-0.9,-2.9){\makebox(0,0)[lb]{\(M^{(3)}_4\)}}

\put(1.1,0.1){\makebox(0,0)[lb]{\(M^{(4)}_1\)}}
\put(1.1,-0.9){\makebox(0,0)[lb]{\(M^{(4)}_2\)}}
\put(1.1,-1.9){\makebox(0,0)[lb]{\(M^{(4)}_3\)}}
\put(1.1,-2.9){\makebox(0,0)[lb]{\(M^{(4)}_4\)}}
\put(1.1,-3.9){\makebox(0,0)[lb]{\(M^{(4)}_5\)}}

\put(3.1,0.1){\makebox(0,0)[lb]{\(M^{(5)}_1\)}}
\put(3.1,-0.9){\makebox(0,0)[lb]{\(M^{(5)}_2\)}}
\put(3.1,-1.9){\makebox(0,0)[lb]{\(M^{(5)}_3\)}}
\put(3.1,-2.9){\makebox(0,0)[lb]{\(M^{(5)}_4\)}}
\put(3.1,-3.9){\makebox(0,0)[lb]{\(M^{(5)}_5\)}}
\put(3.1,-4.9){\makebox(0,0)[lb]{\(M^{(5)}_6\)}}

\put(5.1,0.1){\makebox(0,0)[lb]{\(M^{(6)}_1\)}}
\put(5.1,-0.9){\makebox(0,0)[lb]{\(M^{(6)}_2\)}}
\put(5.1,-1.9){\makebox(0,0)[lb]{\(M^{(6)}_3\)}}
\put(5.1,-2.9){\makebox(0,0)[lb]{\(M^{(6)}_4\)}}
\put(5.1,-3.9){\makebox(0,0)[lb]{\(M^{(6)}_5\)}}
\put(5.1,-4.9){\makebox(0,0)[lb]{\(M^{(6)}_6\)}}
\put(5.1,-5.9){\makebox(0,0)[lb]{\(M^{(6)}_7\)}}

\put(7.1,0.1){\makebox(0,0)[lb]{\(M^{(7)}_1\)}}
\put(7.1,-0.9){\makebox(0,0)[lb]{\(M^{(7)}_2\)}}
\put(7.1,-1.9){\makebox(0,0)[lb]{\(M^{(7)}_3\)}}
\put(7.1,-2.9){\makebox(0,0)[lb]{\(M^{(7)}_4\)}}
\put(7.1,-3.9){\makebox(0,0)[lb]{\(M^{(7)}_5\)}}
\put(7.1,-4.9){\makebox(0,0)[lb]{\(M^{(7)}_6\)}}
\put(7.1,-5.9){\makebox(0,0)[lb]{\(M^{(7)}_7\)}}
\put(7.1,-6.9){\makebox(0,0)[lb]{\(M^{(7)}_8\)}}

\put(9.1,0.1){\makebox(0,0)[lb]{\(M^{(8)}_1\)}}
\put(9.1,-0.9){\makebox(0,0)[lb]{\(M^{(8)}_2\)}}
\put(9.1,-1.9){\makebox(0,0)[lb]{\(M^{(8)}_3\)}}
\put(9.1,-2.9){\makebox(0,0)[lb]{\(M^{(8)}_4\)}}
\put(9.1,-3.9){\makebox(0,0)[lb]{\(M^{(8)}_5\)}}
\put(9.1,-4.9){\makebox(0,0)[lb]{\(M^{(8)}_6\)}}
\put(9.1,-5.9){\makebox(0,0)[lb]{\(M^{(8)}_7\)}}
\put(9.1,-6.9){\makebox(0,0)[lb]{\(M^{(8)}_8\)}}
\put(9.1,-7.9){\makebox(0,0)[lb]{\(M^{(8)}_9\)}}

\end{picture}

}

\end{figure}
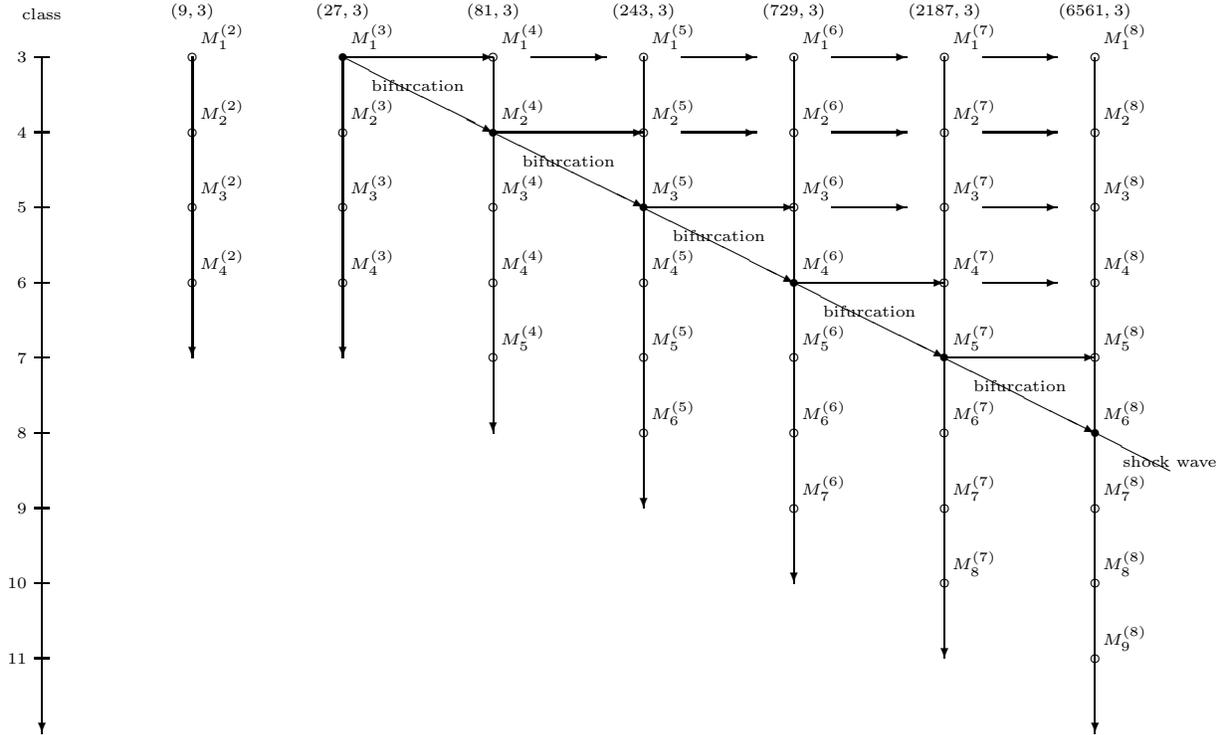

\begin{remark}
\label{rmk:ConstructionCF}
Formula
\eqref{eqn:Regular}
is the \textit{well-known old law} for the construction of the mainline
of coclass trees with elementary commutator quotient \(C_3\times C_3\).
Formulas
\eqref{eqn:Singular}
and
\eqref{eqn:Irregular}
constitute the \textit{new deterministic laws} in the finite region
on and behind the shock wave,
in the case of non-elementary commutator quotients \(C_{3^e}\times C_3\), \(e\ge 4\).
The statements are illuminated graphically in Figure
\ref{fig:MainlinesCF}.
\end{remark}

%\newpage
%--------------------------------------------------------------------------------

%--------------------------------------------------------------------------------

\noindent
In Figure
\ref{fig:MainlinesCF},
the \textit{nilpotency class} is selected as the unifying invariant
on the left hand scale,
since all coclass trees start at a root of class \(\mathrm{cl}(M^{(e)}_1)=3\),
for \(e\ge 2\). The trees are drawn for \(e\le 8\).

The mainline of the leftmost coclass tree \(\mathcal{T}^2(M^{(2)}_1)\) is actually
not involved in the propagation, since it is completely regular and endo-genetic:
see Figure
\ref{fig:Tree21a1AscioneA}.

Exceptionally,
the arrows of directed edges are drawn in \textit{reverse orientation},
in order to point out the ostensive direction of bifurcation and propagation
(irregular exo-genetic propagation in \textit{horizontal} direction,
singular exo-genetic propagation in \textit{diagonal} direction,
and regular endo-genetic propagation in \textit{vertical} direction).

Figure
\ref{fig:MainlinesCF}
impressively shows that the root \(M^{(3)}_1\) of the coclass tree \(\mathcal{T}^3(M^{(3)}_1)\)
is a common ancestor of all mainline vertices of all CF-coclass trees
\(\mathcal{T}^e(M^{(e)}_1)\), \(e\ge 3\), under investigation.
It is clear that \(M^{(3)}_1\) is infinitely capable (root of a coclass tree),
but the aforementioned fact emphasizes that \(M^{(3)}_1=\langle 729,7\rangle\) has the remarkable property
of being \textit{infinitely capable of higher order}.

%--------------------------------------------------------------------------------

\begin{proof} (Theorem
\ref{thm:ConstructionMain})
Let \(G=\langle x,y\rangle\) be a two-generated finite \(3\)-group.
Then we denote the main commutator by
\(s_2=\lbrack y,x\rbrack\)
and higher commutators by
\(\forall_{j=3}^{c+1}\) \(s_{j}=\lbrack s_{j-1},x\rbrack\), \(t_{j}=\lbrack s_{j-1},y\rbrack\).
If the commutator quotient is bicyclic
\(G/G^\prime\simeq C_{3^e}\times C_3\)
with one non-elementary component having logarithmic exponent \(e\ge 2\),
then we assume \(w^3\in G^\prime\) for \(w=x^{3^{e-1}}\) and \(y^3\in G^\prime\).

All mainline vertices involved in Theorem
\ref{thm:ConstructionMain}
possess a parametrized pc-presentation \(M^{(e)}_{i}=\)
\begin{equation}
\label{eqn:PresMain}
\begin{aligned}
\langle x,y \mid\ & x^{3^{e-1}}=w,\ w^3=1,\ y^3=1,\ \forall_{j=2}^{c-3}\ s_j^3=s_{j+2}^2s_{j+3},\ s_{c-2}^3=s_{c}^2,\ s_{c-1}^3=s_{c}^3=1, \\
                  & \forall_{j=3}^c\ s_j=t_j,\ s_{c+1}=t_{c+1}=1\rangle
\end{aligned}
\end{equation}
with two parameters, logarithmic exponent \(e\ge 2\), and nilpotency class \(c=\mathrm{cl}(M^{(e)}_i)=i+2\ge 3\).

Recall that only the \textit{commutator structure} enters
the recursive definition of the descending central series
\(\gamma_1(G)=G\), and \(\gamma_i(G)=\lbrack\gamma_{i-1}(G),G\rbrack\), for \(i\ge 2\),
but also the \textit{power structure} is included in the lower exponent \(p\)-central series
\(P_0(G)=G\), and \(P_i(G)=P_{i-1}(G)^3\cdot\lbrack P_{i-1}(G),G\rbrack\), for \(i\ge 1\).

Generally, for a descendant \(D\),
the parent is \(A=\pi(D)=D/\gamma_{c}(D)\), and
\(A_p=\pi_p(D)=D/P_{c_p-1}(D)\) is the \(p\)-parent,
where \(c=\mathrm{cl}(D)\) is the class,
and \(c_p=\mathrm{cl}_p(D)\) is the \(p\)-class.
Now we put \(D:=M^{(e)}_i\) and consider \textit{three situations}.

1. \textit{Behind} the shock wave:
If \(c<e\), i.e. \(i+2<e\) resp. \(i\le e-3\), then
\(\gamma_{c}(D)=\langle s_c\rangle\) 
and
\(P_{c_p-1}(D)=\langle w\rangle\).
Consequently, we obtain \(s_c=1\) in \(A=\pi(D)\) but \(w\) persists,
that is \(A=M^{(e)}_{i-1}\), if \(i\ge 2\).
However, in \(A_p=\pi_p(D)\), we get \(w=1\) but \(s_c\) persists,
that is \(A_p=M^{(e-1)}_{i}\), provided that \(e\ge 4\)
(for \(e\le 3\), the condition \(3\le c<e\) cannot occur).

2. \textit{On} the shock wave:
If \(c=e\), i.e. \(i+2=e\) resp. \(i=e-2\), then
\(\gamma_{c}(D)=\langle s_c\rangle\) 
and
\(P_{c_p-1}(D)=\langle s_c,w\rangle\) is bicyclic.
Thus, we have \(s_c=1\) but \(w\) persists in \(A\),
that is \(A=M^{(e)}_{i-1}\), if \(e\ge 4\).
However, in \(A_p\) both, \(s_c=1\) and \(w=1\), become trivial,
whence \(A_p=M^{(e-1)}_{i-1}\) with step size \(s=2\) reveals a bifurcation,
provided that \(e\ge 4\) and thus \(i=e-2\ge 2\).

3. \textit{Ahead of} the shock wave:
If \(c>e\), i.e. \(i+2>e\) resp. \(i\ge e-1\), then
\(c_p=c\) and \(\gamma_{c}(D)=P_{c_p-1}(D)=\langle s_c\rangle\).
So we get \(s_c=1\) in \(A=A_p\) but \(w\) persists,
that is \(A=A_p=M^{(e)}_{i-1}\), since \(i\ge 3-1=2\) for each \(e\ge 3\).

Strictly speaking, the preceding considerations only prove that
\(M^{(e)}_i=M^{(e-1)}_i-\#1;k\), resp.
\(M^{(e)}_{i}=M^{(e-1)}_{i-1}-\#2;\ell\), resp.
\(M^{(e)}_{i}=M^{(e)}_{i-1}-\#1;m\),
with positive integers \(k,\ell,m\),
but actual computations with Magma
\cite{MAGMA2021}
show that \(k=\ell=m=1\) for vertices on the mainline.
\end{proof}

%\newpage
%--------------------------------------------------------------------------------

\subsection{Vertices remote from the mainline with depth \(1\)}
\label{ss:OffsideCF}

\noindent
Concerning vertices \(V\)
on the coclass trees \(\mathcal{T}^e(M^{(e)}_1)\), \(e\ge 3\),
which are remote from the mainline,
we restrict ourselves to those with depth \(\mathrm{dp}(V)=1\)
and omit the investigation of others with depth \(2\le\mathrm{dp}(V)\le 3\).
Let \(V^{(e)}_i\) with \(i\ge 2\) be an \textit{offside} immediate descendant of a
mainline vertex \(M^{(e)}_{i-1}\),
and let \(\zeta\) be its centre.

%--------------------------------------------------------------------------------

\begin{theorem}
\label{thm:ConstructionOffside}
The vertices \(V^{(e)}_i\) remote from the mainline
of the coclass tree \(\mathcal{T}^e(M^{(e)}_1)\)
can be constructed recursively,
according to four laws in dependence on nilpotency class and centre \(\zeta\),
\begin{itemize}
\item
by \textbf{irregular exo}-genetic propagation
(behind the shock wave, \textbf{with stable type})
\begin{equation}
\label{eqn:IrregularOffside}
V^{(e)}_i=V^{(e-1)}_i-\#1;1, \text{ for } \zeta \text{ bicyclic},\ e\ge 5,\ 2\le i\le e-3, \text{ i.e. }\mathrm{cl}(V^{(e)}_i)<e,
\end{equation}
\item
by \textbf{singular exo}-genetic propagation
(bifurcation on the shock wave)
\begin{equation}
\label{eqn:SingularOffside}
V^{(e)}_{i}=M^{(e-1)}_{i-1}-\#2;\ell, \text{ for } \zeta \text{ bicyclic},\ e\ge 4,\ i=e-2, \text{ i.e. }\mathrm{cl}(V^{(e)}_i)=e,
\end{equation}
\item
by \textbf{regular endo}-genetic propagation
(ahead of the shock wave)
\begin{equation}
\label{eqn:RegularOffside}
V^{(e)}_{i}=M^{(e)}_{i-1}-\#1;m, \text{ for } \zeta \text{ bicyclic},\ e\ge 3,\ i\ge e-1, \text{ i.e. }\mathrm{cl}(V^{(e)}_i)>e,
\end{equation}
\item
by \textbf{permanent regular endo}-genetic propagation
(independent of the shock wave)
\begin{equation}
\label{eqn:PermanentOffside}
V^{(e)}_{i}=M^{(e)}_{i-1}-\#1;q, \text{ for } \zeta \text{ cyclic},\ e\ge 3,\ i\ge 2.
\end{equation}
\end{itemize}
\end{theorem}

%--------------------------------------------------------------------------------

\begin{proof}
For each \textit{periodic sequence} (also called \textit{coclass family}),
the vertices \(V\) have a parametrized pc-presentation with two parameters \(e\) and \(c\).
By the \textbf{mainline principle},
the generating commutator of the last non-trivial lower central
\(\gamma_c(V)=\langle s_c\rangle\)
does not enter the relations for the mainline,
but enters at least one typical relation, in \textbf{boldface} font, for each vertex off mainline.
The branches of the coclass trees under investigation are periodic with length \(2\). 
On every branch,
there is a unique mainline vertex \(M\) of type \(\mathrm{a}.1\), \(\varkappa(M)=(000;0)\).
We recall its pc-presentation: 
\begin{equation}
\label{eqn:a1MainPres}
\begin{aligned}
\langle x,y \mid\ & x^{3^{e-1}}=w,\ w^3=1,\ y^3=1,\ \forall_{j=2}^{c-3}\ s_j^3=s_{j+2}^2s_{j+3},\ s_{c-2}^3=s_{c}^2,\ s_{c-1}^3=s_{c}^3=1, \\
                  & \forall_{j=3}^c\ s_j=t_j,\ s_{c+1}=t_{c+1}=1\rangle.
\end{aligned}
\end{equation}
Furthermore, there is a unique leaf \(V\) of type \(\mathrm{b}.16\), \(\varkappa(V)\sim (004;0)\):
\begin{equation}
\label{eqn:b16Pres}
\begin{aligned}
\langle x,y \mid\ & x^{3^{e-1}}=w,\ w^3=1,\ \mathbf{y^3=s_c},\ \forall_{j=2}^{c-3}\ s_j^3=s_{j+2}^2s_{j+3},\ s_{c-2}^3=s_{c}^2,\ s_{c-1}^3=s_{c}^3=1, \\
                  & \forall_{j=3}^c\ s_j=t_j,\ s_{c+1}=t_{c+1}=1\rangle.
\end{aligned}
\end{equation}
On odd branches, we have a single leaf,
on even branches, we have two leaves, \(V\) of type \(\mathrm{b}.3\), \(\varkappa(V)\sim (001;0)\),
with \textit{cyclic} centre \(\zeta\) and exponent \(n=1\), resp. \(1\le n\le 2\):
\begin{equation}
\label{eqn:b3Pres}
\begin{aligned}
\langle x,y \mid\ & x^{3^{e-1}}=w,\ \mathbf{w^3=s_c^n},\ y^3=1,\ \forall_{j=2}^{c-3}\ s_j^3=s_{j+2}^2s_{j+3},\ s_{c-2}^3=s_{c}^2,\ s_{c-1}^3=s_{c}^3=1, \\
                  & \forall_{j=3}^c\ s_j=t_j,\ s_{c+1}=t_{c+1}=1\rangle.
\end{aligned}
\end{equation}
On every branch,
there is a unique root \(V\) of type \(\mathrm{a}.1\), \(\varkappa(V)=(000;0)\),
of a \textit{twig} which goes down to depth \(3\)
(we devote our attention to the root alone and abstain from its descendants):
\begin{equation}
\label{eqn:a1TwigPres}
\begin{aligned}
\langle x,y \mid\ & x^{3^{e-1}}=w,\ w^3=1,\ y^3=1,\ \forall_{j=2}^{c-3}\ s_j^3=s_{j+2}^2s_{j+3},\ s_{c-2}^3=s_{c}^2,\ s_{c-1}^3=s_{c}^3=1, \\
                  & \mathbf{t_3=s_3s_c},\ \forall_{j=4}^c\ s_j=t_j,\ s_{c+1}=t_{c+1}=1\rangle.
\end{aligned}
\end{equation}
On every branch,
there are two leaves \(V\) of type \(\mathrm{a}.1\), \(\varkappa(V)=(000;0)\),
with \textit{bicyclic} centre \(\zeta\) and exponent \(1\le n\le 2\):
\begin{equation}
\label{eqn:a1BicycPres}
\begin{aligned}
\langle x,y \mid\ & x^{3^{e-1}}=w,\ w^3=1,\ \mathbf{y^3=s_c^n},\ \forall_{j=2}^{c-3}\ s_j^3=s_{j+2}^2s_{j+3},\ s_{c-2}^3=s_{c}^2,\ s_{c-1}^3=s_{c}^3=1, \\
                  & \mathbf{t_3=s_3s_c},\ \forall_{j=4}^c\ s_j=t_j,\ s_{c+1}=t_{c+1}=1\rangle.
\end{aligned}
\end{equation}
On odd branches, we have a single leaf,
on even branches, we have two leaves, \(V\) of type \(\mathrm{a}.1\), \(\varkappa(V)=(000;0)\),
with \textit{cyclic} centre \(\zeta\) and exponent \(n=1\), resp. \(1\le n\le 2\):
\begin{equation}
\label{eqn:a1CycPres}
\begin{aligned}
\langle x,y \mid\ & x^{3^{e-1}}=w,\ \mathbf{w^3=s_c^n},\ y^3=1,\ \forall_{j=2}^{c-3}\ s_j^3=s_{j+2}^2s_{j+3},\ s_{c-2}^3=s_{c}^2,\ s_{c-1}^3=s_{c}^3=1, \\
                  & \mathbf{t_3=s_3s_c},\ \forall_{j=4}^c\ s_j=t_j,\ s_{c+1}=t_{c+1}=1\rangle.
\end{aligned}
\end{equation}

Similarly as in the proof of Theorem
\ref{thm:ConstructionMain},
for a descendant \(D\),
the parent is \(A=\pi(D)=D/\gamma_{c}(D)\), and
the \(p\)-parent is \(A_p=\pi_p(D)=D/P_{c_p-1}(D)\),
where \(c=\mathrm{cl}(D)\) is the class,
and \(c_p=\mathrm{cl}_p(D)\) is the \(p\)-class.
Now we put \(D:=V^{(e)}_i\) and consider \textit{four situations}.
For the first three items,
let \(D\) be a vertex with \textit{bicyclic centre}, and thus with one of the presentations
\eqref{eqn:b16Pres},
\eqref{eqn:a1TwigPres}
or
\eqref{eqn:a1BicycPres}.

1. \textit{Behind} the shock wave:
If \(c<e\), i.e. \(i+2<e\) resp. \(i\le e-3\), then
\(\gamma_{c}(D)=\langle s_c\rangle\) 
and
\(P_{c_p-1}(D)=\langle w\rangle\).
Consequently, we obtain \(s_c=1\) in \(A=\pi(D)\) but \(w\) persists,
that is \(A=M^{(e)}_{i-1}\), if \(i\ge 2\).
However, in \(A_p=\pi_p(D)\), we get \(w=1\)
but \(s_c\) (and the distinguished relation) persists,
that is \(A_p=V^{(e-1)}_{i}\), same type, provided that \(e\ge 5\)
(for \(e\le 4\), condition \(4\le c<e\) cannot occur).

2. \textit{On} the shock wave:
If \(c=e\), i.e. \(i+2=e\) resp. \(i=e-2\), then
\(\gamma_{c}(D)=\langle s_c\rangle\) 
and
\(P_{c_p-1}(D)=\langle s_c,w\rangle\) is bicyclic.
Thus, we have \(s_c=1\) but \(w\) persists in \(A\),
that is \(A=M^{(e)}_{i-1}\), if \(e\ge 4\).
However, in \(A_p\) both, \(s_c=1\) and \(w=1\), become trivial,
whence \(A_p=M^{(e-1)}_{i-1}\) with step size \(s=2\) reveals a bifurcation,
provided that \(e\ge 4\) and thus \(i=e-2\ge 2\).

3. \textit{Ahead of} the shock wave:
If \(c>e\), i.e. \(i+2>e\) resp. \(i\ge e-1\), then
\(c_p=c\) and \(\gamma_{c}(D)=P_{c_p-1}(D)=\langle s_c\rangle\).
So we get \(s_c=1\) (the distinguished relation degenerates) in \(A=A_p\) but \(w\) persists,
i.e. we get a mainline vertex \(A=A_p=M^{(e)}_{i-1}\), since \(i\ge 3-1=2\) for each \(e\ge 3\).

4.
Let \(D\) be a vertex with \textit{cyclic centre}, and thus with one of the presentations
\eqref{eqn:b3Pres}
or
\eqref{eqn:a1CycPres}.
Although the \(p\)-class \(c_p\) may be bigger than the class \(c\) of \(D\),
nevertheless, the last non-trivial lower centre and lower \(p\)-centre coincide
\(\gamma_{c}(D)=P_{c_p-1}(D)=\langle s_c\rangle\),
due to the \textit{exceptional relation} \(\mathbf{w^3=s_c^n}\).
So we get \(s_c=1\) in \(A=A_p\) and \(w^3=1\) becomes trivial,
that is, the propagation is always regular and endo-genetic
with coinciding parent and \(p\)-parent \(A=A_p=M^{(e)}_{i-1}\) a mainline vertex,
where \(i+2=c\ge 4\), i.e. \(i-1=c-3\ge 1\).

Strictly speaking, the preceding considerations only prove that
\(V^{(e)}_i=V^{(e-1)}_i-\#1;k\), resp.
\(V^{(e)}_{i}=M^{(e-1)}_{i-1}-\#2;\ell\), resp.
\(V^{(e)}_{i}=M^{(e)}_{i-1}-\#1;m\), resp.
\(V^{(e)}_{i}=M^{(e)}_{i-1}-\#1;q\),
with positive integers \(k,\ell,m,q\),
but actual computations with Magma
\cite{MAGMA2021}
show that \(k=1\), and \(\ell,m,q\ge 2\).
\end{proof}

%--------------------------------------------------------------------------------

\begin{remark}
\label{rmk:Polarization}
It should be pointed out that
the descendant vertices \(V\) of each root \(M^{(e)}_1\) share further invariants with the root,
aside from the rank distribution \(\varrho(V)\).
They have closely related transfer kernel types \(\varkappa(V)\)
with three identical components (the \textit{stabilization})
and a single varying component (the \textit{polarization}).
For all CF trees in this paper,
the polarization is located at the \textit{third} component,
and thus distinct from the \textit{puncture},
which is the fourth component, by convention.
\end{remark}

%\newpage
%--------------------------------------------------------------------------------

\section{Laws for coclass trees of BCF-groups}
\label{s:LawsBCF}

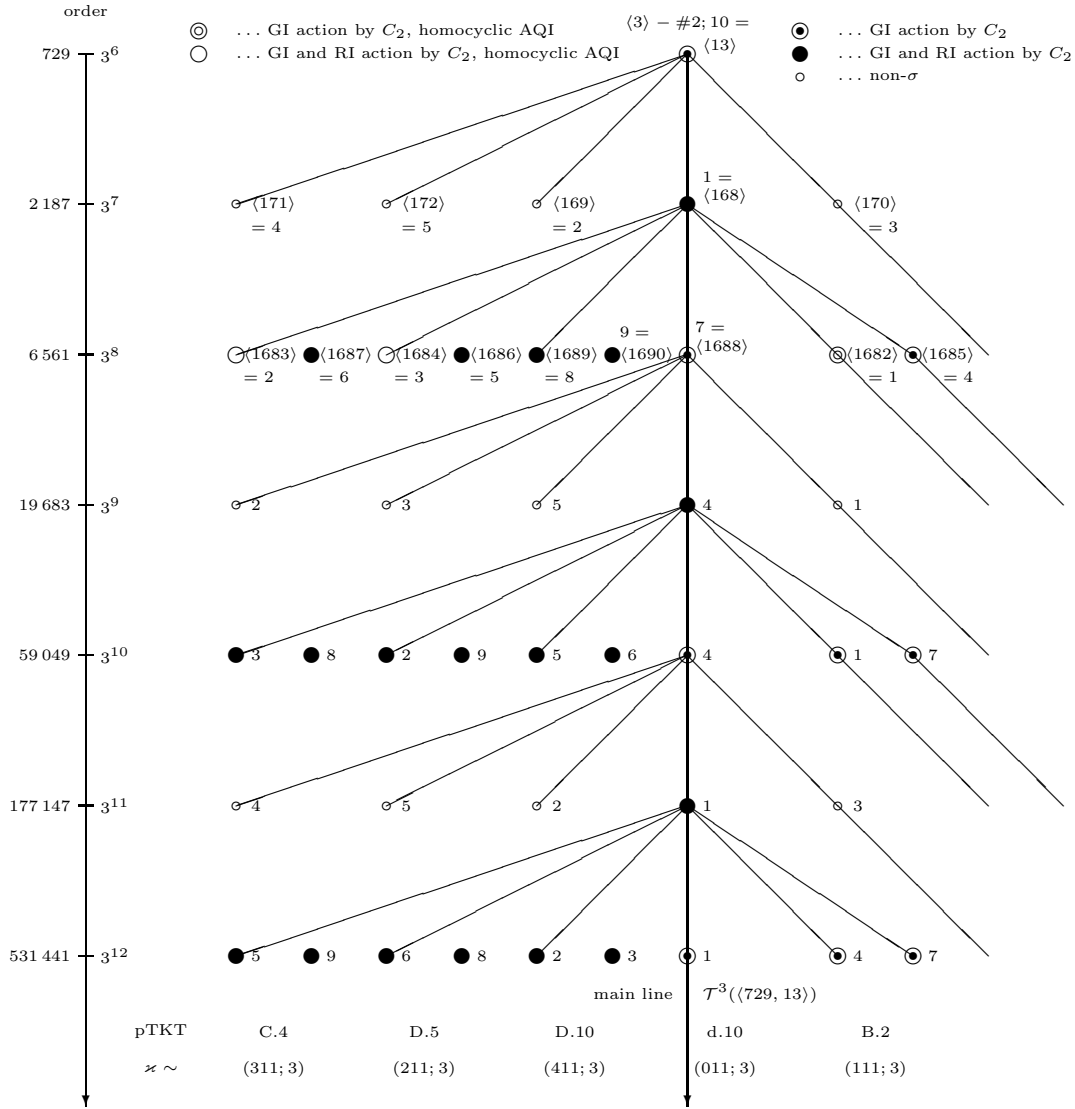
\begin{figure}[ht]
\caption{Depth-pruned metabelian coclass-\(3\) tree of type \(\mathrm{d}.10\) for AQI \((21)\)}
\label{fig:Tree21d10}

{\tiny

\setlength{\unitlength}{1cm}
\begin{picture}(14,15)(-6,-14)

% scale of orders
\put(-5,0.5){\makebox(0,0)[cb]{order}}

\put(-5,0){\line(0,-1){12}}
\multiput(-5.1,0)(0,-2){7}{\line(1,0){0.2}}

\put(-5.2,0){\makebox(0,0)[rc]{\(729\)}}
\put(-4.8,0){\makebox(0,0)[lc]{\(3^6\)}}
\put(-5.2,-2){\makebox(0,0)[rc]{\(2\,187\)}}
\put(-4.8,-2){\makebox(0,0)[lc]{\(3^7\)}}
\put(-5.2,-4){\makebox(0,0)[rc]{\(6\,561\)}}
\put(-4.8,-4){\makebox(0,0)[lc]{\(3^8\)}}
\put(-5.2,-6){\makebox(0,0)[rc]{\(19\,683\)}}
\put(-4.8,-6){\makebox(0,0)[lc]{\(3^9\)}}
\put(-5.2,-8){\makebox(0,0)[rc]{\(59\,049\)}}
\put(-4.8,-8){\makebox(0,0)[lc]{\(3^{10}\)}}
\put(-5.2,-10){\makebox(0,0)[rc]{\(177\,147\)}}
\put(-4.8,-10){\makebox(0,0)[lc]{\(3^{11}\)}}
\put(-5.2,-12){\makebox(0,0)[rc]{\(531\,441\)}}
\put(-4.8,-12){\makebox(0,0)[lc]{\(3^{12}\)}}

\put(-5,-12){\vector(0,-1){2}}

% legend
\put(-3.5,0.3){\circle{0.2}}
\put(-3.5,0.3){\circle{0.1}}
\put(-3,0.3){\makebox(0,0)[lc]{\(\ldots\) GI action by \(C_2\), homocyclic AQI}}
\put(-3.5,0){\circle{0.2}}
\put(-3,0){\makebox(0,0)[lc]{\(\ldots\) GI and RI action by \(C_2\), homocyclic AQI}}

\put(4.5,0.3){\circle{0.2}}
\put(4.5,0.3){\circle*{0.1}}
\put(5,0.3){\makebox(0,0)[lc]{\(\ldots\) GI action by \(C_2\)}}
\put(4.5,0){\circle*{0.2}}
\put(5,0){\makebox(0,0)[lc]{\(\ldots\) GI and RI action by \(C_2\)}}
\put(4.5,-0.3){\circle{0.1}}
\put(5,-0.3){\makebox(0,0)[lc]{\(\ldots\) non-\(\sigma\)}}

% infinite mainline
\put(3,0.3){\makebox(0,0)[cb]{\(\langle 3\rangle-\#2;10=\)}}
\put(3.2,0){\makebox(0,0)[lb]{\(\langle 13\rangle\)}}
\put(3.2,-1.7){\makebox(0,0)[lb]{\(1=\)}}
\put(3.2,-2){\makebox(0,0)[lb]{\(\langle 168\rangle\)}}
\put(3.1,-3.7){\makebox(0,0)[lb]{\(7=\)}}
\put(3.1,-4){\makebox(0,0)[lb]{\(\langle 1688\rangle\)}}
\multiput(3.2,-6)(0,-2){2}{\makebox(0,0)[lc]{\(4\)}}
\multiput(3.2,-10)(0,-2){2}{\makebox(0,0)[lc]{\(1\)}}

\multiput(3,0)(0,-4){4}{\circle{0.2}}
\multiput(3,0)(0,-4){4}{\circle*{0.1}}

\multiput(3,-2)(0,-4){3}{\circle*{0.2}}

\multiput(3,0)(0,-2){6}{\line(0,-1){2}}

\put(3,-12){\vector(0,-1){2}}
\put(2.8,-12.5){\makebox(0,0)[rc]{main line}}
\put(3.2,-12.5){\makebox(0,0)[lc]{\(\mathcal{T}^3(\langle 729,13\rangle)\)}}

% other periodic sequences
\put(-2.8,-2){\makebox(0,0)[lc]{\(\langle 171\rangle\)}}
\put(-2.8,-2.3){\makebox(0,0)[lc]{\(=4\)}}
\put(-0.8,-2){\makebox(0,0)[lc]{\(\langle 172\rangle\)}}
\put(-0.8,-2.3){\makebox(0,0)[lc]{\(=5\)}}
\put(1.2,-2){\makebox(0,0)[lc]{\(\langle 169\rangle\)}}
\put(1.2,-2.3){\makebox(0,0)[lc]{\(=2\)}}
\put(5.2,-2){\makebox(0,0)[lc]{\(\langle 170\rangle\)}}
\put(5.4,-2.3){\makebox(0,0)[lc]{\(=3\)}}

\put(-2.9,-4){\makebox(0,0)[lc]{\(\langle 1683\rangle\)}}
\put(-2.9,-4.3){\makebox(0,0)[lc]{\(=2\)}}
\put(-1.9,-4){\makebox(0,0)[lc]{\(\langle 1687\rangle\)}}
\put(-1.9,-4.3){\makebox(0,0)[lc]{\(=6\)}}
\put(-0.9,-4){\makebox(0,0)[lc]{\(\langle 1684\rangle\)}}
\put(-0.9,-4.3){\makebox(0,0)[lc]{\(=3\)}}
\put(0.1,-4){\makebox(0,0)[lc]{\(\langle 1686\rangle\)}}
\put(0.1,-4.3){\makebox(0,0)[lc]{\(=5\)}}
\put(1.1,-4){\makebox(0,0)[lc]{\(\langle 1689\rangle\)}}
\put(1.1,-4.3){\makebox(0,0)[lc]{\(=8\)}}
\put(2.1,-3.7){\makebox(0,0)[lc]{\(9=\)}}
\put(2.1,-4){\makebox(0,0)[lc]{\(\langle 1690\rangle\)}}
\put(5.1,-4){\makebox(0,0)[lc]{\(\langle 1682\rangle\)}}
\put(5.4,-4.3){\makebox(0,0)[lc]{\(=1\)}}
\put(6.1,-4){\makebox(0,0)[lc]{\(\langle 1685\rangle\)}}
\put(6.4,-4.3){\makebox(0,0)[lc]{\(=4\)}}

\put(-2.8,-6){\makebox(0,0)[lc]{\(2\)}}
\put(-2.8,-8){\makebox(0,0)[lc]{\(3\)}}
\put(-2.8,-10){\makebox(0,0)[lc]{\(4\)}}
\put(-2.8,-12){\makebox(0,0)[lc]{\(5\)}}

\put(-1.8,-8){\makebox(0,0)[lc]{\(8\)}}
\put(-1.8,-12){\makebox(0,0)[lc]{\(9\)}}

\put(-0.8,-6){\makebox(0,0)[lc]{\(3\)}}
\put(-0.8,-8){\makebox(0,0)[lc]{\(2\)}}
\put(-0.8,-10){\makebox(0,0)[lc]{\(5\)}}
\put(-0.8,-12){\makebox(0,0)[lc]{\(6\)}}

\put(0.2,-8){\makebox(0,0)[lc]{\(9\)}}
\put(0.2,-12){\makebox(0,0)[lc]{\(8\)}}

\put(1.2,-6){\makebox(0,0)[lc]{\(5\)}}
\put(1.2,-8){\makebox(0,0)[lc]{\(5\)}}
\put(1.2,-10){\makebox(0,0)[lc]{\(2\)}}
\put(1.2,-12){\makebox(0,0)[lc]{\(2\)}}

\put(2.2,-8){\makebox(0,0)[lc]{\(6\)}}
\put(2.2,-12){\makebox(0,0)[lc]{\(3\)}}

\put(5.2,-6){\makebox(0,0)[lc]{\(1\)}}
\put(5.2,-8){\makebox(0,0)[lc]{\(1\)}}
\put(5.2,-10){\makebox(0,0)[lc]{\(3\)}}
\put(5.2,-12){\makebox(0,0)[lc]{\(4\)}}

\put(6.2,-8){\makebox(0,0)[lc]{\(7\)}}
\put(6.2,-12){\makebox(0,0)[lc]{\(7\)}}

\multiput(-3,-2)(0,-4){3}{\circle{0.1}}

\put(-3,-4){\circle{0.2}}
\multiput(-3,-8)(0,-4){2}{\circle*{0.2}}
\multiput(-2,-4)(0,-4){3}{\circle*{0.2}}

\multiput(-1,-2)(0,-4){3}{\circle{0.1}}

\put(-1,-4){\circle{0.2}}
\multiput(-1,-8)(0,-4){2}{\circle*{0.2}}
\multiput(0,-4)(0,-4){3}{\circle*{0.2}}

\multiput(1,-2)(0,-4){3}{\circle{0.1}}

\multiput(1,-4)(0,-4){3}{\circle*{0.2}}

\multiput(2,-4)(0,-4){3}{\circle*{0.2}}

\multiput(5,-2)(0,-4){3}{\circle{0.1}}

\multiput(5,-4)(0,-4){3}{\circle{0.2}}
\put(5,-4){\circle{0.1}}
\multiput(5,-8)(0,-4){2}{\circle*{0.1}}
\multiput(6,-4)(0,-4){3}{\circle{0.2}}
\multiput(6,-4)(0,-4){3}{\circle*{0.1}}

% directed edges
\multiput(3,0)(0,-2){6}{\line(-3,-1){6}}
\multiput(3,0)(0,-2){6}{\line(-2,-1){4}}
\multiput(3,0)(0,-2){6}{\line(-1,-1){2}}
\multiput(3,0)(0,-2){6}{\line(1,-1){2}}
\multiput(3,-2)(0,-4){3}{\line(3,-2){3}}

% edges to vertices of depth two, three, four
\multiput(5,-2)(0,-2){5}{\line(1,-1){2}}
\multiput(6,-4)(0,-4){2}{\line(1,-1){2}}

% punctured transfer kernel types
\put(-4,-13){\makebox(0,0)[cc]{pTKT}}
\put(-4,-13.5){\makebox(0,0)[cc]{\(\varkappa\sim\)}}
\put(-2.5,-13){\makebox(0,0)[cc]{\(\mathrm{C}.4\)}}
\put(-2.5,-13.5){\makebox(0,0)[cc]{\((311;3)\)}}
\put(-0.5,-13){\makebox(0,0)[cc]{\(\mathrm{D}.5\)}}
\put(-0.5,-13.5){\makebox(0,0)[cc]{\((211;3)\)}}
\put(1.5,-13){\makebox(0,0)[cc]{\(\mathrm{D}.10\)}}
\put(1.5,-13.5){\makebox(0,0)[cc]{\((411;3)\)}}
\put(3.5,-13){\makebox(0,0)[cc]{\(\mathrm{d}.10\)}}
\put(3.5,-13.5){\makebox(0,0)[cc]{\((011;3)\)}}
\put(5.5,-13){\makebox(0,0)[cc]{\(\mathrm{B}.2\)}}
\put(5.5,-13.5){\makebox(0,0)[cc]{\((111;3)\)}}

\end{picture}

}

\end{figure}

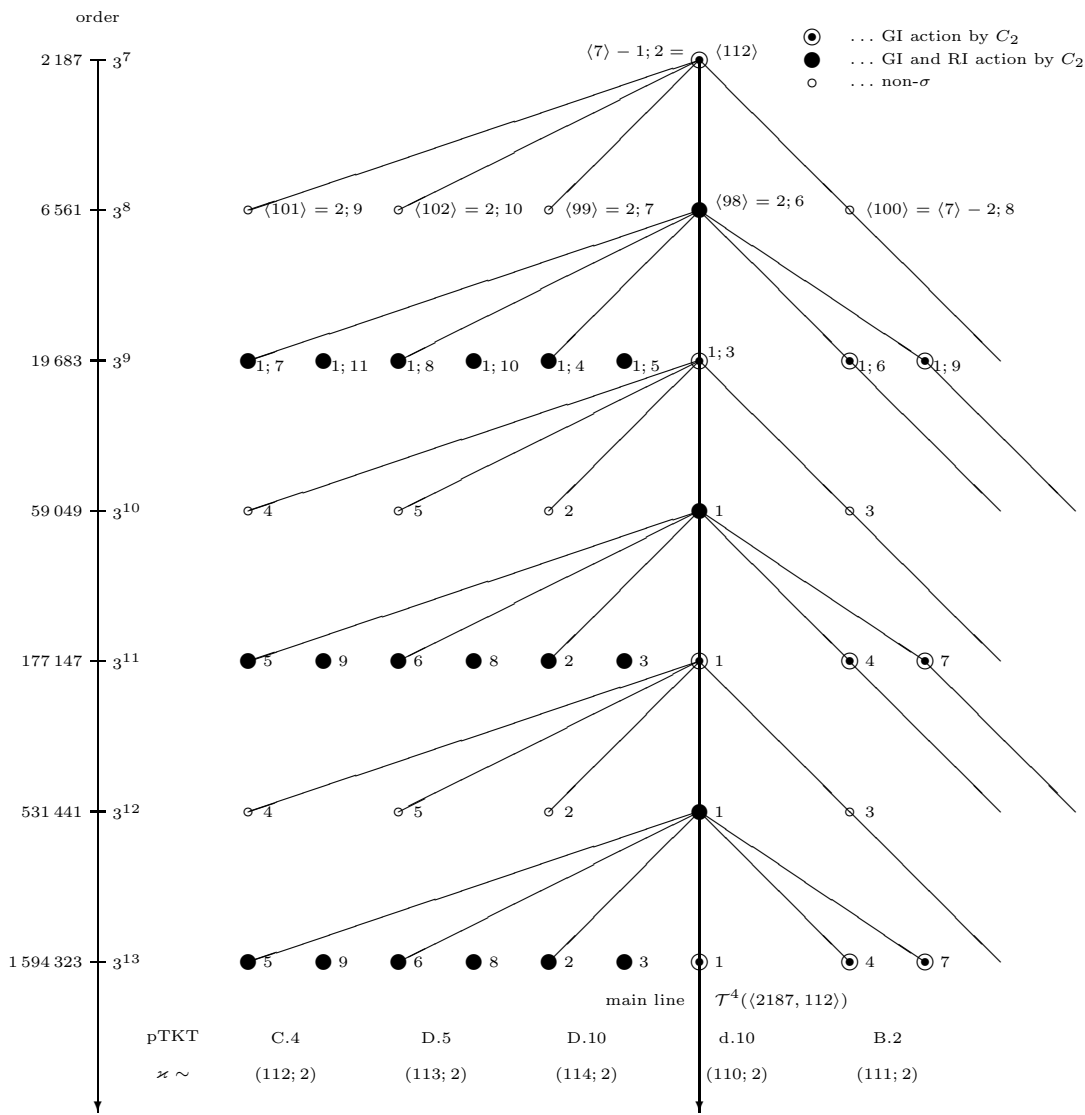
\begin{figure}[ht]
\caption{Depth-pruned metabelian coclass-\(4\) tree of type \(\mathrm{d}.10\) for AQI \((31)\)}
\label{fig:Tree31d10}

{\tiny

\setlength{\unitlength}{1cm}
\begin{picture}(14,15)(-6,-14)

% scale of orders
\put(-5,0.5){\makebox(0,0)[cb]{order}}

\put(-5,0){\line(0,-1){12}}
\multiput(-5.1,0)(0,-2){7}{\line(1,0){0.2}}

\put(-5.2,0){\makebox(0,0)[rc]{\(2\,187\)}}
\put(-4.8,0){\makebox(0,0)[lc]{\(3^7\)}}
\put(-5.2,-2){\makebox(0,0)[rc]{\(6\,561\)}}
\put(-4.8,-2){\makebox(0,0)[lc]{\(3^8\)}}
\put(-5.2,-4){\makebox(0,0)[rc]{\(19\,683\)}}
\put(-4.8,-4){\makebox(0,0)[lc]{\(3^9\)}}
\put(-5.2,-6){\makebox(0,0)[rc]{\(59\,049\)}}
\put(-4.8,-6){\makebox(0,0)[lc]{\(3^{10}\)}}
\put(-5.2,-8){\makebox(0,0)[rc]{\(177\,147\)}}
\put(-4.8,-8){\makebox(0,0)[lc]{\(3^{11}\)}}
\put(-5.2,-10){\makebox(0,0)[rc]{\(531\,441\)}}
\put(-4.8,-10){\makebox(0,0)[lc]{\(3^{12}\)}}
\put(-5.2,-12){\makebox(0,0)[rc]{\(1\,594\,323\)}}
\put(-4.8,-12){\makebox(0,0)[lc]{\(3^{13}\)}}

\put(-5,-12){\vector(0,-1){2}}

% legend
\put(4.5,0.3){\circle{0.2}}
\put(4.5,0.3){\circle*{0.1}}
\put(5,0.3){\makebox(0,0)[lc]{\(\ldots\) GI action by \(C_2\)}}
\put(4.5,0){\circle*{0.2}}
\put(5,0){\makebox(0,0)[lc]{\(\ldots\) GI and RI action by \(C_2\)}}
\put(4.5,-0.3){\circle{0.1}}
\put(5,-0.3){\makebox(0,0)[lc]{\(\ldots\) non-\(\sigma\)}}

% infinite mainline
\put(2.8,0){\makebox(0,0)[rb]{\(\langle 7\rangle-1;2=\)}}
\put(3.2,0){\makebox(0,0)[lb]{\(\langle 112\rangle\)}}
\put(3.2,-2){\makebox(0,0)[lb]{\(\langle 98\rangle=2;6\)}}
\put(3.1,-4){\makebox(0,0)[lb]{\(1;3\)}}
\multiput(3.2,-6)(0,-2){4}{\makebox(0,0)[lc]{\(1\)}}

\multiput(3,0)(0,-4){4}{\circle{0.2}}
\multiput(3,0)(0,-4){4}{\circle*{0.1}}

\multiput(3,-2)(0,-4){3}{\circle*{0.2}}

\multiput(3,0)(0,-2){6}{\line(0,-1){2}}

\put(3,-12){\vector(0,-1){2}}
\put(2.8,-12.5){\makebox(0,0)[rc]{main line}}
\put(3.2,-12.5){\makebox(0,0)[lc]{\(\mathcal{T}^4(\langle 2187,112\rangle)\)}}

% other periodic sequences
\put(-2.8,-2){\makebox(0,0)[lc]{\(\langle 101\rangle=2;9\)}}
\put(-0.8,-2){\makebox(0,0)[lc]{\(\langle 102\rangle=2;10\)}}
\put(1.2,-2){\makebox(0,0)[lc]{\(\langle 99\rangle=2;7\)}}
\put(5.2,-2){\makebox(0,0)[lc]{\(\langle 100\rangle=\langle 7\rangle-2;8\)}}

\put(-2.9,-4){\makebox(0,0)[lt]{\(1;7\)}}
\put(-1.9,-4){\makebox(0,0)[lt]{\(1;11\)}}
\put(-0.9,-4){\makebox(0,0)[lt]{\(1;8\)}}
\put(0.1,-4){\makebox(0,0)[lt]{\(1;10\)}}
\put(1.1,-4){\makebox(0,0)[lt]{\(1;4\)}}
\put(2.1,-4){\makebox(0,0)[lt]{\(1;5\)}}
\put(5.1,-4){\makebox(0,0)[lt]{\(1;6\)}}
\put(6.1,-4){\makebox(0,0)[lt]{\(1;9\)}}

\multiput(-2.8,-6)(0,-4){2}{\makebox(0,0)[lc]{\(4\)}}
\multiput(-0.8,-6)(0,-4){2}{\makebox(0,0)[lc]{\(5\)}}
\multiput(1.2,-6)(0,-4){2}{\makebox(0,0)[lc]{\(2\)}}
\multiput(5.2,-6)(0,-4){2}{\makebox(0,0)[lc]{\(3\)}}

\multiput(-2.8,-8)(0,-4){2}{\makebox(0,0)[lc]{\(5\)}}
\multiput(-1.8,-8)(0,-4){2}{\makebox(0,0)[lc]{\(9\)}}
\multiput(-0.8,-8)(0,-4){2}{\makebox(0,0)[lc]{\(6\)}}
\multiput(0.2,-8)(0,-4){2}{\makebox(0,0)[lc]{\(8\)}}
\multiput(1.2,-8)(0,-4){2}{\makebox(0,0)[lc]{\(2\)}}
\multiput(2.2,-8)(0,-4){2}{\makebox(0,0)[lc]{\(3\)}}
\multiput(5.2,-8)(0,-4){2}{\makebox(0,0)[lc]{\(4\)}}
\multiput(6.2,-8)(0,-4){2}{\makebox(0,0)[lc]{\(7\)}}

\multiput(-3,-2)(0,-4){3}{\circle{0.1}}

%\put(-3,-4){\circle{0.2}}
\multiput(-3,-4)(0,-4){3}{\circle*{0.2}}
\multiput(-2,-4)(0,-4){3}{\circle*{0.2}}

\multiput(-1,-2)(0,-4){3}{\circle{0.1}}

%\put(-1,-4){\circle{0.2}}
\multiput(-1,-4)(0,-4){3}{\circle*{0.2}}
\multiput(0,-4)(0,-4){3}{\circle*{0.2}}

\multiput(1,-2)(0,-4){3}{\circle{0.1}}

\multiput(1,-4)(0,-4){3}{\circle*{0.2}}

\multiput(2,-4)(0,-4){3}{\circle*{0.2}}

\multiput(5,-2)(0,-4){3}{\circle{0.1}}

\multiput(5,-4)(0,-4){3}{\circle{0.2}}
%\put(5,-4){\circle{0.1}}
\multiput(5,-4)(0,-4){3}{\circle*{0.1}}
\multiput(6,-4)(0,-4){3}{\circle{0.2}}
\multiput(6,-4)(0,-4){3}{\circle*{0.1}}

% directed edges
\multiput(3,0)(0,-2){6}{\line(-3,-1){6}}
\multiput(3,0)(0,-2){6}{\line(-2,-1){4}}
\multiput(3,0)(0,-2){6}{\line(-1,-1){2}}
\multiput(3,0)(0,-2){6}{\line(1,-1){2}}
\multiput(3,-2)(0,-4){3}{\line(3,-2){3}}

% edges to vertices of depth two, three, four
\multiput(5,-2)(0,-2){5}{\line(1,-1){2}}
\multiput(6,-4)(0,-4){2}{\line(1,-1){2}}

% punctured transfer kernel types
\put(-4,-13){\makebox(0,0)[cc]{pTKT}}
\put(-4,-13.5){\makebox(0,0)[cc]{\(\varkappa\sim\)}}
\put(-2.5,-13){\makebox(0,0)[cc]{\(\mathrm{C}.4\)}}
\put(-2.5,-13.5){\makebox(0,0)[cc]{\((112;2)\)}}
\put(-0.5,-13){\makebox(0,0)[cc]{\(\mathrm{D}.5\)}}
\put(-0.5,-13.5){\makebox(0,0)[cc]{\((113;2)\)}}
\put(1.5,-13){\makebox(0,0)[cc]{\(\mathrm{D}.10\)}}
\put(1.5,-13.5){\makebox(0,0)[cc]{\((114;2)\)}}
\put(3.5,-13){\makebox(0,0)[cc]{\(\mathrm{d}.10\)}}
\put(3.5,-13.5){\makebox(0,0)[cc]{\((110;2)\)}}
\put(5.5,-13){\makebox(0,0)[cc]{\(\mathrm{B}.2\)}}
\put(5.5,-13.5){\makebox(0,0)[cc]{\((111;2)\)}}

\end{picture}

}

\end{figure}

\begin{proposition}
\label{prp:UniquenessBCF}
For each log exponent \(e\ge 2\),
there exists a unique coclass tree \(\mathcal{T}^{e+1}(\mathbb{M}^{(e+1)}_1)\ni\mathbb{V}\)
with fixed coclass \(\mathrm{cc}(\mathbb{V})=e+1\),
fixed commutator quotient \(\mathbb{V}/\mathbb{V}^\prime\simeq C_{3^e}\times C_3\),
and fixed rank distribution \(\varrho(\mathbb{V})\sim (2,2,3;3)\).
Its mainline \((\mathbb{M}^{(e+1)}_i)_{i\ge 1}\)
is of type \(\mathrm{d}.10\), \(\varkappa(\mathbb{M}^{(e+1)}_i)\sim (110;2)\).
The tree contains metabelian and non-metabelian \(\mathrm{BCF}\)-groups.
The branches are of depth \(4\). \\
(See Figures
\ref{fig:Tree21d10}
--
\ref{fig:Tree41d10}
for the depth-pruned metabelian skeleton, when \(2\le e\le 4\).)
\end{proposition}

\begin{proof}
According to the proof of Proposition
\ref{prp:UniquenessCF},
there are only two coclass trees \(\mathcal{T}^r(R^{(r)}_j)\)
with rank distribution \(\varrho(V)\sim (2,2,3;3)\),
for each \(e\ge 2\), \(e\le r\le e+1\),
a BCF-tree of type \(\mathrm{d}.10\) and a CF-tree of type \(\mathrm{a}.1\).
The former is the \textit{unique} tree with root \(R^{(e+1)}_j=\mathbb{M}^{(e+1)}_1\),
recursively determined by the CF-group \(M^{(3)}_1=\langle 729,7\rangle\)
and the BCF-group \(\mathbb{M}^{(3)}_1=\langle 729,13\rangle\).
Its depth-pruned metabelian branches are periodic of length \(2\) without pre-period,
and all of its vertices are BCF-groups,
since the vertices of the first two branches are BCF-groups.
\end{proof}

%--------------------------------------------------------------------------------

\begin{proposition}
\label{prp:InvariantsBCF}
For \(e\ge 2\),
invariants of
vertices on the mainline \((\mathbb{M}^{(e+1)}_i)_{i\ge 1}\)
of the coclass tree \(\mathcal{T}^{e+1}(\mathbb{M}^{(e+1)}_1)\)
are given as follows:
\begin{equation}
\label{eqn:InvariantsBCF}
\begin{aligned}
\text{log order } \mathrm{lo}(\mathbb{M}^{(e+1)}_i)=e+i+3,
\text{ nilpotency class } \mathrm{cl}(\mathbb{M}^{(e+1)}_i) &= i+2, \text{ for } i\ge 1, \\
p\text{-class } \mathrm{cl}_p(\mathbb{M}^{(e+1)}_i)=
\begin{cases}
i+2 & \text{ if } i>e-1, \\
e+1 & \text{ if } i\le e-1,
\end{cases}
p\text{-coclass } \mathrm{cc}_p(\mathbb{M}^{(e+1)}_i) &=
\begin{cases}
e+1 & \text{ if } i>e-1, \\
i+2 & \text{ if } i\le e-1.
\end{cases}
\end{aligned}
\end{equation}
\end{proposition}

\begin{proof}
Proposition
\ref{prp:InvariantsBCF}
remains true when the mainline vertex \(\mathbb{M}^{(e+1)}_i\)
is replaced by any proper descendant vertex
\(\mathbb{V}^{(e+1)}_i\) with \(i\ge 2\).
All coclass trees under investigation
start at a root of class \(\mathrm{cl}(\mathbb{M}^{(e+1)}_1)=3=1+2\),
for each \(e\ge 2\).
Thus, proper descendants possess nilpotency class \(\mathrm{cl}(\mathbb{V}^{(e+1)}_i)=i+2\ge 4\).
By definition, all vertices \(\mathbb{V}\)
of the coclass tree \(\mathcal{T}^{e+1}(\mathbb{M}^{(e+1)}_1)\)
share the common coclass \(\mathrm{cc}(\mathbb{V})=e+1\).
Consequently, the logarithmic order is the sum
\(\mathrm{lo}(\mathbb{V}^{(e+1)}_i)=\mathrm{cl}(\mathbb{V}^{(e+1)}_i)+\mathrm{cc}(\mathbb{V}^{(e+1)}_i)=i+2+e+1\).
Finally, the \textit{power structure} of all finite \(3\)-groups \(G\)
with commutator quotient \(G/G^\prime\simeq C_{3^e}\times C_3\) is responsible for
the constant \(p\)-class \(\mathrm{cl}_p(\mathbb{V}^{(e+1)}_i)=e+1\),
independently of the class \(\mathrm{cl}(\mathbb{V}^{(e+1)}_i)=i+2\le e+1\),
in the finite region on and behind the shock wave.
\end{proof}

%--------------------------------------------------------------------------------

Concerning vertices \(\mathbb{V}\)
on the coclass trees \(\mathcal{T}^{e+1}(\mathbb{M}^{(e+1)}_1)\), \(e\ge 2\),
which are remote from the mainline,
we restrict ourselves to the metabelian with depth \(\mathrm{dp}(\mathbb{V})=1\),
and we omit the investigation of others with depth \(2\le\mathrm{dp}(\mathbb{V})\le 4\).

Let \(\mathbb{V}^{(e+1)}_i\) with \(i\ge 1\) be
any metabelian vertex \textit{on or remote from} the mainline \((\mathbb{M}^{(e+1)}_{i})_{i\ge 1}\).

\begin{theorem}
\label{thm:ConstructionBCF}
The vertices \(\mathbb{V}^{(e+1)}_i\) on and remote from the mainline
of the coclass tree \(\mathcal{T}^{e+1}(\mathbb{M}^{(e+1)}_1)\)
can be constructed recursively,
according to three laws in dependence on the nilpotency class,
\begin{itemize}
\item
by \textbf{irregular endo}-genetic propagation
(behind the shock wave, \textbf{with type change})
\begin{equation}
\label{eqn:IrregularBCF}
\mathbb{V}^{(e+1)}_i=V^{(e)}_i-\#1;k, \text{ for } e\ge 3,\ 1\le i\le e-2, \text{ i.e. }\mathrm{cl}(\mathbb{V}^{(e+1)}_i)<e+1,
\end{equation}
where \(\mathrm{d}.10\), \(\mathrm{B}.2\), \(\mathrm{D}.10\), \(\mathrm{C}.4\), \(\mathrm{D}.5\) of \(\mathbb{V}^{(e+1)}_i\)
correspond to \(\mathrm{a}.1\), \(\mathrm{a}.1\), \(\mathrm{b}.16\), \(\mathrm{a}.1\), \(\mathrm{a}.1\) of \(V^{(e)}_i\),
\item
by \textbf{singular exo}-genetic propagation
(bifurcation on the shock wave)
\begin{equation}
\label{eqn:SingularBCF}
\mathbb{V}^{(e+1)}_{i}=M^{(e)}_{i-1}-\#2;\ell, \text{ for } e\ge 3,\ i=e-1, \text{ i.e. }\mathrm{cl}(\mathbb{V}^{(e+1)}_i)=e+1,
\end{equation}
\item
by \textbf{regular endo}-genetic propagation
(ahead of the shock wave)
\begin{equation}
\label{eqn:RegularBCF}
\mathbb{V}^{(e+1)}_{i}=\mathbb{M}^{(e+1)}_{i-1}-\#1;m, \text{ for } e\ge 2,\ i\ge e, \text{ i.e. }\mathrm{cl}(\mathbb{V}^{(e+1)}_i)>e+1,
\end{equation}
\end{itemize}
\end{theorem}

\begin{proof}
For each \textit{periodic sequence} (or \textit{coclass family}),
the vertices \(\mathbb{V}\) have a parametrized pc-presentation with two parameters \(e\) and \(c\).
According to the \textbf{mainline principle},
the generating commutator of the last non-trivial lower central
\(\gamma_c(\mathbb{V})=\langle s_c\rangle\)
does not enter the relations for the mainline,
but enters at least one typical relation, in \textbf{boldface} font, for each vertex off mainline.
Branches of coclass trees under investigation are periodic with length \(2\). 
On every branch,
there is a unique mainline vertex \(\mathbb{M}\) of type \(\mathrm{d}.10\), \(\varkappa(\mathbb{M})\sim (110;2)\).
Its pc-presentation is given by: 
\begin{equation}
\label{eqn:d10MainPres}
\begin{aligned}
\langle x,y \mid\ & x^{3^{e}}=w,\ w^3=1,\ y^3=1,\ \forall_{j=2}^{c-3}\ s_j^3=s_{j+2}^2s_{j+3},\ s_{c-2}^3=s_{c}^2,\ s_{c-1}^3=s_{c}^3=1, \\
                  & t_3=s_3w,\ \forall_{j=4}^c\ s_j=t_{j},\ s_{c+1}=t_{c+1}=1\rangle.
\end{aligned}
\end{equation}
On odd branches, we have a single vertex,
on even branches, we have two vertices, \(\mathbb{V}\) of type \(\mathrm{D}.10\), \(\varkappa(\mathbb{V})\sim(114;2)\),
with exponent \(n=1\), resp. \(1\le n\le 2\):
\begin{equation}
\label{eqn:D10Pres}
\begin{aligned}
\langle x,y \mid\ & x^{3^{e}}=w,\ w^3=1,\ \mathbf{y^3=s_c^n},\ \forall_{j=2}^{c-3}\ s_j^3=s_{j+2}^2s_{j+3},\ s_{c-2}^3=s_{c}^2,\ s_{c-1}^3=s_{c}^3=1, \\
                  & t_3=s_3w,\ \forall_{j=4}^c\ s_j=t_{j},\ s_{c+1}=t_{c+1}=1\rangle.
\end{aligned}
\end{equation}
On odd branches, we have a single root,
on even branches, we have two roots, \(\mathbb{V}\) of type \(\mathrm{B}.2\), \(\varkappa(\mathbb{V})\sim (111;2)\),
of a twig, with exponent \(n=1\), resp. \(1\le n\le 2\):
\begin{equation}
\label{eqn:B2Pres}
\begin{aligned}
\langle x,y \mid\ & x^{3^{e}}=w,\ w^3=1,\ y^3=1,\ \forall_{j=2}^{c-3}\ s_j^3=s_{j+2}^2s_{j+3},\ s_{c-2}^3=s_{c}^2,\ s_{c-1}^3=s_{c}^3=1, \\
                  & \mathbf{t_3=s_3s_c^nw},\ \forall_{j=4}^c\ s_j=t_{j},\ s_{c+1}=t_{c+1}=1\rangle.
\end{aligned}
\end{equation}
On odd branches, we have a single vertex,
on even branches, we have two vertices, \(\mathbb{V}\) of type \(\mathrm{C}.4\), \(\varkappa(\mathbb{V})=(112;2)\),
with exponent \(n=1\), resp. \(1\le n\le 2\):
\begin{equation}
\label{eqn:C4Pres}
\begin{aligned}
\langle x,y \mid\ & x^{3^{e}}=w,\ w^3=1,\ \mathbf{y^3=s_c^n},\ \forall_{j=2}^{c-3}\ s_j^3=s_{j+2}^2s_{j+3},\ s_{c-2}^3=s_{c}^2,\ s_{c-1}^3=s_{c}^3=1, \\
                  & \mathbf{t_3=s_3s_c^nw},\ \forall_{j=4}^c\ s_j=t_{j},\ s_{c+1}=t_{c+1}=1\rangle.
\end{aligned}
\end{equation}
On odd branches, we have a single vertex,
on even branches, we have two vertices, \(\mathbb{V}\) of type \(\mathrm{D}.5\), \(\varkappa(\mathbb{V})=(113;2)\),
with exponent \(n=1\), resp. \(1\le n\le 2\):
\begin{equation}
\label{eqn:D5Pres}
\begin{aligned}
\langle x,y \mid\ & x^{3^{e}}=w,\ w^3=1,\ \mathbf{y^3=s_c^n},\ \forall_{j=2}^{c-3}\ s_j^3=s_{j+2}^2s_{j+3},\ s_{c-2}^3=s_{c}^2,\ s_{c-1}^3=s_{c}^3=1, \\
                  & \mathbf{t_3=s_3s_c^{3-n}w},\ \forall_{j=4}^c\ s_j=t_{j},\ s_{c+1}=t_{c+1}=1\rangle.
\end{aligned}
\end{equation}

Similarly as in the proof of Theorem
\ref{thm:ConstructionOffside},
for a descendant \(D\),
the parent is \(A=\pi(D)=D/\gamma_{c}(D)\), and
\(A_p=\pi_p(D)=D/P_{c_p-1}(D)\) is the \(p\)-parent,
where \(c=\mathrm{cl}(D)\) is the class,
and \(c_p=\mathrm{cl}_p(D)\) is the \(p\)-class.
Now we put \(D:=\mathbb{V}^{(e+1)}_i\) and consider \textit{three situations}.

1. \textit{Behind} the shock wave:
If \(c<e+1\), i.e. \(i+2<e+1\) resp. \(i\le e-2\), then
\(\gamma_{c}(D)=\langle s_c,t_c\rangle\), if \(i=1\),
\(\gamma_{c}(D)=\langle s_c\rangle\), if \(i\ge 2\),
and
\(P_{c_p-1}(D)=\langle w\rangle\).
Consequently, if \(i\ge 2\), we obtain \(s_c=1\) in \(A=\pi(D)\) but \(w\) persists, all Formulas
\eqref{eqn:d10MainPres}, \eqref{eqn:D10Pres}, \eqref{eqn:B2Pres}, \eqref{eqn:C4Pres}, \eqref{eqn:D5Pres}
degenerate to
\eqref{eqn:d10MainPres},
that is \(A=\mathbb{M}^{(e+1)}_{i-1}\), if \(i\ge 2\).
However, in \(A_p=\pi_p(D)\), we get \(w=1\)
but \(s_c\) (and the distinguished relation) persists, Formula
\eqref{eqn:d10MainPres}, resp. \eqref{eqn:D10Pres}, \eqref{eqn:B2Pres}, \eqref{eqn:C4Pres}, \eqref{eqn:D5Pres},
becomes Formula
\eqref{eqn:a1MainPres}, resp. \eqref{eqn:b16Pres}, \eqref{eqn:a1TwigPres}, \eqref{eqn:a1BicycPres}, \eqref{eqn:a1BicycPres},
i.e. \(A_p=V^{(e)}_{i}\), provided \(e\ge 4\)
(for \(e\le 3\), condition \(4\le c<e+1\) cannot occur)
or \(e=3\) and \(c=3\).

2. \textit{On} the shock wave:
If \(c=e+1\), i.e. \(i+2=e+1\) resp. \(i=e-1\), then
\(\gamma_{c}(D)=\langle s_c\rangle\) 
and
\(P_{c_p-1}(D)=\langle s_c,w\rangle\) is bicyclic.
Thus, we have \(s_c=1\) but \(w\) persists in \(A\),
that is \(A=\mathbb{M}^{(e+1)}_{i-1}\), if \(e\ge 3\).
However, in \(A_p\) both, \(s_c=1\) and \(w=1\), become trivial, all Formulas
\eqref{eqn:d10MainPres}, \eqref{eqn:D10Pres}, \eqref{eqn:B2Pres}, \eqref{eqn:C4Pres}, \eqref{eqn:D5Pres}
degenerate to
\eqref{eqn:a1MainPres},
whence \(A_p=M^{(e)}_{i-1}\) with step size \(s=2\) reveals a bifurcation,
provided that \(e\ge 3\) and thus \(i=e-1\ge 2\).

3. \textit{Ahead of} the shock wave:
If \(c>e+1\), i.e. \(i+2>e+1\) resp. \(i\ge e\), then
\(c_p=c\) and \(\gamma_{c}(D)=P_{c_p-1}(D)=\langle s_c\rangle\).
So we get \(s_c=1\) in \(A=A_p\) but \(w\) persists, all Formulas
\eqref{eqn:d10MainPres}, \eqref{eqn:D10Pres}, \eqref{eqn:B2Pres}, \eqref{eqn:C4Pres}, \eqref{eqn:D5Pres}
degenerate to
\eqref{eqn:d10MainPres},
that is a mainline vertex \(A=A_p=\mathbb{M}^{(e+1)}_{i-1}\), since \(i\ge 2\) for each \(e\ge 2\).

In the preceding, we have proved that
\(\mathbb{V}^{(e+1)}_i=V^{(e)}_i-\#1;k\), resp.
\(\mathbb{V}^{(e+1)}_{i}=M^{(e)}_{i-1}-\#2;\ell\), resp.
\(\mathbb{V}^{(e+1)}_{i}=\mathbb{M}^{(e+1)}_{i-1}-\#1;m\),
with positive integers \(k,\ell,m\).
\end{proof}

%\newpage
%--------------------------------------------------------------------------------

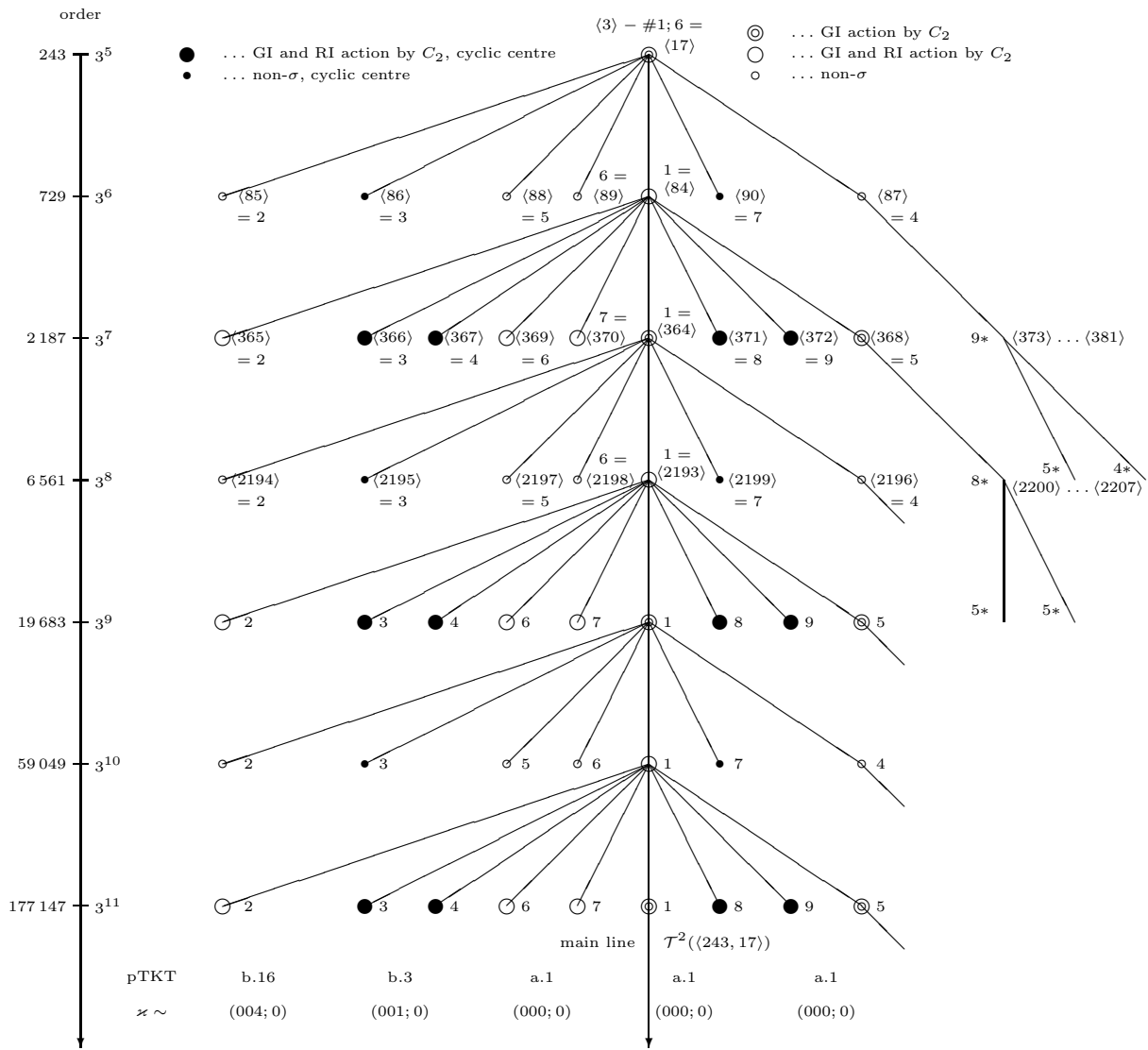
\begin{figure}[ht]
\caption{Coclass-\(2\) tree of type \(\mathrm{a}.1\), rooted in original Ascione \(A\), for AQI \((21)\)}
\label{fig:Tree21a1AscioneA}

{\tiny

\setlength{\unitlength}{1cm}
\begin{picture}(14,14.7)(-5,-13.7)

% scale of orders
\put(-5,0.5){\makebox(0,0)[cb]{order}}

\put(-5,0){\line(0,-1){12}}
\multiput(-5.1,0)(0,-2){7}{\line(1,0){0.2}}

\put(-5.2,0){\makebox(0,0)[rc]{\(243\)}}
\put(-4.8,0){\makebox(0,0)[lc]{\(3^5\)}}
\put(-5.2,-2){\makebox(0,0)[rc]{\(729\)}}
\put(-4.8,-2){\makebox(0,0)[lc]{\(3^6\)}}
\put(-5.2,-4){\makebox(0,0)[rc]{\(2\,187\)}}
\put(-4.8,-4){\makebox(0,0)[lc]{\(3^7\)}}
\put(-5.2,-6){\makebox(0,0)[rc]{\(6\,561\)}}
\put(-4.8,-6){\makebox(0,0)[lc]{\(3^8\)}}
\put(-5.2,-8){\makebox(0,0)[rc]{\(19\,683\)}}
\put(-4.8,-8){\makebox(0,0)[lc]{\(3^9\)}}
\put(-5.2,-10){\makebox(0,0)[rc]{\(59\,049\)}}
\put(-4.8,-10){\makebox(0,0)[lc]{\(3^{10}\)}}
\put(-5.2,-12){\makebox(0,0)[rc]{\(177\,147\)}}
\put(-4.8,-12){\makebox(0,0)[lc]{\(3^{11}\)}}

\put(-5,-12){\vector(0,-1){2}}

% legend
\put(-3.5,0){\circle*{0.2}}
\put(-3,0){\makebox(0,0)[lc]{\(\ldots\) GI and RI action by \(C_2\), cyclic centre}}
\put(-3.5,-0.3){\circle*{0.1}}
\put(-3,-0.3){\makebox(0,0)[lc]{\(\ldots\) non-\(\sigma\), cyclic centre}}

\put(4.5,0.3){\circle{0.2}}
\put(4.5,0.3){\circle{0.1}}
\put(5,0.3){\makebox(0,0)[lc]{\(\ldots\) GI action by \(C_2\)}}
\put(4.5,0){\circle{0.2}}
\put(5,0){\makebox(0,0)[lc]{\(\ldots\) GI and RI action by \(C_2\)}}
\put(4.5,-0.3){\circle{0.1}}
\put(5,-0.3){\makebox(0,0)[lc]{\(\ldots\) non-\(\sigma\)}}

% infinite mainline
\put(3,0.3){\makebox(0,0)[cb]{\(\langle 3\rangle-\#1;6=\)}}
\put(3.2,0){\makebox(0,0)[lb]{\(\langle 17\rangle\)}}
\put(3.2,-1.7){\makebox(0,0)[lb]{\(1=\)}}
\put(3.2,-2){\makebox(0,0)[lb]{\(\langle 84\rangle\)}}
\put(3.2,-3.7){\makebox(0,0)[lb]{\(1=\)}}
\put(3.1,-4){\makebox(0,0)[lb]{\(\langle 364\rangle\)}}
\put(3.2,-5.7){\makebox(0,0)[lb]{\(1=\)}}
\put(3.1,-6){\makebox(0,0)[lb]{\(\langle 2193\rangle\)}}
\multiput(3.2,-8)(0,-2){3}{\makebox(0,0)[lc]{\(1\)}}

\multiput(3,0)(0,-4){4}{\circle{0.2}}
\multiput(3,0)(0,-4){4}{\circle{0.1}}

\multiput(3,-2)(0,-4){3}{\circle{0.2}}

\multiput(3,0)(0,-2){6}{\line(0,-1){2}}

\put(3,-12){\vector(0,-1){2}}
\put(2.8,-12.5){\makebox(0,0)[rc]{main line}}
\put(3.2,-12.5){\makebox(0,0)[lc]{\(\mathcal{T}^2(\langle 243,17\rangle)\)}}

% other periodic sequences
\put(-2.8,-2){\makebox(0,0)[lc]{\(\langle 85\rangle\)}}
\put(-2.8,-2.3){\makebox(0,0)[lc]{\(=2\)}}
\put(-0.8,-2){\makebox(0,0)[lc]{\(\langle 86\rangle\)}}
\put(-0.8,-2.3){\makebox(0,0)[lc]{\(=3\)}}
\put(1.2,-2){\makebox(0,0)[lc]{\(\langle 88\rangle\)}}
\put(1.2,-2.3){\makebox(0,0)[lc]{\(=5\)}}
\put(2.3,-1.7){\makebox(0,0)[lc]{\(6=\)}}
\put(2.2,-2){\makebox(0,0)[lc]{\(\langle 89\rangle\)}}
\put(4.2,-2){\makebox(0,0)[lc]{\(\langle 90\rangle\)}}
\put(4.2,-2.3){\makebox(0,0)[lc]{\(=7\)}}
\put(6.2,-2){\makebox(0,0)[lc]{\(\langle 87\rangle\)}}
\put(6.4,-2.3){\makebox(0,0)[lc]{\(=4\)}}

\put(-2.9,-4){\makebox(0,0)[lc]{\(\langle 365\rangle\)}}
\put(-2.8,-4.3){\makebox(0,0)[lc]{\(=2\)}}
\put(-0.9,-4){\makebox(0,0)[lc]{\(\langle 366\rangle\)}}
\put(-0.8,-4.3){\makebox(0,0)[lc]{\(=3\)}}
\put(0.1,-4){\makebox(0,0)[lc]{\(\langle 367\rangle\)}}
\put(0.2,-4.3){\makebox(0,0)[lc]{\(=4\)}}
\put(1.1,-4){\makebox(0,0)[lc]{\(\langle 369\rangle\)}}
\put(1.2,-4.3){\makebox(0,0)[lc]{\(=6\)}}
\put(2.3,-3.7){\makebox(0,0)[lc]{\(7=\)}}
\put(2.1,-4){\makebox(0,0)[lc]{\(\langle 370\rangle\)}}
\put(4.1,-4){\makebox(0,0)[lc]{\(\langle 371\rangle\)}}
\put(4.2,-4.3){\makebox(0,0)[lc]{\(=8\)}}
\put(5.1,-4){\makebox(0,0)[lc]{\(\langle 372\rangle\)}}
\put(5.2,-4.3){\makebox(0,0)[lc]{\(=9\)}}
\put(6.1,-4){\makebox(0,0)[lc]{\(\langle 368\rangle\)}}
\put(6.4,-4.3){\makebox(0,0)[lc]{\(=5\)}}

\put(-2.9,-6){\makebox(0,0)[lc]{\(\langle 2194\rangle\)}}
\put(-2.8,-6.3){\makebox(0,0)[lc]{\(=2\)}}
\put(-0.9,-6){\makebox(0,0)[lc]{\(\langle 2195\rangle\)}}
\put(-0.8,-6.3){\makebox(0,0)[lc]{\(=3\)}}
\put(1.1,-6){\makebox(0,0)[lc]{\(\langle 2197\rangle\)}}
\put(1.2,-6.3){\makebox(0,0)[lc]{\(=5\)}}
\put(2.3,-5.7){\makebox(0,0)[lc]{\(6=\)}}
\put(2.1,-6){\makebox(0,0)[lc]{\(\langle 2198\rangle\)}}
\put(4.1,-6){\makebox(0,0)[lc]{\(\langle 2199\rangle\)}}
\put(4.2,-6.3){\makebox(0,0)[lc]{\(=7\)}}
\put(6.1,-6){\makebox(0,0)[lc]{\(\langle 2196\rangle\)}}
\put(6.4,-6.3){\makebox(0,0)[lc]{\(=4\)}}

\multiput(-2.7,-8)(0,-2){3}{\makebox(0,0)[lc]{\(2\)}}
\multiput(-0.8,-8)(0,-2){3}{\makebox(0,0)[lc]{\(3\)}}
\multiput(0.2,-8)(0,-4){2}{\makebox(0,0)[lc]{\(4\)}}
\multiput(1.2,-8)(0,-4){2}{\makebox(0,0)[lc]{\(6\)}}
\multiput(1.2,-10)(0,-4){1}{\makebox(0,0)[lc]{\(5\)}}
\multiput(2.2,-8)(0,-4){2}{\makebox(0,0)[lc]{\(7\)}}
\multiput(2.2,-10)(0,-4){1}{\makebox(0,0)[lc]{\(6\)}}
\multiput(4.2,-8)(0,-4){2}{\makebox(0,0)[lc]{\(8\)}}
\multiput(4.2,-10)(0,-4){1}{\makebox(0,0)[lc]{\(7\)}}
\multiput(5.2,-8)(0,-4){2}{\makebox(0,0)[lc]{\(9\)}}
\multiput(6.2,-8)(0,-4){2}{\makebox(0,0)[lc]{\(5\)}}
\multiput(6.2,-10)(0,-4){1}{\makebox(0,0)[lc]{\(4\)}}

% vertices of depth two and three
\put(7.8,-4){\makebox(0,0)[rc]{\(9\ast\)}}
\put(8.1,-4){\makebox(0,0)[lc]{\(\langle 373\rangle\ldots\langle 381\rangle\)}}
\put(8.8,-5.9){\makebox(0,0)[rb]{\(5\ast\)}}
\put(9.8,-5.9){\makebox(0,0)[rb]{\(4\ast\)}}

\put(7.8,-6){\makebox(0,0)[rc]{\(8\ast\)}}
\put(8.1,-6){\makebox(0,0)[lt]{\(\langle 2200\rangle\ldots\langle 2207\rangle\)}}
\put(7.8,-7.9){\makebox(0,0)[rb]{\(5\ast\)}}
\put(8.8,-7.9){\makebox(0,0)[rb]{\(5\ast\)}}

% vertices of depth one
\multiput(-3,-2)(0,-4){3}{\circle{0.1}}
\multiput(-3,-4)(0,-4){3}{\circle{0.2}}

\multiput(-1,-2)(0,-4){3}{\circle*{0.1}}
\multiput(-1,-4)(0,-4){3}{\circle*{0.2}}

\multiput(0,-4)(0,-4){3}{\circle*{0.2}}

\multiput(1,-2)(0,-4){3}{\circle{0.1}}
\multiput(1,-4)(0,-4){3}{\circle{0.2}}

\multiput(2,-2)(0,-4){3}{\circle{0.1}}
\multiput(2,-4)(0,-4){3}{\circle{0.2}}

\multiput(4,-2)(0,-4){3}{\circle*{0.1}}
\multiput(4,-4)(0,-4){3}{\circle*{0.2}}

\multiput(5,-4)(0,-4){3}{\circle*{0.2}}

\multiput(6,-2)(0,-4){3}{\circle{0.1}}
\multiput(6,-4)(0,-4){3}{\circle{0.1}}
\multiput(6,-4)(0,-4){3}{\circle{0.2}}

% directed edges to depth one
\multiput(3,0)(0,-2){6}{\line(-3,-1){6}}
\multiput(3,0)(0,-2){6}{\line(-2,-1){4}}
\multiput(3,-2)(0,-4){3}{\line(-3,-2){3}}
\multiput(3,0)(0,-2){6}{\line(-1,-1){2}}
\multiput(3,0)(0,-2){6}{\line(-1,-2){1}}
\multiput(3,0)(0,-2){6}{\line(1,-2){1}}
\multiput(3,-2)(0,-4){3}{\line(1,-1){2}}
\multiput(3,0)(0,-2){6}{\line(3,-2){3}}

% edges to depth two and three
\multiput(6,-2)(0,-2){2}{\line(1,-1){2}}
\multiput(8,-4)(0,-2){2}{\line(1,-2){1}}
\put(8,-4){\line(1,-1){2}}
\put(8,-6){\line(0,-1){2}}

\multiput(6,-6)(0,-2){4}{\line(1,-1){0.6}}

% punctured transfer kernel types
\put(-4,-13){\makebox(0,0)[cc]{pTKT}}
\put(-4,-13.5){\makebox(0,0)[cc]{\(\varkappa\sim\)}}
\put(-2.5,-13){\makebox(0,0)[cc]{\(\mathrm{b}.16\)}}
\put(-2.5,-13.5){\makebox(0,0)[cc]{\((004;0)\)}}
\put(-0.5,-13){\makebox(0,0)[cc]{\(\mathrm{b}.3\)}}
\put(-0.5,-13.5){\makebox(0,0)[cc]{\((001;0)\)}}
\put(1.5,-13){\makebox(0,0)[cc]{\(\mathrm{a}.1\)}}
\put(1.5,-13.5){\makebox(0,0)[cc]{\((000;0)\)}}
\put(3.5,-13){\makebox(0,0)[cc]{\(\mathrm{a}.1\)}}
\put(3.5,-13.5){\makebox(0,0)[cc]{\((000;0)\)}}
\put(5.5,-13){\makebox(0,0)[cc]{\(\mathrm{a}.1\)}}
\put(5.5,-13.5){\makebox(0,0)[cc]{\((000;0)\)}}

\end{picture}

}

\end{figure}

%--------------------------------------------------------------------------------

\noindent
In Figure
\ref{fig:Tree21a1AscioneA},
the coclass tree
with \textit{Ascione's CF-group} \(A\) as its root
\(\langle 243,17\rangle=\langle 81,3\rangle-\#1;6\)
is drawn up to order \(3^{11}\).
The branches are periodic with length \(2\)
and naturally bounded depth \(3\),
without artificial pruning.
The infinite main line is of type \(\mathrm{a}.1\)
and consists of \(\sigma\)-groups with GI action by \(C_2\).
The other vertices are \(\sigma\)-groups for even branches
and non-\(\sigma\) groups for odd branches.
Vertices with positive depth are of type \(\mathrm{b}.16\)
or \(\mathrm{b}.3\) or \(\mathrm{a}.1\).
All vertices of type \(\mathrm{b}.3\) and some of type \(\mathrm{a}.1\)
have a cyclic centre \(\zeta\simeq C_9\).
Vertices of depth \(\mathrm{dp}\in\lbrace 2,3\rbrace\)
are exclusively of type \(\mathrm{a}.1\) with cyclic centre \(\zeta\simeq C_3\).
They are drawn for the first and second branch only.
On the first branch,
\(\langle 373\rangle\) gives rise to \(\langle 2208\rangle\) \(\ldots\) \(\langle 2212\rangle\), and
\(\langle 378\rangle\) gives rise to \(\langle 2213\rangle\) \(\ldots\) \(\langle 2216\rangle\).
On the second branch,
\(\langle 2200\rangle\) gives rise to \(5\) non-\(\sigma\)-descendants, and
\(\langle 2205\rangle\) gives rise to \(5\) descendants with generator inverting (GI) and relator inverting (RI) action by \(C_2\).

%--------------------------------------------------------------------------------

\begin{remark}
\label{rmk:Tree21a1AscioneA}
The purely graph theoretic structure of the coclass tree in Figure
\ref{fig:Tree21a1AscioneA}
was indicated in
\cite[Tbl. 6, p. 272]{AHL1977}
with much less details and only up to order \(3^8\).
The root \(A=\langle 243,17\rangle\) was described ten years earlier by James
\cite{Jm1968}
as a member of the first branch \(\Phi_3(1)\) of Hall's isoclinism class \(\Phi_3\).
See the pc-presentation for \(\Phi_3(2111)e\simeq A\) in
\cite[p. 620]{Jm1980}.
Another pc-presentation for \(A\) can be extracted from the microfiches in
\cite[p. 320, folio I02]{AHL1977}.
\end{remark}

%\newpage
%--------------------------------------------------------------------------------

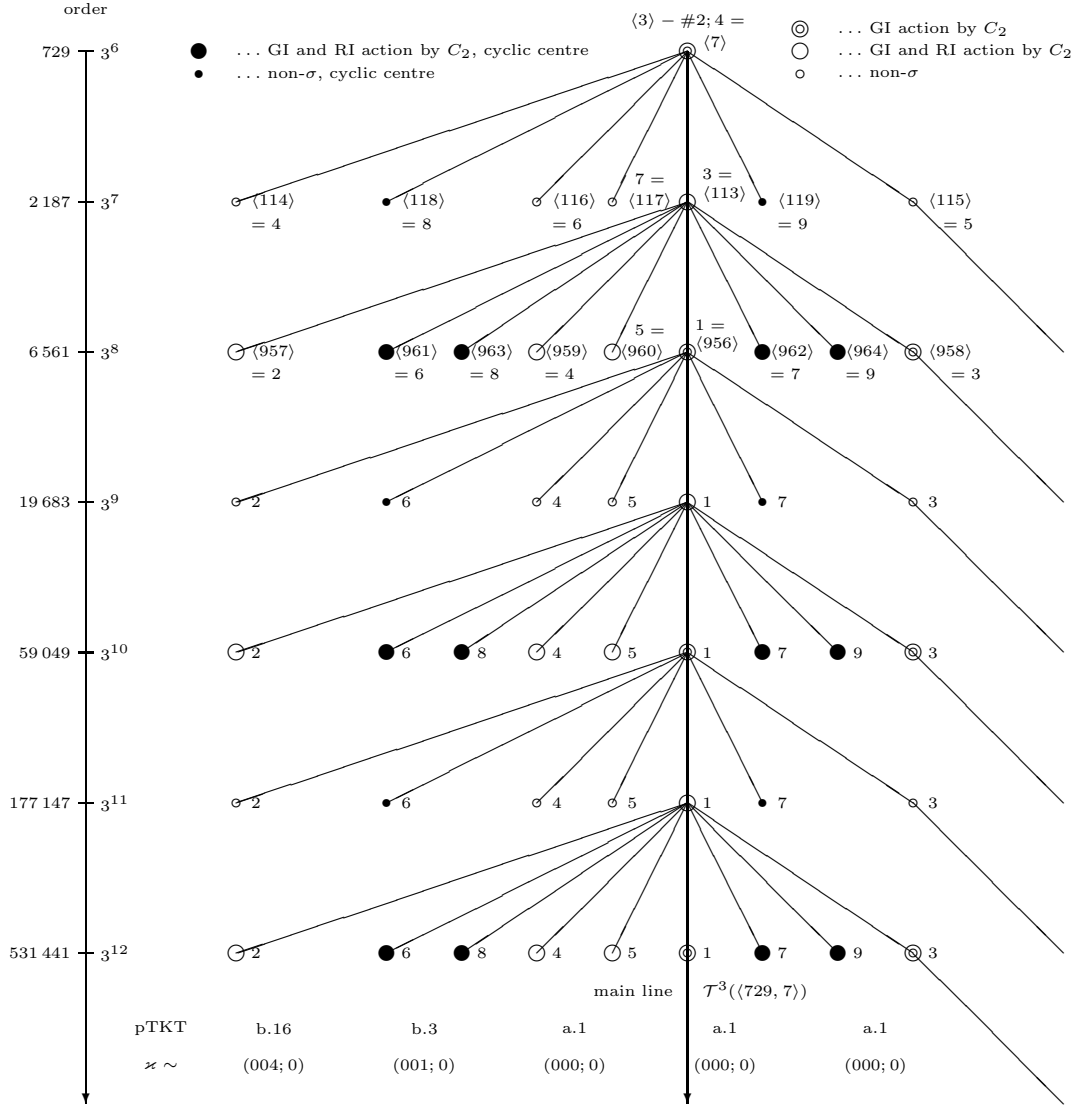
\begin{figure}[ht]
\caption{Coclass-\(3\) tree of type \(\mathrm{a}.1\), rooted in generalized Ascione A, for AQI \((31)\)}
\label{fig:Tree31a1AscioneA}

{\tiny

\setlength{\unitlength}{1cm}
\begin{picture}(14,14.7)(-6,-13.7)

% scale of orders
\put(-5,0.5){\makebox(0,0)[cb]{order}}

\put(-5,0){\line(0,-1){12}}
\multiput(-5.1,0)(0,-2){7}{\line(1,0){0.2}}

\put(-5.2,0){\makebox(0,0)[rc]{\(729\)}}
\put(-4.8,0){\makebox(0,0)[lc]{\(3^6\)}}
\put(-5.2,-2){\makebox(0,0)[rc]{\(2\,187\)}}
\put(-4.8,-2){\makebox(0,0)[lc]{\(3^7\)}}
\put(-5.2,-4){\makebox(0,0)[rc]{\(6\,561\)}}
\put(-4.8,-4){\makebox(0,0)[lc]{\(3^8\)}}
\put(-5.2,-6){\makebox(0,0)[rc]{\(19\,683\)}}
\put(-4.8,-6){\makebox(0,0)[lc]{\(3^9\)}}
\put(-5.2,-8){\makebox(0,0)[rc]{\(59\,049\)}}
\put(-4.8,-8){\makebox(0,0)[lc]{\(3^{10}\)}}
\put(-5.2,-10){\makebox(0,0)[rc]{\(177\,147\)}}
\put(-4.8,-10){\makebox(0,0)[lc]{\(3^{11}\)}}
\put(-5.2,-12){\makebox(0,0)[rc]{\(531\,441\)}}
\put(-4.8,-12){\makebox(0,0)[lc]{\(3^{12}\)}}

\put(-5,-12){\vector(0,-1){2}}

% legend
\put(-3.5,0){\circle*{0.2}}
\put(-3,0){\makebox(0,0)[lc]{\(\ldots\) GI and RI action by \(C_2\), cyclic centre}}
\put(-3.5,-0.3){\circle*{0.1}}
\put(-3,-0.3){\makebox(0,0)[lc]{\(\ldots\) non-\(\sigma\), cyclic centre}}

\put(4.5,0.3){\circle{0.2}}
\put(4.5,0.3){\circle{0.1}}
\put(5,0.3){\makebox(0,0)[lc]{\(\ldots\) GI action by \(C_2\)}}
\put(4.5,0){\circle{0.2}}
\put(5,0){\makebox(0,0)[lc]{\(\ldots\) GI and RI action by \(C_2\)}}
\put(4.5,-0.3){\circle{0.1}}
\put(5,-0.3){\makebox(0,0)[lc]{\(\ldots\) non-\(\sigma\)}}

% infinite mainline
\put(3,0.3){\makebox(0,0)[cb]{\(\langle 3\rangle-\#2;4=\)}}
\put(3.2,0){\makebox(0,0)[lb]{\(\langle 7\rangle\)}}
\put(3.2,-1.7){\makebox(0,0)[lb]{\(3=\)}}
\put(3.2,-2){\makebox(0,0)[lb]{\(\langle 113\rangle\)}}
\put(3.1,-3.7){\makebox(0,0)[lb]{\(1=\)}}
\put(3.1,-4){\makebox(0,0)[lb]{\(\langle 956\rangle\)}}
\multiput(3.2,-6)(0,-2){4}{\makebox(0,0)[lc]{\(1\)}}

\multiput(3,0)(0,-4){4}{\circle{0.2}}
\multiput(3,0)(0,-4){4}{\circle{0.1}}

\multiput(3,-2)(0,-4){3}{\circle{0.2}}

\multiput(3,0)(0,-2){6}{\line(0,-1){2}}

\put(3,-12){\vector(0,-1){2}}
\put(2.8,-12.5){\makebox(0,0)[rc]{main line}}
\put(3.2,-12.5){\makebox(0,0)[lc]{\(\mathcal{T}^3(\langle 729,7\rangle)\)}}

% other periodic sequences
\put(-2.8,-2){\makebox(0,0)[lc]{\(\langle 114\rangle\)}}
\put(-2.8,-2.3){\makebox(0,0)[lc]{\(=4\)}}
\put(-0.8,-2){\makebox(0,0)[lc]{\(\langle 118\rangle\)}}
\put(-0.8,-2.3){\makebox(0,0)[lc]{\(=8\)}}
\put(1.2,-2){\makebox(0,0)[lc]{\(\langle 116\rangle\)}}
\put(1.2,-2.3){\makebox(0,0)[lc]{\(=6\)}}
\put(2.3,-1.7){\makebox(0,0)[lc]{\(7=\)}}
\put(2.2,-2){\makebox(0,0)[lc]{\(\langle 117\rangle\)}}
\put(4.2,-2){\makebox(0,0)[lc]{\(\langle 119\rangle\)}}
\put(4.2,-2.3){\makebox(0,0)[lc]{\(=9\)}}
\put(6.2,-2){\makebox(0,0)[lc]{\(\langle 115\rangle\)}}
\put(6.4,-2.3){\makebox(0,0)[lc]{\(=5\)}}

\put(-2.8,-4){\makebox(0,0)[lc]{\(\langle 957\rangle\)}}
\put(-2.8,-4.3){\makebox(0,0)[lc]{\(=2\)}}
\put(-0.9,-4){\makebox(0,0)[lc]{\(\langle 961\rangle\)}}
\put(-0.9,-4.3){\makebox(0,0)[lc]{\(=6\)}}
\put(0.1,-4){\makebox(0,0)[lc]{\(\langle 963\rangle\)}}
\put(0.1,-4.3){\makebox(0,0)[lc]{\(=8\)}}
\put(1.1,-4){\makebox(0,0)[lc]{\(\langle 959\rangle\)}}
\put(1.1,-4.3){\makebox(0,0)[lc]{\(=4\)}}
\put(2.3,-3.7){\makebox(0,0)[lc]{\(5=\)}}
\put(2.1,-4){\makebox(0,0)[lc]{\(\langle 960\rangle\)}}
\put(4.1,-4){\makebox(0,0)[lc]{\(\langle 962\rangle\)}}
\put(4.1,-4.3){\makebox(0,0)[lc]{\(=7\)}}
\put(5.1,-4){\makebox(0,0)[lc]{\(\langle 964\rangle\)}}
\put(5.1,-4.3){\makebox(0,0)[lc]{\(=9\)}}
\put(6.2,-4){\makebox(0,0)[lc]{\(\langle 958\rangle\)}}
\put(6.5,-4.3){\makebox(0,0)[lc]{\(=3\)}}

\multiput(-2.8,-6)(0,-2){4}{\makebox(0,0)[lc]{\(2\)}}
\multiput(-0.8,-6)(0,-2){4}{\makebox(0,0)[lc]{\(6\)}}
\multiput(0.2,-8)(0,-4){2}{\makebox(0,0)[lc]{\(8\)}}
\multiput(1.2,-6)(0,-2){4}{\makebox(0,0)[lc]{\(4\)}}
\multiput(2.2,-6)(0,-2){4}{\makebox(0,0)[lc]{\(5\)}}
\multiput(4.2,-6)(0,-2){4}{\makebox(0,0)[lc]{\(7\)}}
\multiput(5.2,-8)(0,-4){2}{\makebox(0,0)[lc]{\(9\)}}
\multiput(6.2,-6)(0,-2){4}{\makebox(0,0)[lc]{\(3\)}}

% vertices of depth one
\multiput(-3,-2)(0,-4){3}{\circle{0.1}}
\multiput(-3,-4)(0,-4){3}{\circle{0.2}}

\multiput(-1,-2)(0,-4){3}{\circle*{0.1}}
\multiput(-1,-4)(0,-4){3}{\circle*{0.2}}

\multiput(0,-4)(0,-4){3}{\circle*{0.2}}

\multiput(1,-2)(0,-4){3}{\circle{0.1}}
\multiput(1,-4)(0,-4){3}{\circle{0.2}}

\multiput(2,-2)(0,-4){3}{\circle{0.1}}
\multiput(2,-4)(0,-4){3}{\circle{0.2}}

\multiput(4,-2)(0,-4){3}{\circle*{0.1}}
\multiput(4,-4)(0,-4){3}{\circle*{0.2}}

\multiput(5,-4)(0,-4){3}{\circle*{0.2}}

\multiput(6,-2)(0,-4){3}{\circle{0.1}}
\multiput(6,-4)(0,-4){3}{\circle{0.1}}
\multiput(6,-4)(0,-4){3}{\circle{0.2}}

% directed edges to depth one
\multiput(3,0)(0,-2){6}{\line(-3,-1){6}}
\multiput(3,0)(0,-2){6}{\line(-2,-1){4}}
\multiput(3,-2)(0,-4){3}{\line(-3,-2){3}}
\multiput(3,0)(0,-2){6}{\line(-1,-1){2}}
\multiput(3,0)(0,-2){6}{\line(-1,-2){1}}
\multiput(3,0)(0,-2){6}{\line(1,-2){1}}
\multiput(3,-2)(0,-4){3}{\line(1,-1){2}}
\multiput(3,0)(0,-2){6}{\line(3,-2){3}}

% edges to depth two and three
\multiput(6,-2)(0,-2){6}{\line(1,-1){2}}

% punctured transfer kernel types
\put(-4,-13){\makebox(0,0)[cc]{pTKT}}
\put(-4,-13.5){\makebox(0,0)[cc]{\(\varkappa\sim\)}}
\put(-2.5,-13){\makebox(0,0)[cc]{\(\mathrm{b}.16\)}}
\put(-2.5,-13.5){\makebox(0,0)[cc]{\((004;0)\)}}
\put(-0.5,-13){\makebox(0,0)[cc]{\(\mathrm{b}.3\)}}
\put(-0.5,-13.5){\makebox(0,0)[cc]{\((001;0)\)}}
\put(1.5,-13){\makebox(0,0)[cc]{\(\mathrm{a}.1\)}}
\put(1.5,-13.5){\makebox(0,0)[cc]{\((000;0)\)}}
\put(3.5,-13){\makebox(0,0)[cc]{\(\mathrm{a}.1\)}}
\put(3.5,-13.5){\makebox(0,0)[cc]{\((000;0)\)}}
\put(5.5,-13){\makebox(0,0)[cc]{\(\mathrm{a}.1\)}}
\put(5.5,-13.5){\makebox(0,0)[cc]{\((000;0)\)}}

\end{picture}

}

\end{figure}

%--------------------------------------------------------------------------------

\noindent
Similarly as in Figure
\ref{fig:Tree21a1AscioneA},
the \textit{propagation} in Figure
\ref{fig:Tree31a1AscioneA}
is also completely \textit{regular},
and purely \textit{endo}-genetic
with respect to the coclass tree \(\mathcal{T}^3(M^{(3)}_{1})\).
However,
the \textit{generalized CF-group \(A\) of Ascione}, \(M^{(3)}_{1}\simeq\langle 729,7\rangle\),
is located \textit{on} the shock wave,
and consequently also possesses additional \textit{exo}-genetic \(p\)-descendants
with step size \(s=1\),
namely \(M^{(3)}_{1}-\#1;i\) with \(i\in\lbrace 1,2\rbrace\).
This leads to \textit{exceptional relative identifiers}
for the vertices of the first branch \(\mathcal{B}(M^{(3)}_{1})\),
beginning with \(M^{(3)}_{2}=M^{(3)}_{1}-\#1;3\simeq\langle 2187,113\rangle\).
(The mainline vertex \(M^{(3)}_{i}\) of all other branches with \(i\ge 3\)
has relative identifier \(M^{(3)}_{i-1}-\#1;1\).)

\newpage
%--------------------------------------------------------------------------------

\begin{figure}[ht]
\caption{Coclass-\(4\) tree of type \(\mathrm{a}.1\), rooted in generalized Ascione A, for AQI \((41)\)}
\label{fig:Tree41a1AscioneA}

{\tiny

\setlength{\unitlength}{1cm}
\begin{picture}(14,14.7)(-6,-13.7)

% scale of orders
\put(-5,0.5){\makebox(0,0)[cb]{order}}

\put(-5,0){\line(0,-1){12}}
\multiput(-5.1,0)(0,-2){7}{\line(1,0){0.2}}

\put(-5.2,0){\makebox(0,0)[rc]{\(2\,187\)}}
\put(-4.8,0){\makebox(0,0)[lc]{\(3^7\)}}
\put(-5.2,-2){\makebox(0,0)[rc]{\(6\,561\)}}
\put(-4.8,-2){\makebox(0,0)[lc]{\(3^8\)}}
\put(-5.2,-4){\makebox(0,0)[rc]{\(19\,683\)}}
\put(-4.8,-4){\makebox(0,0)[lc]{\(3^9\)}}
\put(-5.2,-6){\makebox(0,0)[rc]{\(59\,049\)}}
\put(-4.8,-6){\makebox(0,0)[lc]{\(3^{10}\)}}
\put(-5.2,-8){\makebox(0,0)[rc]{\(177\,147\)}}
\put(-4.8,-8){\makebox(0,0)[lc]{\(3^{11}\)}}
\put(-5.2,-10){\makebox(0,0)[rc]{\(531\,441\)}}
\put(-4.8,-10){\makebox(0,0)[lc]{\(3^{12}\)}}
\put(-5.2,-12){\makebox(0,0)[rc]{\(1\,594\,323\)}}
\put(-4.8,-12){\makebox(0,0)[lc]{\(3^{13}\)}}

\put(-5,-12){\vector(0,-1){2}}

% legend
\put(-3.5,0){\circle*{0.2}}
\put(-3,0){\makebox(0,0)[lc]{\(\ldots\) GI and RI action by \(C_2\), cyclic centre}}
\put(-3.5,-0.3){\circle*{0.1}}
\put(-3,-0.3){\makebox(0,0)[lc]{\(\ldots\) non-\(\sigma\), cyclic centre}}

\put(4.5,0.3){\circle{0.2}}
\put(4.5,0.3){\circle{0.1}}
\put(5,0.3){\makebox(0,0)[lc]{\(\ldots\) GI action by \(C_2\)}}
\put(4.5,0){\circle{0.2}}
\put(5,0){\makebox(0,0)[lc]{\(\ldots\) GI and RI action by \(C_2\)}}
\put(4.5,-0.3){\circle{0.1}}
\put(5,-0.3){\makebox(0,0)[lc]{\(\ldots\) non-\(\sigma\)}}

% infinite mainline
\put(3,0.3){\makebox(0,0)[cb]{\(\langle 7\rangle-\#1;1=\)}}
\put(3.2,0){\makebox(0,0)[lb]{\(\langle 111\rangle\)}}
\put(3.1,-1.7){\makebox(0,0)[lb]{\(2;1=\)}}
\put(3.2,-2){\makebox(0,0)[lb]{\(\langle 93\rangle\)}}
\put(3.2,-4){\makebox(0,0)[lb]{\(3\)}}
\multiput(3.2,-6)(0,-2){4}{\makebox(0,0)[lc]{\(1\)}}

\multiput(3,0)(0,-4){4}{\circle{0.2}}
\multiput(3,0)(0,-4){4}{\circle{0.1}}

\multiput(3,-2)(0,-4){3}{\circle{0.2}}

\multiput(3,0)(0,-2){6}{\line(0,-1){2}}

\put(3,-12){\vector(0,-1){2}}
\put(2.8,-12.5){\makebox(0,0)[rc]{main line}}
\put(3.2,-12.5){\makebox(0,0)[lc]{\(\mathcal{T}^4(\langle 2187,111\rangle)\)}}

% other periodic sequences
\put(-2.8,-2){\makebox(0,0)[lc]{\(\langle 94\rangle\)}}
\put(-2.8,-2.3){\makebox(0,0)[lc]{\(=\langle 7\rangle-2;2\)}}
\put(-0.8,-2){\makebox(0,0)[lc]{\(\langle 954\rangle\)}}
\put(-0.8,-2.3){\makebox(0,0)[lc]{\(=\langle 111\rangle-1;3\)}}
\put(1.2,-2){\makebox(0,0)[lc]{\(\langle 96\rangle\)}}
\put(1.2,-2.3){\makebox(0,0)[lc]{\(=2;4\)}}
\put(2.3,-1.7){\makebox(0,0)[lc]{\(2;5=\)}}
\put(2.2,-2){\makebox(0,0)[lc]{\(\langle 97\rangle\)}}
\put(4.2,-2){\makebox(0,0)[lc]{\(\langle 955\rangle\)}}
\put(4.2,-2.3){\makebox(0,0)[lc]{\(=1;4\)}}
\put(6.2,-2){\makebox(0,0)[lc]{\(\langle 95\rangle\)}}
\put(6.4,-2.3){\makebox(0,0)[lc]{\(=2;3\)}}

\put(-2.8,-4){\makebox(0,0)[lc]{\(4\)}}
\put(-0.8,-4){\makebox(0,0)[lc]{\(8\)}}
\put(0.2,-4){\makebox(0,0)[lc]{\(10\)}}
\put(1.2,-4){\makebox(0,0)[lc]{\(6\)}}
\put(2.2,-4){\makebox(0,0)[lc]{\(7\)}}
\put(4.2,-4){\makebox(0,0)[lc]{\(9\)}}
\put(5.2,-4){\makebox(0,0)[lc]{\(11\)}}
\put(6.2,-4){\makebox(0,0)[lc]{\(5\)}}

\multiput(-2.8,-6)(0,-2){4}{\makebox(0,0)[lc]{\(2\)}}
\multiput(-0.8,-6)(0,-2){4}{\makebox(0,0)[lc]{\(6\)}}
\multiput(0.2,-8)(0,-4){2}{\makebox(0,0)[lc]{\(8\)}}
\multiput(1.2,-6)(0,-2){4}{\makebox(0,0)[lc]{\(4\)}}
\multiput(2.2,-6)(0,-2){4}{\makebox(0,0)[lc]{\(5\)}}
\multiput(4.2,-6)(0,-2){4}{\makebox(0,0)[lc]{\(7\)}}
\multiput(5.2,-8)(0,-4){2}{\makebox(0,0)[lc]{\(9\)}}
\multiput(6.2,-6)(0,-2){4}{\makebox(0,0)[lc]{\(3\)}}

% vertices of depth one
\multiput(-3,-2)(0,-4){3}{\circle{0.1}}
\multiput(-3,-4)(0,-4){3}{\circle{0.2}}

\multiput(-1,-2)(0,-4){3}{\circle*{0.1}}
\multiput(-1,-4)(0,-4){3}{\circle*{0.2}}

\multiput(0,-4)(0,-4){3}{\circle*{0.2}}

\multiput(1,-2)(0,-4){3}{\circle{0.1}}
\multiput(1,-4)(0,-4){3}{\circle{0.2}}

\multiput(2,-2)(0,-4){3}{\circle{0.1}}
\multiput(2,-4)(0,-4){3}{\circle{0.2}}

\multiput(4,-2)(0,-4){3}{\circle*{0.1}}
\multiput(4,-4)(0,-4){3}{\circle*{0.2}}

\multiput(5,-4)(0,-4){3}{\circle*{0.2}}

\multiput(6,-2)(0,-4){3}{\circle{0.1}}
\multiput(6,-4)(0,-4){3}{\circle{0.1}}
\multiput(6,-4)(0,-4){3}{\circle{0.2}}

% directed edges to depth one
\multiput(3,0)(0,-2){6}{\line(-3,-1){6}}
\multiput(3,0)(0,-2){6}{\line(-2,-1){4}}
\multiput(3,-2)(0,-4){3}{\line(-3,-2){3}}
\multiput(3,0)(0,-2){6}{\line(-1,-1){2}}
\multiput(3,0)(0,-2){6}{\line(-1,-2){1}}
\multiput(3,0)(0,-2){6}{\line(1,-2){1}}
\multiput(3,-2)(0,-4){3}{\line(1,-1){2}}
\multiput(3,0)(0,-2){6}{\line(3,-2){3}}

% edges to depth two and three
\multiput(6,-2)(0,-2){6}{\line(1,-1){2}}

% punctured transfer kernel types
\put(-4,-13){\makebox(0,0)[cc]{pTKT}}
\put(-4,-13.5){\makebox(0,0)[cc]{\(\varkappa\sim\)}}
\put(-2.5,-13){\makebox(0,0)[cc]{\(\mathrm{b}.16\)}}
\put(-2.5,-13.5){\makebox(0,0)[cc]{\((004;0)\)}}
\put(-0.5,-13){\makebox(0,0)[cc]{\(\mathrm{b}.3\)}}
\put(-0.5,-13.5){\makebox(0,0)[cc]{\((001;0)\)}}
\put(1.5,-13){\makebox(0,0)[cc]{\(\mathrm{a}.1\)}}
\put(1.5,-13.5){\makebox(0,0)[cc]{\((000;0)\)}}
\put(3.5,-13){\makebox(0,0)[cc]{\(\mathrm{a}.1\)}}
\put(3.5,-13.5){\makebox(0,0)[cc]{\((000;0)\)}}
\put(5.5,-13){\makebox(0,0)[cc]{\(\mathrm{a}.1\)}}
\put(5.5,-13.5){\makebox(0,0)[cc]{\((000;0)\)}}

\end{picture}

}

\end{figure}
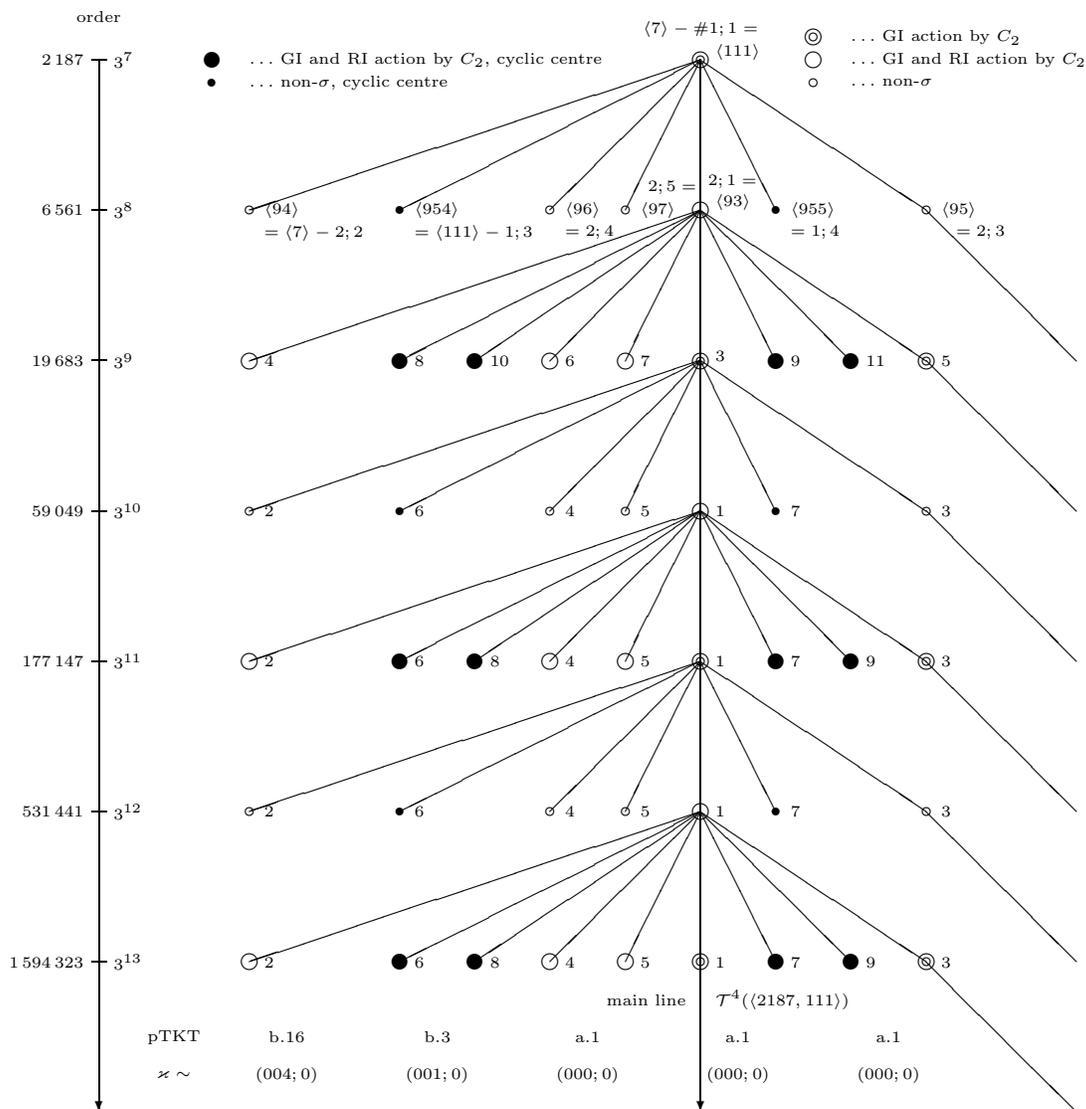

%--------------------------------------------------------------------------------

\noindent
As opposed to the Figures
\ref{fig:Tree21a1AscioneA}
and
\ref{fig:Tree31a1AscioneA},
the entire \textit{first branch} is \textit{exceptional} in Figure
\ref{fig:Tree41a1AscioneA}.
The root \(M^{(4)}_{1}\simeq\langle 2187,111\rangle\),
that is the \textit{generalized CF-group \(A\) of Ascione},
arises by \textit{irregular exo}-genetic propagation
as \(p\)-descendant \(M^{(3)}_{1}-\#1;1\) with step size \(s=1\)
of the periodic root \(M^{(3)}_{1}\simeq\langle 729,7\rangle\) of bifurcations,
according to Formula
\eqref{eqn:Irregular}.
The vertices \(V^{(4)}_{2,i}=M^{(3)}_{1}-\#2;i\) with \(1\le i\le 5\) and \textit{bicyclic} centre
arise by \textit{singular exo}-genetic propagation
as \(p\)-descendants with step size \(s=2\) of \(M^{(3)}_{1}\),
according to Formulas
\eqref{eqn:Singular}
and
\eqref{eqn:SingularOffside}.
The vertices \(V^{(4)}_{2,i+3}=M^{(4)}_{1}-\#1;i\) with \(i\in\lbrace 3,4\rbrace\) and \textit{cyclic} centre,
however, arise by \textit{permanent regular endo}-genetic propagation
as \(p\)-descendants with step size \(s=1\) of the root \(M^{(4)}_{1}\) 
of a \textit{periodic chain} with step size \(s=1\),
according to Formula
\eqref{eqn:PermanentOffside}.
See also the Exhaustion Theorem
\ref{thm:Exhaustion}.

\newpage
%--------------------------------------------------------------------------------

%--------------------------------------------------------------------------------

\noindent
In Figure
\ref{fig:Tree21d10},
the \textit{depth-pruned metabelian skeleton} of the coclass tree with root
\(\langle 729,13\rangle=\langle 81,3\rangle-\#2;10\)
is drawn up to order \(3^{12}\).
Without artificial pruning,
the branches have naturally bounded depth \(4\) (not drawn).
The depth-pruned metabelian branches are periodic with length \(2\).
The infinite main line is of type \(\mathrm{d}.10\)
and consists of \(\sigma\)-groups with GI action by \(C_2\),
every other even with RI action.
The offside vertices are \(\sigma\)-groups for even branches
and non-\(\sigma\) groups for odd branches.
Vertices with positive depth are of type \(\mathrm{D}.10\), \(\mathrm{D}.5\),
\(\mathrm{C}.4\), \(\mathrm{B}.2\) or \(\mathrm{d}.10\).
Metabelian vertices of depth \(\mathrm{dp}\ge 2\)
are exclusively of type \(\mathrm{B}.2\),
for instance,
\(\langle 2187,170\rangle\) has \(105\) descendants with step size \(s=1\),
whereas \(\langle 6561,1682\rangle\) and \(\langle 6561,1685\rangle\) have \(58\) descendants each.

The polarization of all CF-descendants of \(\langle 243,17\rangle=\langle 81,3\rangle-\#1;6\)
in Figure
\ref{fig:Tree21a1AscioneA}
was in the \textit{third} component.
In Figure
\ref{fig:Tree21d10},
we see that, astonishingly,
the polarization of all BCF-descendants of \(\langle 729,13\rangle=\langle 81,3\rangle-\#2;10\)
is in the \textit{first} component,
although both roots are constructed as immediate descendants of the same parent.
However, this discrepancy does not really matter,
since the coclass trees \(\mathcal{T}^2(\langle 243,17\rangle)\) and \(\mathcal{T}^3(\langle 729,13\rangle)\)
are actually completely independent,
in contrast to the following coclass trees in Figure
\ref{fig:Tree31d10}
and
\ref{fig:Tree41d10}.

\newpage
%--------------------------------------------------------------------------------

\noindent
In Figure
\ref{fig:Tree31d10},
the root \(\langle 2187,112\rangle=\langle 729,7\rangle-\#1;2\) is \(p\)-terminal,
and all its immediate descendants on the \textit{first} branch are
step size-\(2\) \(p\)-descendants of \(\langle 729,7\rangle\).
This is the prototype of an application of Formula
\eqref{eqn:d10ExplicitOdd},
\(\mathbb{M}^{(4)}_{1}=M^{(3)}_{1}-\#1;2\), and
\eqref{eqn:D10B2C4D5ExplicitOdd},
\(\mathbb{V}^{(4)}_{2,t-5}=M^{(3)}_{1}-\#2;t\),
both for \(e=e_0=3\), the latter for \(t\in\lbrace 6,7,8,9,10\rbrace\).
The \textit{second} branch is regular with root \(\langle 6561,98\rangle=\langle 729,7\rangle-\#2;6\),
but its relative identifiers are exceptional,
since \(\langle 6561,98\rangle-\#1;i\) with \(1\le i\le 2\) are non-metabelian
and thus remain hidden in the metabelian skeleton.
Beginning with the \textit{third} branch, all \textit{odd} branches are regular
with regular relative identifiers, according to Formula
\eqref{eqn:OddBranchBCF}.
Beginning with the \textit{fourth} branch, all \textit{even} branches are regular
with regular relative identifiers, according to Formula
\eqref{eqn:EvenBranchBCF}.
Less explicitly, Formula
\eqref{eqn:d10ExplicitOdd}
can be expressed by Formula
\eqref{eqn:IrregularBCF},
and Formula
\eqref{eqn:D10B2C4D5ExplicitOdd}
by Formula
\eqref{eqn:SingularBCF}.
The propagation in regular branches is covered by Formula
\eqref{eqn:RegularBCF}.

\newpage
%--------------------------------------------------------------------------------

\begin{figure}[ht]
\caption{Depth-pruned metabelian coclass-\(5\) tree of type \(\mathrm{d}.10\) for AQI \((41)\)}
\label{fig:Tree41d10}

{\tiny

\setlength{\unitlength}{1cm}
\begin{picture}(14,15)(-6,-14)

% scale of orders
\put(-5,0.5){\makebox(0,0)[cb]{order}}

\put(-5,0){\line(0,-1){12}}
\multiput(-5.1,0)(0,-2){7}{\line(1,0){0.2}}

\put(-5.2,0){\makebox(0,0)[rc]{\(6\,561\)}}
\put(-4.8,0){\makebox(0,0)[lc]{\(3^8\)}}
\put(-5.2,-2){\makebox(0,0)[rc]{\(19\,683\)}}
\put(-4.8,-2){\makebox(0,0)[lc]{\(3^9\)}}
\put(-5.2,-4){\makebox(0,0)[rc]{\(59\,049\)}}
\put(-4.8,-4){\makebox(0,0)[lc]{\(3^{10}\)}}
\put(-5.2,-6){\makebox(0,0)[rc]{\(177\,147\)}}
\put(-4.8,-6){\makebox(0,0)[lc]{\(3^{11}\)}}
\put(-5.2,-8){\makebox(0,0)[rc]{\(531\,441\)}}
\put(-4.8,-8){\makebox(0,0)[lc]{\(3^{12}\)}}
\put(-5.2,-10){\makebox(0,0)[rc]{\(1\,594\,323\)}}
\put(-4.8,-10){\makebox(0,0)[lc]{\(3^{13}\)}}
\put(-5.2,-12){\makebox(0,0)[rc]{\(4\,782\,969\)}}
\put(-4.8,-12){\makebox(0,0)[lc]{\(3^{14}\)}}

\put(-5,-12){\vector(0,-1){2}}

% legend
\put(4.5,0.3){\circle{0.2}}
\put(4.5,0.3){\circle*{0.1}}
\put(5,0.3){\makebox(0,0)[lc]{\(\ldots\) GI action by \(C_2\)}}
\put(4.5,0){\circle*{0.2}}
\put(5,0){\makebox(0,0)[lc]{\(\ldots\) GI and RI action by \(C_2\)}}
\put(4.5,-0.3){\circle{0.1}}
\put(5,-0.3){\makebox(0,0)[lc]{\(\ldots\) non-\(\sigma\)}}

% infinite mainline
\put(2.8,0){\makebox(0,0)[rb]{\(\langle 111\rangle-1;2=\)}}
\put(3.2,0){\makebox(0,0)[lb]{\(\langle 953\rangle\)}}
\put(3.2,-2){\makebox(0,0)[lb]{\(\langle 93\rangle-1;2\)}}
\put(3.1,-4){\makebox(0,0)[lb]{\(\langle 93\rangle-2;6\)}}
\put(3.2,-6){\makebox(0,0)[lc]{\(2\)}}
\multiput(3.2,-8)(0,-2){3}{\makebox(0,0)[lc]{\(1\)}}

\multiput(3,0)(0,-4){4}{\circle{0.2}}
\multiput(3,0)(0,-4){4}{\circle*{0.1}}

\multiput(3,-2)(0,-4){3}{\circle*{0.2}}

\multiput(3,0)(0,-2){6}{\line(0,-1){2}}

\put(3,-12){\vector(0,-1){2}}
\put(2.8,-12.5){\makebox(0,0)[rc]{main line}}
\put(3.2,-12.5){\makebox(0,0)[lc]{\(\mathcal{T}^5(\langle 6561,953\rangle)\)}}

% other periodic sequences
\put(-2.8,-2){\makebox(0,0)[lc]{\(\langle 96\rangle-1;2\)}}
\put(-0.8,-2){\makebox(0,0)[lc]{\(\langle 97\rangle-1;2\)}}
\put(1.2,-2){\makebox(0,0)[lc]{\(\langle 94\rangle-1;2\)}}
\put(5.2,-2){\makebox(0,0)[lc]{\(\langle 95\rangle-1;2\)}}

\put(-2.9,-4){\makebox(0,0)[lt]{\(2;10\)}}
\put(-1.9,-4){\makebox(0,0)[lt]{\(2;14\)}}
\put(-0.9,-4){\makebox(0,0)[lt]{\(2;11\)}}
\put(0.1,-4){\makebox(0,0)[lt]{\(2;13\)}}
\put(1.1,-4){\makebox(0,0)[lt]{\(2;7\)}}
\put(2.1,-4){\makebox(0,0)[lt]{\(2;8\)}}
\put(5.1,-4){\makebox(0,0)[lt]{\(2;9\)}}
\put(6.1,-4){\makebox(0,0)[lt]{\(\langle 93\rangle-2;12\)}}

\put(-2.8,-6){\makebox(0,0)[lc]{\(5\)}}
\put(-0.8,-6){\makebox(0,0)[lc]{\(6\)}}
\put(1.2,-6){\makebox(0,0)[lc]{\(3\)}}
\put(5.2,-6){\makebox(0,0)[lc]{\(4\)}}

\put(-2.8,-10){\makebox(0,0)[lc]{\(4\)}}
\put(-0.8,-10){\makebox(0,0)[lc]{\(5\)}}
\put(1.2,-10){\makebox(0,0)[lc]{\(2\)}}
\put(5.2,-10){\makebox(0,0)[lc]{\(3\)}}

\multiput(-2.8,-8)(0,-4){2}{\makebox(0,0)[lc]{\(5\)}}
\multiput(-1.8,-8)(0,-4){2}{\makebox(0,0)[lc]{\(9\)}}
\multiput(-0.8,-8)(0,-4){2}{\makebox(0,0)[lc]{\(6\)}}
\multiput(0.2,-8)(0,-4){2}{\makebox(0,0)[lc]{\(8\)}}
\multiput(1.2,-8)(0,-4){2}{\makebox(0,0)[lc]{\(2\)}}
\multiput(2.2,-8)(0,-4){2}{\makebox(0,0)[lc]{\(3\)}}
\multiput(5.2,-8)(0,-4){2}{\makebox(0,0)[lc]{\(4\)}}
\multiput(6.2,-8)(0,-4){2}{\makebox(0,0)[lc]{\(7\)}}

\multiput(-3,-2)(0,-4){3}{\circle{0.1}}

%\put(-3,-4){\circle{0.2}}
\multiput(-3,-4)(0,-4){3}{\circle*{0.2}}
\multiput(-2,-4)(0,-4){3}{\circle*{0.2}}

\multiput(-1,-2)(0,-4){3}{\circle{0.1}}

%\put(-1,-4){\circle{0.2}}
\multiput(-1,-4)(0,-4){3}{\circle*{0.2}}
\multiput(0,-4)(0,-4){3}{\circle*{0.2}}

\multiput(1,-2)(0,-4){3}{\circle{0.1}}

\multiput(1,-4)(0,-4){3}{\circle*{0.2}}

\multiput(2,-4)(0,-4){3}{\circle*{0.2}}

\multiput(5,-2)(0,-4){3}{\circle{0.1}}

\multiput(5,-4)(0,-4){3}{\circle{0.2}}
%\put(5,-4){\circle{0.1}}
\multiput(5,-4)(0,-4){3}{\circle*{0.1}}
\multiput(6,-4)(0,-4){3}{\circle{0.2}}
\multiput(6,-4)(0,-4){3}{\circle*{0.1}}

% directed edges
\multiput(3,0)(0,-2){6}{\line(-3,-1){6}}
\multiput(3,0)(0,-2){6}{\line(-2,-1){4}}
\multiput(3,0)(0,-2){6}{\line(-1,-1){2}}
\multiput(3,0)(0,-2){6}{\line(1,-1){2}}
\multiput(3,-2)(0,-4){3}{\line(3,-2){3}}

% edges to vertices of depth two, three, four
\multiput(5,-2)(0,-2){5}{\line(1,-1){2}}
\multiput(6,-4)(0,-4){2}{\line(1,-1){2}}

% punctured transfer kernel types
\put(-4,-13){\makebox(0,0)[cc]{pTKT}}
\put(-4,-13.5){\makebox(0,0)[cc]{\(\varkappa\sim\)}}
\put(-2.5,-13){\makebox(0,0)[cc]{\(\mathrm{C}.4\)}}
\put(-2.5,-13.5){\makebox(0,0)[cc]{\((112;2)\)}}
\put(-0.5,-13){\makebox(0,0)[cc]{\(\mathrm{D}.5\)}}
\put(-0.5,-13.5){\makebox(0,0)[cc]{\((113;2)\)}}
\put(1.5,-13){\makebox(0,0)[cc]{\(\mathrm{D}.10\)}}
\put(1.5,-13.5){\makebox(0,0)[cc]{\((114;2)\)}}
\put(3.5,-13){\makebox(0,0)[cc]{\(\mathrm{d}.10\)}}
\put(3.5,-13.5){\makebox(0,0)[cc]{\((110;2)\)}}
\put(5.5,-13){\makebox(0,0)[cc]{\(\mathrm{B}.2\)}}
\put(5.5,-13.5){\makebox(0,0)[cc]{\((111;2)\)}}

\end{picture}

}

\end{figure}

\noindent
In Figure
\ref{fig:Tree41d10},
the root and the initial branches \(1\) and \(2\) of the BCF-tree \(\mathcal{T}^5\) 
are constituted by immediate descendants of CF-groups.
The \textit{root} \(\langle 6561,953\rangle=\langle 2187,111\rangle-\#1;2\)
propagates endo-genetically from the root of the CF-tree \(\mathcal{T}^4\)
with (logarithmic) commutator quotient \((41)\),
according to Formula
\eqref{eqn:IrregularBCF}.
The \textit{first branch} comes from distinct mainline and offside vertices,
\(\langle 6561,i\rangle-\#1;2\) with \(93\le i\le 97\),
on the first branch of \(\mathcal{T}^4\),
according to the same Formula
\eqref{eqn:IrregularBCF}
with type change. 
The entire \textit{second branch} uniformly propagates from the second mainline vertex of \(\mathcal{T}^4\),
according to Formula
\eqref{eqn:SingularBCF}.
Although branch \(3\) is regular, 
its relative identifiers are exceptional,
since \(\langle 6561,93\rangle-\#2;6-\#1;1\) is non-metabelian
and thus does not show up in the metabelian skeleton.
\textit{Regular branches} (third, etc.) are constructed according to Formula
\eqref{eqn:RegularBCF}.

\newpage
%--------------------------------------------------------------------------------

\section{Periodic bifurcations and periodic chains}
\label{s:Periodicity}

\noindent
The statements in this section
exhibit several new kinds of \textit{periodicities} in \(p\)-descendant trees.
The notations are based on both preceding sections, \S\
\ref{s:LawsCF}
on CF-groups, and \S\
\ref{s:LawsBCF}
on BCF-groups.

%--------------------------------------------------------------------------------

Generally, it is convenient to view a coclass tree \(\mathcal{T}^r\)
as union of a finite \textit{pre-period} \(\mathcal{V}\) and an infinite disjoint union
\(\mathcal{T}^r=\mathcal{V}\dot{\cup}\bigl(\dot{\bigcup}_{k\ge 0}\,\mathcal{P}_k\bigr)\)
of copies \(\mathcal{P}_k=\dot{\bigcup}_{i=0}^{\ell-1}\,\mathcal{B}_{p+k\ell+i}\)
of a collection of finitely many branches
\(\dot{\bigcup}_{i=0}^{\ell-1}\,\mathcal{B}_{p+i}\),
the \textit{period} with \textit{length} \(\ell\ge 1\) and \textit{starting subscript} \(p\ge 1\),
such that the branches
\((\forall_{k=0}^{\infty})\,(\forall_{i=0}^{\ell-1})\,\mathcal{B}_{p+k\ell+i}\simeq\mathcal{B}_{p+i}\)
are isomorphic as finite graphs.

In the present article,
all coclass trees are \textit{depth-\(1\) pruned metabelian skeletons}
without pre-period, \(\mathcal{V}=\emptyset\), minimal starting subscript \(p=1\),
and period length \(\ell=2\), that is, we have
\(\mathcal{T}^r=\dot{\bigcup}_{k\ge 0}\,\mathcal{P}_k\)
with
\(\mathcal{P}_k=\mathcal{B}_{1+2k}\dot{\cup}\mathcal{B}_{2+2k}\).

%--------------------------------------------------------------------------------

\begin{definition}
\label{dfn:Branches}
For each integer \(i\ge 1\),
the finite subtree
\(\mathcal{B}(M^{(e)}_{i})=\mathcal{T}^{e}(M^{(e)}_{i})\setminus\mathcal{T}^{e}(M^{(e)}_{i+1})\)
of the depth-pruned CF-coclass tree \(\mathcal{T}^{e}(M^{(e)}_{1})\)
is called \textit{\(i\)-th depth-\(1\) pruned branch},
and the finite subtree
\(\mathcal{B}(\mathbb{M}^{(e+1)}_{i})=\mathcal{T}^{e+1}(\mathbb{M}^{(e+1)}_{i})\setminus\mathcal{T}^{e+1}(\mathbb{M}^{(e+1)}_{i+1})\)
of the depth-pruned metabelian BCF-coclass tree \(\mathcal{T}^{e+1}(\mathbb{M}^{(e+1)}_{1})\)
is called \textit{\(i\)-th depth-\(1\) pruned metabelian branch}.
\end{definition}

%--------------------------------------------------------------------------------

From now on, we omit the phrase \lq\lq depth-\(1\) pruned metabelian\rq\rq.
The precise constitution of the branches in Definition
\ref{dfn:Branches}
by CF-vertices, respectively BCF-vertices, is given experimentally:

\begin{proposition}
\label{prp:OddBranches}
(\textbf{Odd branches})
Let \(e\ge 2\) be a logarithmic integer exponent. \\
For each \textbf{odd} integer \(i\ge 1\), the \(i\)-th CF-branch
\(\mathcal{B}(M^{(e)}_{i})=\lbrace M^{(e)}_{i}\rbrace\cup\lbrace (V^{(e)}_{i+1,j})_{j=2}^7\rbrace\)
consists of the mainline vertex \(M^{(e)}_{i}\) (branch root)
and its immediate step size-\(1\) offside \textbf{descendants}
\begin{equation}
\label{eqn:OddBranchCF}
\begin{aligned}
V^{(e)}_{i+1,2} & \text{ of type } \mathrm{b}.16,\ \varkappa\sim (004;0), \\
V^{(e)}_{i+1,3} & \text{ of type } \mathrm{a}.1,\ \varkappa=(000;0), \textit{ capable twig root}, \\
V^{(e)}_{i+1,4} & \text{ of type } \mathrm{a}.1,\ \varkappa=(000;0), \text{ bicyclic centre } (e-1,1), \\
V^{(e)}_{i+1,5} & \text{ of type } \mathrm{a}.1,\ \varkappa=(000;0), \text{ bicyclic centre } (e-1,1), \\
V^{(e)}_{i+1,6} & \text{ of type } \mathrm{b}.3,\ \varkappa\sim (001;0), \text{ cyclic centre } (e), \\
V^{(e)}_{i+1,7} & \text{ of type } \mathrm{a}.1,\ \varkappa=(000;0), \text{ cyclic centre } (e).
\end{aligned}
\end{equation}
For each \textbf{odd} integer \(i\ge 1\), the \(i\)-th BCF-branch
\(\mathcal{B}(\mathbb{M}^{(e+1)}_{i})=\lbrace\mathbb{M}^{(e+1)}_{i}\rbrace\cup\lbrace (\mathbb{V}^{(e+1)}_{i+1,j})_{j=2}^5\rbrace\)
consists of the mainline vertex \(\mathbb{M}^{(e+1)}_{i}\) (branch root)
and its immediate step size-\(1\) offside \textbf{descendants}
\begin{equation}
\label{eqn:OddBranchBCF}
\begin{aligned}
\mathbb{V}^{(e+1)}_{i+1,2} & \text{ of type } \mathrm{D}.10,\ \varkappa\sim (114;2), \\
\mathbb{V}^{(e+1)}_{i+1,3} & \text{ of type } \mathrm{B}.2,\ \varkappa\sim (111;2), \textit{ capable twig root}, \\
\mathbb{V}^{(e+1)}_{i+1,4} & \text{ of type } \mathrm{C}.4,\ \varkappa\sim (112;2), \\
\mathbb{V}^{(e+1)}_{i+1,5} & \text{ of type } \mathrm{D}.5,\ \varkappa\sim (113;2).
\end{aligned}
\end{equation}
\end{proposition}

%--------------------------------------------------------------------------------

\begin{remark}
\label{rmk:Branches}
1. In order to be able to include the mainline vertex,
we always assume \(V^{(e)}_{i+1,1}=M^{(e)}_{i+1}\) of type \(\mathrm{a}.1\) for CF-groups
and \(\mathbb{V}^{(e+1)}_{i+1,1}=\mathbb{M}^{(e+1)}_{i+1}\) of type \(\mathrm{d}.10\) for BCF-groups
(Propositions
\ref{prp:OddBranches}
--
\ref{prp:EvenBranches},
where the ordering of the offside vertices usually coincides with the ordering in Figures
\ref{fig:Tree31a1AscioneA}
--
\ref{fig:Tree41d10}).

2. The branches for even subscripts \(i\ge 2\) have bigger cardinality
(\(9\) for CF, \(9\) for BCF in Proposition
\ref{prp:EvenBranches})
than those for odd subscripts \(i\ge 1\)
(\(7\) for CF, \(5\) for BCF in Proposition
\ref{prp:OddBranches}).

3. It must be emphasized very clearly
that the (abstract) \textbf{descendants} in Propositions
\ref{prp:OddBranches}
--
\ref{prp:EvenBranches}
are \textbf{not} \(p\)-descendants
in the region \textbf{behind and on the shock wave}.
\end{remark}

\newpage
%--------------------------------------------------------------------------------

\begin{proposition}
\label{prp:EvenBranches}
(\textbf{Even branches})
Let \(e\ge 2\) be a logarithmic integer exponent. \\
For each \textbf{even} integer \(i\ge 2\), the \(i\)-th CF-branch
\(\mathcal{B}(M^{(e)}_{i})=\lbrace M^{(e)}_{i}\rbrace\cup\lbrace (V^{(e)}_{i+1,j})_{j=2}^9\rbrace\)
consists of the mainline vertex \(M^{(e)}_{i}\) (branch root)
and its immediate step size-\(1\) offside \textbf{descendants}
\begin{equation}
\label{eqn:EvenBranchCF}
\begin{aligned}
V^{(e)}_{i+1,2} & \text{ of type } \mathrm{b}.16,\ \varkappa\sim (004;0), \\
V^{(e)}_{i+1,3} & \text{ of type } \mathrm{a}.1,\ \varkappa=(000;0), \textit{ capable twig root}, \\
V^{(e)}_{i+1,4} & \text{ of type } \mathrm{a}.1,\ \varkappa=(000;0), \text{ bicyclic centre } (e-1,1), \\
V^{(e)}_{i+1,5} & \text{ of type } \mathrm{a}.1,\ \varkappa=(000;0), \text{ bicyclic centre } (e-1,1), \\
V^{(e)}_{i+1,6} & \text{ of type } \mathrm{b}.3,\ \varkappa\sim (001;0), \text{ cyclic centre } (e), \\
V^{(e)}_{i+1,7} & \text{ of type } \mathrm{a}.1,\ \varkappa=(000;0), \text{ cyclic centre } (e), \\
V^{(e)}_{i+1,8} & \text{ of type } \mathrm{b}.3,\ \varkappa\sim (001;0), \text{ cyclic centre } (e), \\
V^{(e)}_{i+1,9} & \text{ of type } \mathrm{a}.1,\ \varkappa=(000;0), \text{ cyclic centre } (e).
\end{aligned}
\end{equation}
For each \textbf{even} integer \(i\ge 2\), the \(i\)-th BCF-branch
\(\mathcal{B}(\mathbb{M}^{(e+1)}_{i})=\lbrace\mathbb{M}^{(e+1)}_{i}\rbrace\cup\lbrace (\mathbb{V}^{(e+1)}_{i+1,j})_{j=2}^9\rbrace\)
consists of the mainline vertex \(\mathbb{M}^{(e+1)}_{i}\) (branch root)
and its immediate step size-\(1\) offside \textbf{descendants}
\begin{equation}
\label{eqn:EvenBranchBCF}
\begin{aligned}
\mathbb{V}^{(e+1)}_{i+1,2} & \text{ of type } \mathrm{D}.10,\ \varkappa\sim (114;2), \\
\mathbb{V}^{(e+1)}_{i+1,3} & \text{ of type } \mathrm{D}.10,\ \varkappa\sim (114;2), \\
\mathbb{V}^{(e+1)}_{i+1,4} & \text{ of type } \mathrm{B}.2,\ \varkappa\sim (111;2), \textit{ capable twig root}, \\
\mathbb{V}^{(e+1)}_{i+1,5} & \text{ of type } \mathrm{C}.4,\ \varkappa\sim (112;2), \\
\mathbb{V}^{(e+1)}_{i+1,6} & \text{ of type } \mathrm{D}.5,\ \varkappa\sim (113;2), \\
\mathbb{V}^{(e+1)}_{i+1,7} & \text{ of type } \mathrm{B}.2,\ \varkappa\sim (111;2), \textit{ capable twig root}, \\
\mathbb{V}^{(e+1)}_{i+1,8} & \text{ of type } \mathrm{D}.5,\ \varkappa\sim (113;2), \\
\mathbb{V}^{(e+1)}_{i+1,9} & \text{ of type } \mathrm{C}.4,\ \varkappa\sim (112;2).
\end{aligned}
\end{equation}
\end{proposition}

%--------------------------------------------------------------------------------

%\noindent
We are now in the position to clarify in depth the structure of \textit{periodic bifurcations},
that is, periodic chains with constant step size \(s=2\).
This phenomenon concerns only the links between CF-coclass trees.

\begin{proposition}
\label{prp:Size2Chain}
(\textbf{Periodic size-\(2\) chain})
For each commutator quotient \(C_{3^e}\times C_3\) with \(e\ge 3\),
the mainline of the unique CF-coclass tree \(\mathcal{T}^{e}(M^{(e)}_{1})\),
which starts at the root \(M^{(e)}_{1}\) with \(\mathrm{cl}=3\), \(\mathrm{cc}=e\), \(\mathrm{lo}=3+e\),
contains a unique vertex \(M^{(e)}_{e-2}\) with bifurcation due to nuclear rank \(n=2\),
and with \(\mathrm{cl}=e\), \(\mathrm{cc}=e\), \(\mathrm{lo}=2e\).
The complete periodic size-\(2\) chain is given by \((M^{(e)}_{e-2})_{e\ge 3}\).
\end{proposition}

\begin{proof}
By definition of the coclass tree \(\mathcal{T}^{e}(M^{(e)}_{1})\),
all of its vertices share the common coclass \(\mathrm{cc}=e\).
Each tree arises from a root \(M^{(e)}_{1}\) with class \(\mathrm{cl}=3=1+2\),
whence generally each mainline vertex
\(M^{(e)}_{i}\) is of class \(\mathrm{cl}=i+2\), for \(i\ge 1\).
In particular, the distinguished vertex \(M^{(e)}_{e-2}\) has class \(\mathrm{cc}=(e-2)+2=e\).
Evidence of its elevated nuclear rank \(n=2\) will be provided by the following Theorem
\ref{thm:Bifurcations}.
The logarithmic order is always the sum \(\mathrm{lo}=\mathrm{cl}+\mathrm{cc}\)
of class and coclass.
\end{proof}

%--------------------------------------------------------------------------------

%\noindent
Each of the step sizes,
\(s=1\) and \(s=2\),
of the bifurcations
generates both, exo-genetic and endo-genetic \(p\)-descendants.
Whereas Propositions
\ref{prp:OddBranches}
--
\ref{prp:EvenBranches}
are \textit{not constructive},
the following Theorems
\ref{thm:Bifurcations}
and
\ref{thm:Size1Chains}
can be viewed as \textit{deterministic laws for the construction} of vertices
with the aid of the \(p\)-group generation algorithm
\cite{Nm1977,Ob1990,HEO2005,GNO2006}.
Both theorems are experimental.

\newpage
%--------------------------------------------------------------------------------

\begin{theorem}
\label{thm:Bifurcations}
(\textbf{Structure of bifurcations})
For each integer \(e\ge 3\),
the distinguished CF-mainline vertex \(B:=M^{(e)}_{e-2}\), which possesses bifurcation, gives rise to
\begin{itemize}
\item
\(1\) exo-, and \(8\), respectively \(10\), endo-genetic propagations with step size \(s=1\),
\begin{equation}
\label{eqn:StepSize1}
\begin{aligned}
M^{(e+1)}_{e-2}          &= B-\#1;1 \textbf{ (exo-genetic)}, \\
\mathbb{M}^{(e+1)}_{e-2} &= B-\#1;2, \\
M^{(e)}_{e-1}            &= B-\#1;3, \\
V^{(e)}_{e-1,i-2}        &= B-\#1;i \text{ with } 4\le i\le 9\ (e \text{ odd), respectively } \le 11\ (e \text{ even)};
\end{aligned}
\end{equation}
\item
\(5\) exo-, and \(5\), respectively \(9\), endo-genetic propagations with step size \(s=2\),
\begin{equation}
\label{eqn:StepSize2}
\begin{aligned}
M^{(e+1)}_{e-1}              &= B-\#2;1 \textbf{ (exo-genetic recursion)}, \\
V^{(e+1)}_{e-1,i}            &= B-\#2;i \text{ with } 2\le i\le 5 \textbf{ (exo-genetic)}, \\
\mathbb{M}^{(e+1)}_{e-1}     &= B-\#2;6, \\
\mathbb{V}^{(e+1)}_{e-1,i-5} &= B-\#2;i\text{ with } 7\le i\le 10\ (e \text{ odd), respectively } \le 14\ (e \text{ even)}.
\end{aligned}
\end{equation}
\end{itemize}
In particular,
by endo-genetic propagations,
\(B\) generates
the complete branch \(\mathcal{B}(M^{(e)}_{e-2})\) of CF-groups (with \(s=1\)), and
the complete branch \(\mathcal{B}(\mathbb{M}^{(e+1)}_{e-2})\) of BCF-groups (with \(s=2\)),
with exception of the branch roots, \(M^{(e)}_{e-2}\) and \(\mathbb{M}^{(e+1)}_{e-2}\)
(i.e., only vertices with depth \(\mathrm{dp}=1\)).
\end{theorem}

%--------------------------------------------------------------------------------

%\noindent
In the following Theorem
\ref{thm:Size1Chains},
we abstain from vertices with brushwood type \(\mathrm{B}.2\) of high complexity,
and we restrict our attention to vertices with types \(\mathrm{d}.10\), \(\mathrm{D}.10\), \(\mathrm{C}.4\), and \(\mathrm{D}.5\).
Again, each member of the chains generates both, exo-genetic and endo-genetic \(p\)-descendants.

\begin{theorem}
\label{thm:Size1Chains}
(\textbf{Periodic size-\(1\) chains})
For each integer \(e\ge 3\),
\begin{itemize}
\item
\(M:=M^{(e+1)}_{e-2}=B-\#1;1\) of type \(\mathrm{a}.1\) gives rise to
\begin{equation}
\label{eqn:d10Period}
\begin{aligned}
M^{(e+2)}_{e-2}          &= M-\#1;1 \text{ of type } \mathrm{a}.1 \textbf{ (exo-genetic recursion)}, \\
\mathbb{M}^{(e+2)}_{e-1} &= M-\#1;2 \text{ of type } \mathrm{d}.10, \\
V^{(e+1)}_{e-1,i+3}      &= M-\#1;i \text{ of types } \mathrm{b}.3,\mathrm{a}.1, \text{ with } 3\le i\le 4, \\
V^{(e+1)}_{e-1,i+3}      &= M-\#1;i \text{ of types } \mathrm{b}.3,\mathrm{a}.1, \text{ with } 5\le i\le 6, \text{ only for even } e.
\end{aligned}
\end{equation}
\item
\(V:=V^{(e+1)}_{e-1,2}=B-\#2;2\) of type \(\mathrm{b}.16\) gives rise to
\begin{equation}
\label{eqn:D10Period}
\begin{aligned}
V^{(e+2)}_{e-1,2}          &= V-\#1;1 \text{ of type } \mathrm{b}.16 \textbf{ (exo-genetic recursion)}, \\
\mathbb{V}^{(e+2)}_{e-1,2} &= V-\#1;2 \text{ of type } \mathrm{D}.10, \\
\mathbb{V}^{(e+2)}_{e-1,3} &= V-\#1;3 \text{ of type } \mathrm{D}.10, \text{ only for even } e.
\end{aligned}
\end{equation}
\item
\(V:=V^{(e+1)}_{e-1,4}=B-\#2;4\) of type \(\mathrm{a}.1\) gives rise to
\begin{equation}
\label{eqn:C4Period}
\begin{aligned}
V^{(e+2)}_{e-1,4} &= V-\#1;1 \text{ of type } \mathrm{a}.1 \textbf{ (exo-genetic recursion)}, \\
\mathbb{V}^{(e+2)}_{e-1,i} &= V-\#1;2 \text{ of type } \mathrm{C}.4, \text{ with } i=4\ (e \text{ odd) or } i=5\ (e \text{ even),} \\
\mathbb{V}^{(e+2)}_{e-1,9} &= V-\#1;3 \text{ of type } \mathrm{C}.4, \text{ only for even } e.
\end{aligned}
\end{equation}
\item
\(V:=V^{(e+1)}_{e-1,5}=B-\#2;5\) of type \(\mathrm{a}.1\) gives rise to
\begin{equation}
\label{eqn:D5Period}
\begin{aligned}
V^{(e+2)}_{e-1,5} &= V-\#1;1 \text{ of type } \mathrm{a}.1 \textbf{ (exo-genetic recursion)}, \\
\mathbb{V}^{(e+2)}_{e-1,i} &= V-\#1;2 \text{ of type } \mathrm{D}.5, \text{ with } i=5\ (e \text{ odd) or } i=6\ (e \text{ even),} \\
\mathbb{V}^{(e+2)}_{e-1,8} &= V-\#1;3 \text{ of type } \mathrm{D}.5, \text{ only for even } e.
\end{aligned}
\end{equation}
\end{itemize}
\end{theorem}

\newpage
%--------------------------------------------------------------------------------

\begin{figure}[ht]
\caption{Details of exo- and endo-genetic propagation at a bifurcation (\(e\ge 3\) odd)}
\label{fig:BifurcationOdd}

{\tiny

\setlength{\unitlength}{1cm}
\begin{picture}(14,9)(-5.5,-8)

% scale of orders
\put(-5,0.5){\makebox(0,0)[cb]{lo}}

\put(-5,0){\line(0,-1){6}}
\multiput(-5.1,0)(0,-2){4}{\line(1,0){0.2}}

\put(-5.2,0){\makebox(0,0)[rc]{\(2e\)}}
\put(-5.2,-2){\makebox(0,0)[rc]{\(2e+1\)}}
\put(-5.2,-4){\makebox(0,0)[rc]{\(2e+2\)}}
\put(-5.2,-6){\makebox(0,0)[rc]{\(2e+3\)}}

\put(-5,-6){\vector(0,-1){2}}

% commutator quotients
\put(-4.5,0.5){\makebox(0,0)[lc]{\((3^{e+2},3)\)}}
\put(-4.5,0.0){\makebox(0,0)[lc]{(exo)}}

\put(-2.25,0.7){\line(0,-1){8.5}}

\put(-2.0,0.5){\makebox(0,0)[lc]{\((3^{e+1},3)\)}}
\put(-2.0,0.0){\makebox(0,0)[lc]{(exo)}}

\put(3,0.7){\line(0,-1){8.5}}

\put(6.5,0.5){\makebox(0,0)[lc]{\((3^e,3)\)}}
\put(6.5,0.0){\makebox(0,0)[lc]{(endo)}}

% root
\put(3.5,0.0){\circle{0.2}}
\put(3.8,0.0){\line(0,-1){2}}
\put(3.8,0.0){\line(2,-1){4}}
\put(3.8,0.0){\line(-5,-2){5}}
\put(3.3,0.0){\line(0,-1){4}}
\put(3.3,0.0){\line(1,-1){4}}
\put(3.3,0.0){\line(-3,-4){3}}

% step size s=1, endo-genetic
\put(3.5,-2.0){\circle*{0.2}}
\multiput(4.0,-2.0)(0.5,0){3}{\circle{0.2}}
\multiput(6.5,-2.0)(0.5,0){4}{\circle{0.2}}

% step size s=2, endo-genetic
\multiput(3.5,-4.0)(1.0,0){5}{\circle*{0.2}}
\put(3.5,-4.0){\line(0,-1){2}}
\put(3.5,-4.0){\line(1,-1){2}}

% step size s=1, exo-genetic
\put(-1.5,-2.0){\circle{0.2}}
\put(-1.5,-2.0){\line(-3,-2){3}}
\put(-1.5,-2.0){\line(0,-1){2}}
\put(-1.5,-2.0){\line(1,-4){0.5}}

% step size s=2, exo-genetic
\multiput(0.5,-4.0)(0.5,0){5}{\circle{0.2}}
\multiput(1.0,-4.0)(0.5,0){4}{\line(0,-1){2}}
\multiput(1.0,-4.0)(0.5,0){4}{\line(-5,-2){5}}

% secondary descendants, endo-genetic
\put(-2.0,-4.0){\circle*{0.1}}
\multiput(-1.5,-4.0)(0.5,0){2}{\circle{0.1}}
\multiput(1.0,-6.0)(0.5,0){4}{\circle*{0.1}}
\multiput(3.5,-6.0)(0.25,0){9}{\circle*{0.1}}

% secondary descendants, exo-genetic
\put(-4.5,-4.0){\circle{0.1}}
\multiput(-4.0,-6.0)(0.5,0){4}{\circle{0.1}}

% punctured transfer kernel types
\put(3.5,0.1){\makebox(0,0)[lb]{\(\mathrm{a}.1\) (I)}}

\put(-1.5,-1.8){\makebox(0,0)[cb]{\(\mathrm{a}.1\)}}

\put(3.5,-2.2){\makebox(0,0)[ct]{\(\mathrm{d}.10\)}}
\multiput(4.0,-2.2)(1.0,0){2}{\makebox(0,0)[ct]{\(\mathrm{a}.1\)}}
\multiput(4.0,-2.5)(1.0,0){2}{\line(0,-1){0.7}}
\put(4.5,-2.2){\makebox(0,0)[ct]{\(\mathrm{b}.16\)}}
\put(6.5,-2.2){\makebox(0,0)[ct]{\(\mathrm{a}.1\)}}
\multiput(7.0,-2.2)(1.0,0){2}{\makebox(0,0)[ct]{\(\mathrm{a}.1\)}}
\put(7.5,-2.2){\makebox(0,0)[ct]{\(\mathrm{b}.3\)}}

\put(-4.5,-4.2){\makebox(0,0)[ct]{\(\mathrm{a}.1\)}}

\put(-2.0,-4.2){\makebox(0,0)[ct]{\(\mathrm{d}.10\)}}
\put(-1.5,-4.2){\makebox(0,0)[ct]{\(\mathrm{b}.3\)}}
\put(-1.0,-4.2){\makebox(0,0)[ct]{\(\mathrm{a}.1\)}}

\put(0.4,-3.9){\makebox(0,0)[rb]{(II)}}
\multiput(0.5,-3.8)(1.0,0){3}{\makebox(0,0)[lb]{\(\mathrm{a}.1\)}}
\put(0.9,-3.8){\makebox(0,0)[lb]{\(\mathrm{b}.16\)}}
\put(2.0,-3.8){\makebox(0,0)[lb]{\(\mathrm{a}.1\)}}

\put(3.5,-3.8){\makebox(0,0)[cb]{\(\mathrm{d}.10\)}}
\put(4.5,-4.2){\makebox(0,0)[ct]{\(\mathrm{D}.10\)}}
\put(5.5,-3.8){\makebox(0,0)[cb]{\(\mathrm{B}.2\)}}
\put(5.5,-4.2){\line(0,-1){0.7}}
\put(6.5,-4.2){\makebox(0,0)[ct]{\(\mathrm{C}.4\)}}
\put(7.5,-4.2){\makebox(0,0)[ct]{\(\mathrm{D}.5\)}}

\put(-4.0,-6.2){\makebox(0,0)[ct]{\(\mathrm{b}.16\)}}
\multiput(-3.5,-6.2)(0.5,0){3}{\makebox(0,0)[ct]{\(\mathrm{a}.1\)}}
\put(1.0,-6.2){\makebox(0,0)[ct]{\(\mathrm{D}.10\)}}
\put(1.6,-6.2){\makebox(0,0)[ct]{\(\mathrm{B}.2\)}}
\put(2.2,-6.2){\makebox(0,0)[ct]{\(\mathrm{C}.4\)}}
\put(2.75,-6.2){\makebox(0,0)[ct]{\(\mathrm{D}.5\)}}

% legend
\put(4.5,-6.7){\circle{0.2}}
\put(5,-6.7){\makebox(0,0)[lc]{\(\ldots\) primary CF-group}}
\put(4.5,-7.0){\circle*{0.2}}
\put(5,-7.0){\makebox(0,0)[lc]{\(\ldots\) primary BCF-group}}
\put(4.5,-7.3){\circle{0.1}}
\put(5,-7.3){\makebox(0,0)[lc]{\(\ldots\) secondary CF-group}}
\put(4.5,-7.6){\circle*{0.1}}
\put(5,-7.6){\makebox(0,0)[lc]{\(\ldots\) secondary BCF-group}}

\end{picture}

}

\end{figure}
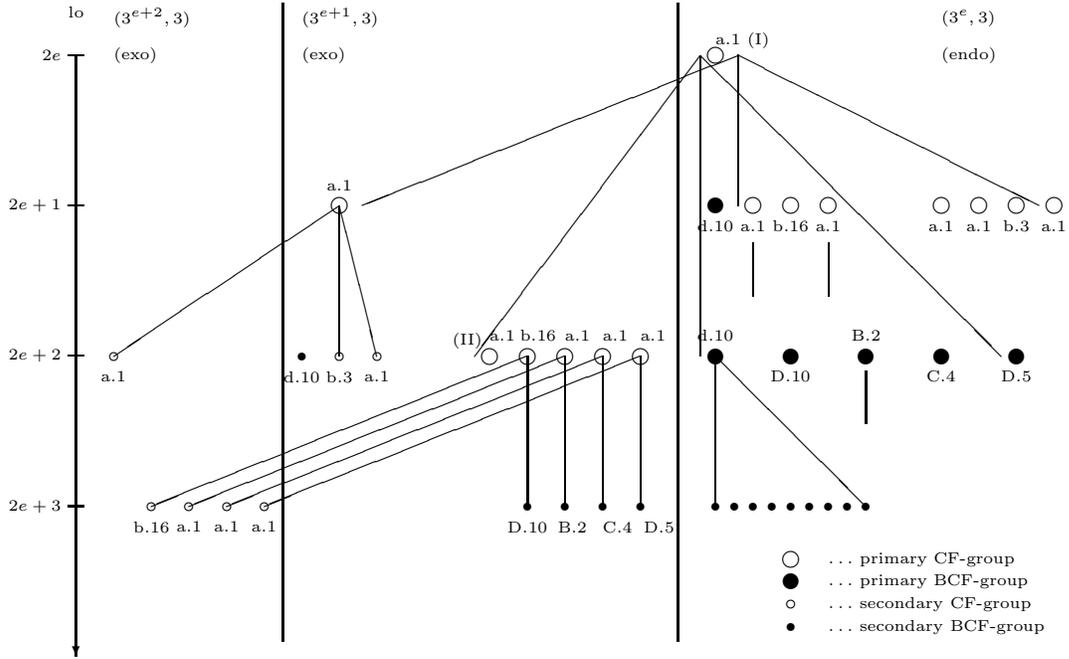

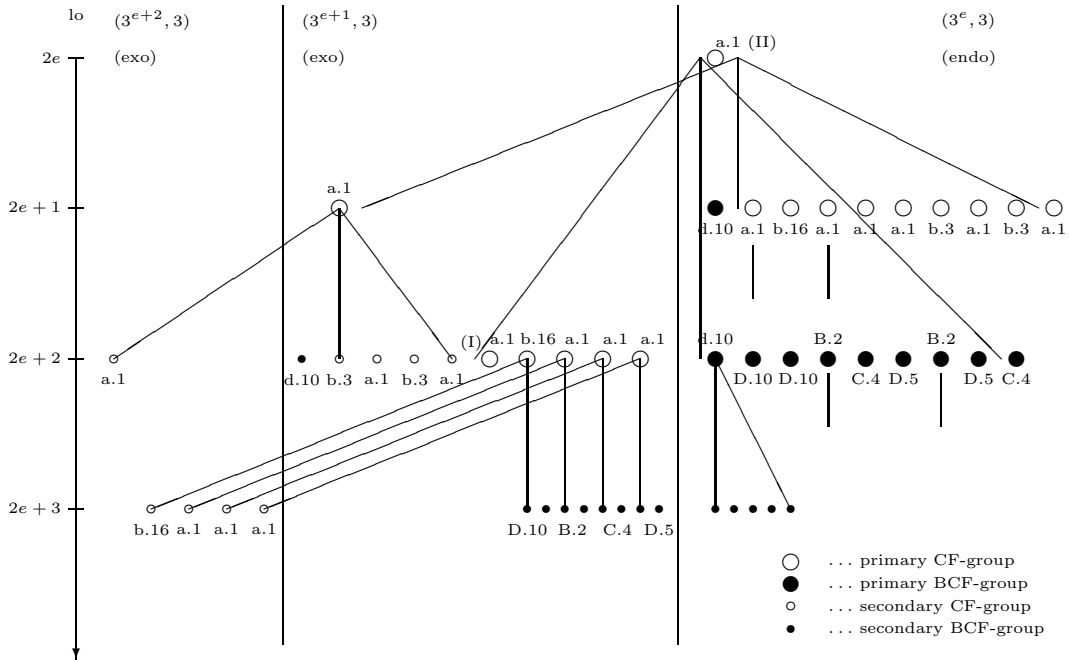
\begin{figure}[ht]
\caption{Details of exo- and endo-genetic propagation at a bifurcation (\(e\ge 4\) even)}
\label{fig:BifurcationEven}

{\tiny

\setlength{\unitlength}{1cm}
\begin{picture}(14,9)(-5.5,-8)

% scale of orders
\put(-5,0.5){\makebox(0,0)[cb]{lo}}

\put(-5,0){\line(0,-1){6}}
\multiput(-5.1,0)(0,-2){4}{\line(1,0){0.2}}

\put(-5.2,0){\makebox(0,0)[rc]{\(2e\)}}
\put(-5.2,-2){\makebox(0,0)[rc]{\(2e+1\)}}
\put(-5.2,-4){\makebox(0,0)[rc]{\(2e+2\)}}
\put(-5.2,-6){\makebox(0,0)[rc]{\(2e+3\)}}

\put(-5,-6){\vector(0,-1){2}}

% commutator quotients
\put(-4.5,0.5){\makebox(0,0)[lc]{\((3^{e+2},3)\)}}
\put(-4.5,0.0){\makebox(0,0)[lc]{(exo)}}

\put(-2.25,0.7){\line(0,-1){8.5}}

\put(-2.0,0.5){\makebox(0,0)[lc]{\((3^{e+1},3)\)}}
\put(-2.0,0.0){\makebox(0,0)[lc]{(exo)}}

\put(3,0.7){\line(0,-1){8.5}}

\put(6.5,0.5){\makebox(0,0)[lc]{\((3^e,3)\)}}
\put(6.5,0.0){\makebox(0,0)[lc]{(endo)}}

% root
\put(3.5,0.0){\circle{0.2}}
\put(3.8,0.0){\line(0,-1){2}}
\put(3.8,0.0){\line(2,-1){4}}
\put(3.8,0.0){\line(-5,-2){5}}
\put(3.3,0.0){\line(0,-1){4}}
\put(3.3,0.0){\line(1,-1){4}}
\put(3.3,0.0){\line(-3,-4){3}}

% step size s=1, endo-genetic
\put(3.5,-2.0){\circle*{0.2}}
\multiput(4.0,-2.0)(0.5,0){9}{\circle{0.2}}

% step size s=2, endo-genetic
\multiput(3.5,-4.0)(0.5,0){9}{\circle*{0.2}}
\put(3.5,-4.0){\line(0,-1){2}}
\put(3.5,-4.0){\line(1,-2){1}}

% step size s=1, exo-genetic
\put(-1.5,-2.0){\circle{0.2}}
\put(-1.5,-2.0){\line(-3,-2){3}}
\put(-1.5,-2.0){\line(0,-1){2}}
\put(-1.5,-2.0){\line(3,-4){1.5}}

% step size s=2, exo-genetic
\multiput(0.5,-4.0)(0.5,0){5}{\circle{0.2}}
\multiput(1.0,-4.0)(0.5,0){4}{\line(0,-1){2}}
\multiput(1.0,-4.0)(0.5,0){4}{\line(-5,-2){5}}

% secondary descendants, endo-genetic
\put(-2.0,-4.0){\circle*{0.1}}
\multiput(-1.5,-4.0)(0.5,0){4}{\circle{0.1}}
\multiput(1.0,-6.0)(0.25,0){8}{\circle*{0.1}}
\multiput(3.5,-6.0)(0.25,0){5}{\circle*{0.1}}

% secondary descendants, exo-genetic
\put(-4.5,-4.0){\circle{0.1}}
\multiput(-4.0,-6.0)(0.5,0){4}{\circle{0.1}}

% punctured transfer kernel types
\put(3.5,0.1){\makebox(0,0)[lb]{\(\mathrm{a}.1\) (II)}}

\put(-1.5,-1.8){\makebox(0,0)[cb]{\(\mathrm{a}.1\)}}

\put(3.5,-2.2){\makebox(0,0)[ct]{\(\mathrm{d}.10\)}}
\multiput(4.0,-2.2)(1.0,0){5}{\makebox(0,0)[ct]{\(\mathrm{a}.1\)}}
\multiput(4.0,-2.5)(1.0,0){2}{\line(0,-1){0.7}}
\put(4.5,-2.2){\makebox(0,0)[ct]{\(\mathrm{b}.16\)}}
\put(5.5,-2.2){\makebox(0,0)[ct]{\(\mathrm{a}.1\)}}
\put(6.5,-2.2){\makebox(0,0)[ct]{\(\mathrm{b}.3\)}}
\put(7.5,-2.2){\makebox(0,0)[ct]{\(\mathrm{b}.3\)}}

\put(-4.5,-4.2){\makebox(0,0)[ct]{\(\mathrm{a}.1\)}}

\put(-2.0,-4.2){\makebox(0,0)[ct]{\(\mathrm{d}.10\)}}
\multiput(-1.5,-4.2)(1.0,0){2}{\makebox(0,0)[ct]{\(\mathrm{b}.3\)}}
\multiput(-1.0,-4.2)(1.0,0){2}{\makebox(0,0)[ct]{\(\mathrm{a}.1\)}}

\put(0.4,-3.9){\makebox(0,0)[rb]{(I)}}
\multiput(0.5,-3.8)(1.0,0){3}{\makebox(0,0)[lb]{\(\mathrm{a}.1\)}}
\put(0.9,-3.8){\makebox(0,0)[lb]{\(\mathrm{b}.16\)}}
\put(2.0,-3.8){\makebox(0,0)[lb]{\(\mathrm{a}.1\)}}

\put(3.5,-3.8){\makebox(0,0)[cb]{\(\mathrm{d}.10\)}}
\multiput(4.0,-4.2)(0.6,0){2}{\makebox(0,0)[ct]{\(\mathrm{D}.10\)}}
\multiput(5.0,-3.8)(1.5,0){2}{\makebox(0,0)[cb]{\(\mathrm{B}.2\)}}
\multiput(5.0,-4.2)(1.5,0){2}{\line(0,-1){0.7}}
\multiput(5.5,-4.2)(2.0,0){2}{\makebox(0,0)[ct]{\(\mathrm{C}.4\)}}
\multiput(6.0,-4.2)(1.0,0){2}{\makebox(0,0)[ct]{\(\mathrm{D}.5\)}}

\put(-4.0,-6.2){\makebox(0,0)[ct]{\(\mathrm{b}.16\)}}
\multiput(-3.5,-6.2)(0.5,0){3}{\makebox(0,0)[ct]{\(\mathrm{a}.1\)}}
\put(1.0,-6.2){\makebox(0,0)[ct]{\(\mathrm{D}.10\)}}
\put(1.6,-6.2){\makebox(0,0)[ct]{\(\mathrm{B}.2\)}}
\put(2.2,-6.2){\makebox(0,0)[ct]{\(\mathrm{C}.4\)}}
\put(2.75,-6.2){\makebox(0,0)[ct]{\(\mathrm{D}.5\)}}

% legend
\put(4.5,-6.7){\circle{0.2}}
\put(5,-6.7){\makebox(0,0)[lc]{\(\ldots\) primary CF-group}}
\put(4.5,-7.0){\circle*{0.2}}
\put(5,-7.0){\makebox(0,0)[lc]{\(\ldots\) primary BCF-group}}
\put(4.5,-7.3){\circle{0.1}}
\put(5,-7.3){\makebox(0,0)[lc]{\(\ldots\) secondary CF-group}}
\put(4.5,-7.6){\circle*{0.1}}
\put(5,-7.6){\makebox(0,0)[lc]{\(\ldots\) secondary BCF-group}}

\end{picture}

}

\end{figure}

Figures
\ref{fig:BifurcationOdd}
and
\ref{fig:BifurcationEven}
provide a graphical illumination of the statements in Theorem
\ref{thm:Bifurcations}
and
\ref{thm:Size1Chains}.
The periodicity of length two is indicated by
roman numerals (I) and (II) at identification points.

%--------------------------------------------------------------------------------

\newpage
%--------------------------------------------------------------------------------

Eventually, we state the main theorem
as the coronation of the present article.

\begin{theorem}
\label{thm:Exhaustion}
\textbf{(Exhaustion theorem)}
Due to an infinite chain of periodic bifurcations,
the \(p\)-descendant tree \(\mathcal{T}_p(R)\) of the metabelian root
\(R=\langle 729,7\rangle\) with abelianization \(R/R^\prime\simeq C_{27}\times C_{3}\)
includes as \textbf{subsets},
for every commutator quotient \(C_{3^e}\times C_{3}\)
with logarithmic exponent \(e\ge 3\),
all depth-pruned coclass trees \(\mathcal{T}^{e}(M^{(e)}_{1})\)
of CF-groups with rank distribution \(\varrho\sim (223;3)\) and
all metabelian skeletons of depth-pruned coclass trees \(\mathcal{T}^{e+1}(\mathbb{M}^{(e+1)}_{1})\)
of BCF-groups with rank distribution \(\varrho\sim (223;3)\).
The former are of type \(\mathrm{a}.1\), \(\varkappa=(000;0)\),
the latter of type \(\mathrm{d}.10\), \(\varkappa=(110;2)\).
The depth-pruning process eliminates all vertices with depth \(\mathrm{dp}\ge 2\).
\end{theorem}

We point out that we cannot speak about subtrees,
because the coclass trees are \textbf{completely disconnected}
as subgraphs of \(p\)-descendants
in the finite region behind the shock wave.
The coclass trees are \textbf{not subtrees} of \(\mathcal{T}_p(R)\)
(the problem are the different edges, not the vertices).

%\newpage
%--------------------------------------------------------------------------------

\begin{proof}
%(Another proof of the Exhaustion Theorem
%\ref{thm:Exhaustion})
Let \(e\ge 3\) be the logarithmic exponent
of an assigned non-elementary bicyclic commutator quotient \(C_{3^e}\times C_{3}\).

First, we show that all vertices of the CF coclass tree \(\mathcal{T}^{e}(M^{(e)}_{1})\)
are \(p\)-descendants of the root \(M^{(3)}_{1}\simeq\langle 729,7\rangle\).
\begin{itemize}
\item
Vertices ahead of the shock wave, with class \(c>e\),
are constructed as regular descendants with endo-genetic propagation by iteration of Formula
\eqref{eqn:Regular}
along the main line, and a single application of Formula
\eqref{eqn:RegularOffside}
for vertices off main line.
\item
For \(e\ge 4\), vertices on the shock wave, with class \(c=e\),
are constructed as singular \(p\)-descendants with exo-genetic propagation by a single application of Formula
\eqref{eqn:Singular},
if they are main line, and Formula
\eqref{eqn:SingularOffside},
if they are offside with bicyclic centre.
If they are offside with cyclic centre,
they are constructed as regular \(p\)-descendants with endo-genetic propagation by Formula
\eqref{eqn:PermanentOffside}.
\item
For \(e\ge 4\), all roots of CF coclass trees, with class \(c=3\),
are constructed as irregular \(p\)-descendants with exo-genetic propagation by iteration of Formula
\eqref{eqn:Irregular}.
\item
In the case \(e\ge 5\), vertices behind the shock wave, with class \(3<c<e\),
are constructed as irregular \(p\)-descendants with exo-genetic propagation by iteration of Formula
\eqref{eqn:Irregular},
if they are main line, and Formula
\eqref{eqn:IrregularOffside},
if they are offside with bicyclic centre.
If they are offside with cyclic centre,
they are constructed as regular \(p\)-descendants with endo-genetic propagation by Formula
\eqref{eqn:PermanentOffside}.
\end{itemize}

Second, we show that all vertices of the BCF coclass tree \(\mathcal{T}^{e+1}(\mathbb{M}^{(e+1)}_{1})\)
are also \(p\)-descendants of the same root \(M^{(3)}_{1}\simeq\langle 729,7\rangle\).
\begin{itemize}
\item
Vertices ahead of the shock wave, with class \(c>e+1\),
are constructed as regular descendants with endo-genetic propagation by iteration of Formula
\eqref{eqn:RegularBCF}
along the main line, followed by a single application of the same Formula
\eqref{eqn:RegularBCF}
for vertices off main line.
\item
Vertices on the shock wave, with class \(c=e+1\),
are constructed as singular \(p\)-descendants with exo-genetic propagation by a single application of Formula
\eqref{eqn:SingularBCF}.
\item
All roots of BCF coclass trees, with class \(c=3\),
are constructed as irregular \(p\)-descendants with exo-genetic propagation by a single application of Formula
\eqref{eqn:IrregularBCF}.
\item
In the case \(e\ge 4\), vertices behind the shock wave, with class \(3<c<e+1\),
are constructed as irregular \(p\)-descendants with exo-genetic propagation by a single application of Formula
\eqref{eqn:IrregularBCF}.
\end{itemize}

\noindent
By the preceding distiction of cases, all claimed metabelian depth-pruned vertices are exhausted.
\end{proof}

\newpage
%--------------------------------------------------------------------------------

%\noindent
The \textbf{Exhaustion Theorem}
can be viewed from another perspective:
instead of recursion formulas,
completely explicit instructions
are given for the construction of
vertices on coclass trees
of CF-groups and BCF-groups.
Assume \(e_0\) is a starting exponent
and \(e\ge e_0\) is a variable exponent.

\medskip
\noindent
For \(e_0\ge 3\) \textbf{odd},

\medskip
CF-groups are constructed as \\
vertices on mainlines of type \(\mathrm{a}.1\),
\begin{equation}
\label{eqn:a1ExplicitOdd}
M^{(e)}_{e_0-2}=M^{(e_0)}_{e_0-2}\lbrack-\#1;1\rbrack^{e-e_0},
\end{equation}
offside vertices of types \(\mathrm{b}.3\) and \(\mathrm{a}.1\) with cyclic centre,
\begin{equation}
\label{eqn:b3a1ExplicitOdd}
V^{(e)}_{e_0-1,j+3}=M^{(e_0)}_{e_0-2}\lbrack-\#1;1\rbrack^{e-e_0}-\#1;j,\ e\ge e_0+1,\ j\in\lbrace\mathbf{3,4}\rbrace,
\end{equation}
and offside vertices of types \(\mathrm{b}.16\), \(\mathrm{a}.1\) twig, and two \(\mathrm{a}.1\) with bicyclic centre,
\begin{equation}
\label{eqn:b16a1ExplicitOdd}
V^{(e)}_{e_0-1,t}=M^{(e_0)}_{e_0-2}-\#2;t\lbrack-\#1;1\rbrack^{e-(e_0+1)},\ e\ge e_0+1,\ t\in\lbrace\mathbf{2,3,4,5}\rbrace;
\end{equation}

\medskip
BCF-groups are constructed as  \\
vertices on mainlines of type \(\mathrm{d}.10\),
\begin{equation}
\label{eqn:d10ExplicitOdd}
\mathbb{M}^{(e+1)}_{e_0-2}=M^{(e_0)}_{e_0-2}\lbrack-\#1;1\rbrack^{e-e_0}-\#1;\mathbf{2},
\end{equation}
and offside vertices of types \(\mathrm{D}.10\), \(\mathrm{B}.2\), \(\mathrm{C}.4\), \(\mathrm{D}.5\), 
\begin{equation}
\label{eqn:D10B2C4D5ExplicitOdd}
\mathbb{V}^{(e+1)}_{e_0-1,t}=M^{(e_0)}_{e_0-2}-\#2;t\lbrack-\#1;1\rbrack^{e-(e_0+1)}-\#1;\mathbf{2},\ e\ge e_0+1,\ t\in\lbrace\mathbf{2,3,4,5}\rbrace.
\end{equation}

\medskip
\noindent
For \(e_0\ge 4\) \textbf{even},

\medskip
CF-groups are constructed as \\
vertices on mainlines of type \(\mathrm{a}.1\),
\begin{equation}
\label{eqn:a1ExplicitEven}
M^{(e)}_{e_0-2}=M^{(e_0)}_{e_0-2}\lbrack-\#1;1\rbrack^{e-e_0},
\end{equation}
offside vertices of types \(\mathrm{b}.3\), \(\mathrm{a}.1\), \(\mathrm{b}.3\), and \(\mathrm{a}.1\) with cyclic centre,
\begin{equation}
\label{eqn:b3a1ExplicitEven}
V^{(e)}_{e_0-1,j+3}=M^{(e_0)}_{e_0-2}\lbrack-\#1;1\rbrack^{e-e_0}-\#1;j,\ e\ge e_0+1,\ j\in\lbrace\mathbf{3,4,5,6}\rbrace,
\end{equation}
and offside vertices of types \(\mathrm{b}.16\), \(\mathrm{a}.1\) twig, and two \(\mathrm{a}.1\) with bicyclic centre,
\begin{equation}
\label{eqn:b16a1ExplicitEven}
V^{(e)}_{e_0-1,t}=M^{(e_0)}_{e_0-2}-\#2;t\lbrack-\#1;1\rbrack^{e-(e_0+1)},\ e\ge e_0+1,\ t\in\lbrace\mathbf{2,3,4,5}\rbrace;
\end{equation}

\medskip
BCF-groups are constructed as \\
vertices on mainlines of type \(\mathrm{d}.10\),
\begin{equation}
\label{eqn:d10ExplicitEven}
\mathbb{M}^{(e+1)}_{e_0-2}=M^{(e_0)}_{e_0-2}\lbrack-\#1;1\rbrack^{e-e_0}-\#1;\mathbf{2},
\end{equation}
and offside vertices of types \(\mathrm{D}.10\), \(\mathrm{B}.2\), \(\mathrm{C}.4\), \(\mathrm{D}.5\), 
\begin{equation}
\label{eqn:D10B2C4D5ExplicitEven}
\mathbb{V}^{(e+1)}_{e_0-1,k(j)}=M^{(e_0)}_{e_0-2}-\#2;t\lbrack-\#1;1\rbrack^{e-(e_0+1)}-\#1;j,\ e\ge e_0+1,\ j\in\lbrace\mathbf{2,3}\rbrace,\ t\in\lbrace\mathbf{2,3,4,5}\rbrace,
\end{equation}
where
\begin{equation}
\label{eqn:Subscripts}
k(j)=
\begin{cases}
j & \text{ if } t=2, \\
4 & \text{ if } t=3,\ j=2, \\
7 & \text{ if } t=3,\ j=3, \\
5 & \text{ if } t=4,\ j=2, \\
9 & \text{ if } t=4,\ j=3, \\
6 & \text{ if } t=5,\ j=2, \\
8 & \text{ if } t=5,\ j=3.
\end{cases}
\end{equation}

\newpage
%--------------------------------------------------------------------------------

\section{Extension and unification of excited states}
\label{s:Unification}

%\noindent
The results concerning periodic non-metabelian Schur \(\sigma\)-groups \(G\)
with moderate rank distribution \(\varrho(G)\) and
types \(\mathrm{D}.10\), \(\mathrm{C}.4\), \(\mathrm{D}.5\) in
\cite{Ma2021a,Ma2021c,Ma2021d}
can be restated, extended, and unified
in the terminology and notation of the present article.
Periodic chains of both step sizes \(s\in\lbrace 1,2\rbrace\) must be employed,
bifurcations with step size \(s=2\) for the selection of \textit{excited states} \(n\ge 0\),
and chains with step size \(s=1\) for growing commutator quotients with logarithmic exponents \(e\ge 2\).

%--------------------------------------------------------------------------------

\subsection{Ground state}
\label{ss:Ground}

%\noindent
Pairs of periodic Schur \(\sigma\)-groups for
the \textit{ground state}, \(n=0\), were discovered in
\cite[\S\ 9, Thm. 12, Eqn. (9.1)--(9.3)]{Ma2021a}.
For each of the types \(\mathrm{D}.10\), \(\mathrm{C}.4\), \(\mathrm{D}.5\),
determined by the fixed parameter \(t\in\lbrace 2,4,5\rbrace\),
they were given by the sequence of doublets
\(G=G(e,i)=\langle 3^8,93\rangle-\#2;t\lbrack-\#1;1\rbrack^{e-5}-\#1;i-\#1;1\)
with running parameter \(e\ge 5\) and selector \(i\in\lbrace 2,3\rbrace\).

The constitution by an infinite main trunk
and finite twigs was illuminated more closely in
\cite[\S\ 4, Thm. 3--6, Eqn. (11)--(22)]{Ma2021d}.
For each \(t\in\lbrace 2,4,5\rbrace\), a periodic chain of CF-groups
\(T_e=\langle 3^8,93\rangle-\#2;t\lbrack-\#1;1\rbrack^{e-5}\) with \(e\ge 5\)
forms the \textit{trunk} of type \(\mathrm{b}.16\) for \(t=2\),
and of type \(\mathrm{a}.1\) for \(t\in\lbrace 4,5\rbrace\).
Each of these vertices \(T_e\) gives rise to a finite \textit{double twig of depth two},
consisting of BCF-groups, the metabelianizations
\(M_{e,i}=T_e-\#1;i\simeq G_{e,i}/G_{e,i}^{\prime\prime}\) with \(i\in\lbrace 2,3\rbrace\) in depth one,
and the Schur \(\sigma\)-groups
\(G_{e,i}=M_{e,i}-\#1;1\) in depth two.
The type is \(\mathrm{D}.10\) for \(t=2\),
\(\mathrm{C}.4\) for \(t=4\), and
\(\mathrm{D}.5\) for \(t=5\).

%--------------------------------------------------------------------------------

In both previous papers
\cite{Ma2021a,Ma2021d},
a connection between the ground state and branches of coclass trees is missing.
The completely explicit notation of the present article
admits the following restatement of all facts concerning the \textit{ground state}.

\begin{theorem}
\label{thm:Ground}
The metabelianizations of the ground state of Schur \(\sigma\)-groups with type
\(\mathrm{D}.10\) for \(t=2\),
\(\mathrm{C}.4\) for \(t=4\), and
\(\mathrm{D}.5\) for \(t=5\)
are given by
\begin{equation}
\label{eqn:Ground}
\mathbb{V}^{(e+1)}_{3,k(j)}=V^{(e)}_{3,t}-\#1;j \quad \text{ with } \quad j\in\lbrace 2,3\rbrace, \quad V^{(e)}_{3,t}=\langle 3^8,93\rangle-\#2;t\lbrack-\#1;1\rbrack^{e-5},
\end{equation}
for each \(e\ge 5\). The subscript \(k(j)\) is given by Formula
\eqref{eqn:Subscripts}.
\(V^{(e)}_{3,t}\) belongs to
the \textit{second} branch \(\mathcal{B}(M^{(e)}_{2})\) of the CF-coclass tree \(\mathcal{T}^{e}(M^{(e)}_{1})\),
and \(\mathbb{V}^{(e+1)}_{3,k(j)}\) belongs to
the \textit{second} branch \(\mathcal{B}(\mathbb{M}^{(e+1)}_{2})\) of the BCF-coclass tree \(\mathcal{T}^{e+1}(\mathbb{M}^{(e+1)}_{1})\).
The Schur \(\sigma\)-group \(\mathbb{V}^{(e+1)}_{3,k(j)}-\#1;1\) has soluble length three.
\end{theorem}

\begin{remark}
\label{rmk:Ground}
So the new insight in comparison to
\cite{Ma2021a,Ma2021d}
is that the endo-genetic propagation behind the shock wave establishes
a branchwise mapping \(V^{(e)}_{3,t}\mapsto\bigl(\mathbb{V}^{(e+1)}_{3,k(2)},\mathbb{V}^{(e+1)}_{3,k(3)}\bigr)\)
from the CF-coclass tree \(\mathcal{T}^{e}(M^{(e)}_{1})\)
to the BCF-coclass tree \(\mathcal{T}^{e+1}(\mathbb{M}^{(e+1)}_{1})\),
for each \(e\ge 5\).
\end{remark}

\begin{proof}
With respect to Schur \(\sigma\)-groups as possible descendants,
only distinguished CF-mainline vertices \(M^{(e)}_{e-2}\) with \textit{even} coclass \(e\ge 4\) are relevant.
For the \textit{ground} state,
we need the smallest even bifurcation \(M^{(4)}_{2}\) with \(e=4\) and exo-genetic offside \(p\)-descendants
\(V^{(5)}_{3,t}=M^{(4)}_{2}-\#2;t\)
with types \(\mathrm{b}.16\) for \(t=2\), and \(\mathrm{a}.1\) for \(t\in\lbrace 4,5\rbrace\),
each of them root of a periodic chain with step size \(s=1\), namely
\(V^{(e)}_{3,t}=M^{(4)}_{2}-\#2;t\lbrack-\#1;1\rbrack^{e-5}\) for \(e\ge 5\),
according to Formula
\eqref{eqn:b16a1ExplicitEven}.
These CF-groups give rise to pairs of BCF-groups as endo-genetic \(p\)-descendants
\(\mathbb{V}^{(6)}_{3,k(j)}=V^{(5)}_{3,t}-\#1;j\), and more generally, for \(e\ge 5\),
\(\mathbb{V}^{(e+1)}_{3,k(j)}=V^{(e)}_{3,t}-\#1;j\) with \(j\in\lbrace 2,3\rbrace\),
according to Formula
\eqref{eqn:D10B2C4D5ExplicitEven}.
In the SmallGroups library
\cite{BEO2005},
\(M^{(4)}_{2}\) has the absolute identifier \(\langle 3^8,93\rangle\),
which completes the proof.
\end{proof}

\newpage
%--------------------------------------------------------------------------------

\subsection{First excited state}
\label{ss:FirstExcited}

%\noindent
Pairs of periodic Schur \(\sigma\)-groups for
the \textit{first excited state}, \(n=1\), were discovered in
\cite[\S\ 2, Thm. 2, Eqn. (2)--(4)]{Ma2021c}.
For each of the types \(\mathrm{D}.10\), \(\mathrm{C}.4\), \(\mathrm{D}.5\),
determined by the fixed parameter \(\ell\in\lbrace 2,4,5\rbrace\),
they were given by the sequence of doublets
\(G=G(e,i)\simeq M(e,i)\lbrack-\#1;1\rbrack^2\)
with metabelianization
\(M=M(e,i)=G(e,i)/G(e,i)^{\prime\prime}\simeq W_\ell\lbrack-\#1;1\rbrack^{e-7}-\#1;\ell\),
where \(e\ge 7\), \(i\in\lbrace 2,3\rbrace\), and
\(W_\ell=\langle 3^8,93\rangle\lbrack-\#2;1\rbrack^2-\#2;\ell\).

The constitution by an infinite main trunk
and finite twigs was illuminated more closely in
\cite[\S\ 5, Thm. 8--10, Eqn. (27)--(38)]{Ma2021d}.
For each \(t\in\lbrace 2,4,5\rbrace\), a periodic chain of CF-groups
\(T_e=\langle 3^8,93\rangle\lbrack-\#2;1\rbrack^{2}-\#2;t\lbrack-\#1;1\rbrack^{e-7}\) with \(e\ge 7\)
forms the \textit{trunk} of type \(\mathrm{b}.16\) for \(t=2\),
and of type \(\mathrm{a}.1\) for \(t\in\lbrace 4,5\rbrace\).
Each of these vertices \(T_e\) gives rise to a finite \textit{double twig of depth three},
consisting of BCF-groups, the metabelianizations
\(M_{e,i}=T_e-\#1;i\simeq G_{e,i}/G_{e,i}^{\prime\prime}\) with \(i\in\lbrace 2,3\rbrace\) in depth one,
and the Schur \(\sigma\)-groups
\(G_{e,i}=M_{e,i}\lbrack-\#1;1\rbrack^{2}\) in depth three.
The type is \(\mathrm{D}.10\) for \(t=2\),
\(\mathrm{C}.4\) for \(t=4\), and
\(\mathrm{D}.5\) for \(t=5\).

%--------------------------------------------------------------------------------

As before, in both previous papers
\cite{Ma2021c,Ma2021d},
a connection between the first excited state and branches of coclass trees is missing.
Again, the completely explicit notation of the present article
admits the following restatement of all facts concerning the \textit{first excited state}.

\begin{theorem}
\label{thm:FirstExcited}
The metabelianizations of the first excited state of Schur \(\sigma\)-groups with type
\(\mathrm{D}.10\) for \(t=2\),
\(\mathrm{C}.4\) for \(t=4\), and
\(\mathrm{D}.5\) for \(t=5\)
are given by
\begin{equation}
\label{eqn:FirstExcited}
\mathbb{V}^{(e+1)}_{5,k(j)}=V^{(e)}_{5,t}-\#1;j \quad \text{ with } \quad j\in\lbrace 2,3\rbrace,
\quad V^{(e)}_{5,t}=\langle 3^8,93\rangle\lbrack-\#2;1\rbrack^{2}-\#2;t\lbrack-\#1;1\rbrack^{e-7},
\end{equation}
for each \(e\ge 7\). The subscript \(k(j)\) is given by Formula
\eqref{eqn:Subscripts}.
\(V^{(e)}_{5,t}\) belongs to
the \textit{fourth} branch \(\mathcal{B}(M^{(e)}_{4})\) of the CF-coclass tree \(\mathcal{T}^{e}(M^{(e)}_{1})\),
and \(\mathbb{V}^{(e+1)}_{5,k(j)}\) belongs to
the \textit{fourth} branch \(\mathcal{B}(\mathbb{M}^{(e+1)}_{4})\) of the BCF-coclass tree \(\mathcal{T}^{e+1}(\mathbb{M}^{(e+1)}_{1})\).
The Schur \(\sigma\)-group \(\mathbb{V}^{(e+1)}_{5,k(j)}\lbrack-\#1;1\rbrack^{2}\) has soluble length three.
\end{theorem}

\begin{remark}
\label{rmk:FirstExcited}
Again, the new insight in comparison to
\cite{Ma2021c,Ma2021d}
is that the endo-genetic propagation behind the shock wave establishes
a branchwise mapping \(V^{(e)}_{5,t}\mapsto\bigl(\mathbb{V}^{(e+1)}_{5,k(2)},\mathbb{V}^{(e+1)}_{5,k(3)}\bigr)\)
from the CF-coclass tree \(\mathcal{T}^{e}(M^{(e)}_{1})\)
to the BCF-coclass tree \(\mathcal{T}^{e+1}(\mathbb{M}^{(e+1)}_{1})\),
for each \(e\ge 7\).
\end{remark}

\begin{proof}
For the \textit{first excited} state,
we need the next even bifurcation \(M^{(6)}_{4}\) with \(e=6\) and exo-genetic offside \(p\)-descendants
\(V^{(7)}_{5,t}=M^{(6)}_{4}-\#2;t\)
with types \(\mathrm{b}.16\) for \(t=2\), and \(\mathrm{a}.1\) for \(t\in\lbrace 4,5\rbrace\),
according to Formula
\eqref{eqn:b16a1ExplicitEven}.
These CF-groups give rise to pairs of BCF-groups as endo-genetic \(p\)-descendants
\(\mathbb{V}^{(8)}_{5,k(j)}=V^{(7)}_{5,t}-\#1;j\) with \(j\in\lbrace 2,3\rbrace\),
according to Formula
\eqref{eqn:D10B2C4D5ExplicitEven}.
In order to start within the SmallGroups database
\cite{BEO2005},
we observe that \(M^{(6)}_{4}=M^{(4)}_{2}\lbrack-\#2;1\rbrack^2\),
where \(M^{(4)}_{2}\simeq\langle 3^8,93\rangle\).
\end{proof}

%--------------------------------------------------------------------------------

\subsection{\(n\)th excited state}
\label{ss:nthExcited}

%\noindent
Now we can easily extend the previous results
by generalization to the \(n\)-th excited state for \(n\ge 2\). 
For the sake of completeness, we include \(n=0\) and \(n=1\).

\begin{theorem}
\label{thm:nthExcited}
The metabelianizations of the \(n\)-th excited state of Schur \(\sigma\)-groups with type
\(\mathrm{D}.10\) for \(t=2\),
\(\mathrm{C}.4\) for \(t=4\), and
\(\mathrm{D}.5\) for \(t=5\)
are given by
\begin{equation}
\label{eqn:nthExcited}
\begin{aligned}
\mathbb{V}^{(e+1)}_{3+2n,k(j)} &= V^{(e)}_{3+2n,t}-\#1;j \text{ with } j\in\lbrace 2,3\rbrace, \\
V^{(e)}_{3+2n,t} &= \langle 3^8,93\rangle\lbrack-\#2;1\rbrack^{2n}-\#2;t\lbrack-\#1;1\rbrack^{e-(5+2n)},
\end{aligned}
\end{equation}
for each \(e\ge 5+2n\). The subscript \(k(j)\) is given by Formula
\eqref{eqn:Subscripts}.
\(V^{(e)}_{3+2n,t}\) belongs to
the \(2(n+1)\)-\textit{th} branch \(\mathcal{B}(M^{(e)}_{2(n+1)})\) of the CF-coclass tree \(\mathcal{T}^{e}(M^{(e)}_{1})\),
and \(\mathbb{V}^{(e+1)}_{3+2n,k(j)}\) belongs to
the \(2(n+1)\)-\textit{th} branch \(\mathcal{B}(\mathbb{M}^{(e+1)}_{2(n+1)})\) of the BCF-coclass tree \(\mathcal{T}^{e+1}(\mathbb{M}^{(e+1)}_{1})\).
The Schur \(\sigma\)-group \(\mathbb{V}^{(e+1)}_{3+2n,k(j)}\lbrack-\#1;1\rbrack^{n+1}\) has soluble length three.
\end{theorem}

\begin{proof}
By induction with respect to the excited state \(n\ge 2\),
using Theorem
\ref{thm:Ground}
for \(n=0\) and Theorem
\ref{thm:FirstExcited}
for \(n=1\) as induction hypothesis.
\end{proof}

\newpage
%--------------------------------------------------------------------------------

\section{Parents of Class Two}
\label{s:Class2}

\noindent
In the proofs of Theorem
\ref{thm:ConstructionMain}
and
\ref{thm:ConstructionBCF},
we had to exclude the investigation of parents
\(A=\pi(D)=D/\gamma_c(D)\)
of the roots
\(D=M^{(e)}_{1}\) respectively \(D=\mathbb{M}^{(e+1)}_{1}\)
of coclass trees for \(e\ge 2\).
In Lemma
\ref{lem:Class2},
we construct a periodic chain with step size \(s=1\),
which consists precisely of these parents.
Since the distinction between CF- and BCF-groups begins with class three,
the parents \(\pi(D)\) are neither CF nor BCF
but simply \textit{class two}.

\begin{lemma}
\label{lem:Class2}
\textbf{(Unboundedly extensible \(3\)-groups of class \(2\))} \\
For each logarithmic exponent \(e\ge 2\),
the \textbf{unique} infinitely capable \(3\)-group \(B_e\) of class \(2\)
and type \(\mathrm{a}.1\), \(\varkappa=(000;0)\),
with commutator quotient \(C_{3^e}\times C_{3}\)
is given as the following member of a periodic chain with step size \(s=1\).
It is parent of both, \(M^{(e)}_{1}\) and \(\mathbb{M}^{(e+1)}_{1}\).
\begin{equation}
\label{eqn:Class2}
B_e:=B\lbrack -\#1;1\rbrack^{e-2}, \quad \pi(M^{(e)}_{1})=B_e, \quad \pi(\mathbb{M}^{(e+1)}_{1})=B_e,
\end{equation}
where \(B=B_2\simeq \mathrm{SmallGroup}(81,3)\) denotes the root of the chain.
\end{lemma}

\begin{example}
\label{exm:Class2}
Aside from the root
\(B_2\simeq\langle 81,3\rangle\),
the SmallGroups database
\cite{BEO2005}
also contains
\(B_3\simeq\langle 243,12\rangle\),
\(B_4\simeq\langle 729,61\rangle\),
\(B_5\simeq\langle 2187,315\rangle\), and
\(B_6\simeq\langle 6561,2063\rangle\).
As the elementary analogue, we can view
the extra special group
\(B_1\simeq\langle 27,3\rangle\)
with commutator quotient \(C_{3}\times C_{3}\).
\end{example}

\begin{proof}
With our usual convention \(s_2=\lbrack y,x\rbrack\) for the main commutator
of a finite two-generated \(3\)-group \(G=\langle x,y\rangle\),
a parametrized pc-presentation of all members of the chain is given by
\begin{equation}
\label{eqn:PresClass2}
B_e=\langle x,y\mid x^{3^{e-1}}=w,\ w^3=1,\ y^3=1,\ s_2^3=1\rangle.
\end{equation}
Whereas the nilpotency class of all members is constant \(\mathrm{cl}(B_e)=2\),
the \(p\)-class \(c_p=\mathrm{cl}_p(B_e)=e\) depends on the logarithmic exponent \(e\).
Since the last non-trivial lower exponent-\(p\) central is
\(P_{c_p-1}(B_e)=\langle w\rangle\),
it follows that \(\pi_p(B_e)=B_e/P_{e-1}(B_e)\simeq B_{e-1}\) for \(e\ge 3\).
Actual computation with Magma
\cite{MAGMA2021}
shows that
\(B_e=B_{e-1}-\#1;1\) for \(e\ge 3\),
and thus by induction
\(B_e=B_2(-\#1;1)^{e-2}\) for \(e\ge 2\).

Now we come to the justification of the parent relations.
First observe that Formula
\eqref{eqn:a1MainPres}
degenerates to 
\begin{equation}
\label{eqn:a1RootPres}
M^{(e)}_{1}=\langle x,y \mid\ x^{3^{e-1}}=w,\ w^3=1,\ y^3=1,\ s_{2}^3=s_{3}^3=1,\ s_3=t_3,\ s_{4}=t_{4}=1\ \rangle.
\end{equation}
in the special case of the root with class \(c=3\), for each \(e\ge 2\).
We put \(D=M^{(e)}_{1}\).
For \(e\ge 4\), we are in the irregular region behind the shock wave, and we have \(c=3\), \(c_p=e\),
\(\gamma_{3}(D)=\langle s_3\rangle\) and \(P_{e-1}(D)=\langle w\rangle\),
whence \(A=\pi(D)=D/\gamma_{3}(D)\simeq B_e\), as claimed,
and \(A_p=\pi_p(D)=D/P_{e-1}(D)\simeq M^{(e-1)}_{1}\), as known from Formula
\eqref{eqn:Irregular}.
For \(e=3\), the behavior on the shock wave is singular, i.e.
\(c=c_p=3\), but \(\gamma_{3}(D)=\langle s_3\rangle\)
as opposed to \(P_{2}(D)=\langle s_3,w\rangle\).
Thus \(A=\pi(D)=D/\gamma_{3}(D)\simeq B_3\),
but \(A_p=\pi_p(D)=D/P_{2}(D)\simeq B_2\),
due to bifurcation.
For \(e=2\), the situation is regular (ahead of the shock wave), i.e.
\(c=c_p=3\), \(\gamma_{3}(D)=P_{2}(D)=\langle s_3\rangle\) and
\(A_p=\pi_p(D)=A=\pi(D)=D/\gamma_{3}(D)\simeq B_2\).

On the other hand, note that Formula
\eqref{eqn:d10MainPres}
degenerates to
\begin{equation}
\label{eqn:d10RootPres}
\mathbb{M}^{(e+1)}_{1}=\langle x,y \mid\ x^{3^{e}}=w,\ w^3=1,\ y^3=1,\ s_{2}^3=s_{3}^3=1,\ t_3=s_3w,\ s_{4}=t_{4}=1\rangle.
\end{equation}
in the special case of the root with class \(c=3\), for each \(e\ge 2\).
We put \(D=\mathbb{M}^{(e+1)}_{1}\).
Then we have \(c=3\), \(c_p=e+1\),
\(\gamma_{3}(D)=\langle s_3,t_3\rangle=\langle s_3,w\rangle\) and \(P_{e}(D)=\langle w\rangle\),
whence \(A=\pi(D)=D/\gamma_{3}(D)\simeq B_e\), as claimed,
and \(A_p=\pi_p(D)=D/P_{e}(D)\simeq M^{(e)}_{1}\), as known from Formula
\eqref{eqn:IrregularBCF}.
\end{proof}

\newpage
%--------------------------------------------------------------------------------

\section{Conclusion}
\label{s:Conclusion}

\noindent
In a series of preceding papers
\cite{Ma2021a,Ma2021b,Ma2021c,Ma2021d},
we have developed a new theory of finite \(3\)-groups \(G\)
with bicyclic commutator quotient \(G/G^\prime\simeq C_{3^e}\times C_3\)
having one non-elementary component with logarithmic exponent \(e\ge 2\).
Theoretical foundations were based on two invariants of \(G\)
with respect to its four maximal subgroups \(H_1,\ldots,H_3;H_4\)
(with distinguished \(H_4\)),
the abelian quotient invariants (AQI) \(\alpha(G)=(H_i/H_i^\prime)_{i=1}^4\) and
the punctured transfer kernel type (pTKT) \(\varkappa(G)=(\ker(T_i))_{i=1}^4\),
combined in the \textit{Artin pattern} \(\mathrm{AP}(G)=(\alpha(G),\varkappa(G))\).

The primary motivation for these works
was the application to possible automorphism groups \(\mathrm{Gal}(\mathrm{F}_3^\infty(K)/K)\)
of \(3\)-class field towers over \textit{imaginary} quadratic number fields
\(K=\mathbb{Q}(\sqrt{d})\), \(d<0\),
which must be \textit{Schur \(\sigma\)-groups}
(with balanced presentation and generator inverting (GI) automorphism).
In the justification of newly discovered \textit{periodicities} among such groups,
two strange phenomena attracted our vigilance and attention:
\begin{itemize}
\item
cumbersome difficulties in the construction of groups with small nilpotency class \(\mathrm{cl}(G)\le e\),
\item
unexpected connections and relationships between CF-groups
\cite{AHL1977}
and BCF-groups
\cite{Ne1989}.
\end{itemize}

In the present article,
we abandoned all motivations by algebraic number theory and class field theory,
we removed the focus on Schur groups and even on \(\sigma\)-groups (except in \S\
\ref{s:Unification}),
and we solved the above mentioned two problems completely
for two infinite families of \textit{coclass trees}
\cite{Ma2015a,Ma2016a},
one, \(\mathcal{T}^{e}(M^{(e)}_1)\), consisting of CF-groups and mainline of type \(\mathrm{a}.1\),
the other, \(\mathcal{T}^{e+1}(\mathbb{M}^{(e+1)}_1)\), consisting of BCF-groups and mainline of type \(\mathrm{d}.10\),
and unbounded \(e\ge 3\) in both situations.

The first difficulty is explained by shedding new light on
the \textit{commutator structure} and \textit{power structure}
and their impact on the descending central series, the lower exponent-\(p\) central series,
and the \(p\)-group generation algorithm
\cite{Nm1977,Ob1990,HEO2005,GNO2006}
(also called \textit{extension algorithm} in
\cite{AHL1977}).

The second phenomenon is due to closely related \textit{power-commutator-presentations}
for certain CF-groups and BCF-groups,
the \textit{mainline principle} for the generator of the
last non-trivial lower central \(\gamma_c(G)=\langle s_c\rangle\),
and peculiarities of the last non-trivial lower \(p\)-central \(P_{c_p-1}(G)=\langle w\rangle\).

The marvellous and astonishing statement of Theorems
\ref{thm:ConstructionMain}
and
\ref{thm:ConstructionOffside}
is the constructibility of all vertices \(V^{(e)}_i\), \(i\ge 2\),
on infinitely many CF coclass trees \(\mathcal{T}^e(M^{(e)}_1)\), \(e\ge 3\),
of type \(\mathrm{a}.1\), \(\varkappa=(0,0,0;0)\), with rank distribution \(\varrho\sim (2,2,3;3)\),
as descendants of a single root \(M^{(3)}_1=\langle 729,7\rangle\),
which is the analogue of \textbf{Ascione's CF-group} \(\mathbf{A}\) for the commutator quotient \((27,3)\).
The highlight of this work, completely unexpected up to now,
asserts the constructibility of all vertices \(\mathbb{V}^{(e+1)}_i\), \(i\ge 2\),
on infinitely many BCF coclass trees \(\mathcal{T}^{e+1}(\mathbb{M}^{(e+1)}_1)\), \(e\ge 3\),
of type \(\mathrm{d}.10\), \(\varkappa=(1,1,0;2)\), also with rank distribution \(\varrho\sim (2,2,3;3)\),
as descendants of the same CF-root \(M^{(3)}_1=\langle 729,7\rangle\),
according to Theorem
\ref{thm:ConstructionBCF}.

%\newpage
%--------------------------------------------------------------------------------

\section{Outlook}
\label{s:Outlook}

\noindent
In view of future research,
it should be pointed out that three similar theorems can be proved
for the root \(\dot{M}^{(3)}_1=\langle 729,6\rangle\),
the analogue of \textbf{Ascione's CF-group} \(\mathbf{G}\) for the commutator quotient \((27,3)\),
which gives rise to infinitely many CF coclass trees \(\mathcal{T}^e(\dot{M}^{(e)}_1)\), \(e\ge 3\),
of the same type \(\mathrm{a}.1\), \(\varkappa=(0,0,0;0)\),
and to infinitely many \textbf{pairs of BCF coclass trees}
\(\mathcal{T}^{e+1}(\dot{\mathbb{M}}^{(e+1)}_1)\) and \(\mathcal{T}^{e+1}(\ddot{\mathbb{M}}^{(e+1)}_1)\), \(e\ge 3\),
of type \(\mathrm{e}.14\), \(\varkappa=(1,2,3;0)\),
all three with a distinct rank distribution \(\varrho\sim (2,2,2;3)\).
As opposed to the trees in the present article,
the \textit{polarization} for all these trees coincides with the \textit{puncture} at the fourth component.

Since the main line of a coclass tree,
\(\mathcal{T}^e(M^{(e)}_1)\) respectively \(\mathcal{T}^{e+1}(\mathbb{M}^{(e+1)}_1)\),
gives rise to an infinite projective limit of the same coclass,
\(M^{(e)}_\infty=\lim_i M^{(e)}_i\) respectively \(\mathbb{M}^{(e+1)}_\infty=\lim_i \mathbb{M}^{(e+1)}_i\),
it would be interesting to investigate
whether \(M^{(3)}_\infty\) \lq\lq generates\rq\rq\ all limit groups
\(M^{(e)}_\infty\) with \(e\ge 4\) and \(\mathbb{M}^{(e+1)}_\infty\) with \(e\ge 3\),
in some sense.

\newpage
%--------------------------------------------------------------------------------

%--------------------------------------------------------------------------------

\end{document}